\documentclass{article}

\usepackage{amsthm, amscd, fixmath, mathrsfs, txfonts}
\usepackage[all]{xy}\CompileMatrices\SelectTips{cm}{10}

\DeclareSymbolFont{bbold}{U}{dsrom}{m}{n}
\DeclareSymbolFontAlphabet{\mathbb}{bbold}

\renewcommand{\P}{{\Bbb P}}
\renewcommand{\rho}{\varrho}
\newcommand{\A}{{\cal A}}

\newcommand{\Bbb}{\mathbb}
\newcommand{\frak}{\mathfrak}
\renewcommand{\cal}{\mathscr}
\newcommand{\catqot}{/\hskip-3pt/}
\newcommand{\C}{{\Bbb C}}

\newcommand{\E}{{\cal E}}

\newcommand{\F}{{\cal F}}

\newcommand{\GL}{\mathop{\rm GL}}

\newcommand{\Hom}{\mathop{{\cal H}\mathit om}}
\newcommand{\HOM}{\mathop{\rm Hom}}

\newcommand{\id}{\mathop{\rm id}}

\newcommand{\K}{{\Bbb K}}
\renewcommand{\L}{{\cal L}}

\newcommand{\Oh}{{\cal O}}
\newcommand{\Pe}{{\Bbb P}}

\newcommand{\PGl}{\mathop{\rm PGL}}
\newcommand{\Proj}{\mathop{\rm Proj}}

\newcommand{\Q}{{\Bbb Q}}
\newcommand{\R}{{\Bbb R}}

\newcommand{\SL}{\mathop{\rm SL}}

\newcommand{\Spec}{\mathop{{\cal S}\mathit pec}}

\newcommand{\Sym}{\mathop{{\cal S}\mathit ym}}

\renewcommand{\tilde}{\widetilde}

\newcommand{\Z}{{\Bbb Z}}

\newcommand{\la}{\lambda}
\newcommand{\lra}{\longrightarrow}
\newcommand{\ra}{\rightarrow}

\newcommand{\lma}{\longmapsto}
\newcommand{\p}{\prime}
\newcommand{\q}{\quad}

\renewcommand{\phi}{\varphi}
\newcommand{\rk}{\mathop{\rm rk}}

\renewcommand{\theta}{\vartheta}
\newcommand{\ul}{\underline}
\newcommand{\ol}{\overline}

\theoremstyle{plain}
\newtheorem{Thm}{\sc Theorem}[subsection]
\newtheorem*{FThm*}{\sc Main Theorem}
\newtheorem*{SThm*}{\sc Second Main Theorem}
\newtheorem*{TThm*}{\sc Third Main Theorem}
\newtheorem*{Thm*}{\sc Theorem}
\newtheorem{Cor}[Thm]{\sc Corollary}
\newtheorem*{Cor*}{\sc Corollary}
\newtheorem{Fact}[Thm]{Facts}
\newtheorem{Prop}[Thm]{\sc Proposition}
\newtheorem*{Prop*}{\sc Proposition}
\newtheorem{Lem}[Thm]{\sc Lemma}
\newtheorem*{Lem*}{\sc Lemma}

\theoremstyle{definition}

\theoremstyle{remark}

\newtheorem{Rem}[Thm]{Remark}
\newtheorem{Ex}[Thm]{Example}
\newtheorem*{Ex*}{Example}
\newtheorem*{Rem*}{Remark}

\newcommand{\Gm}{\mathbb{G}_m}

\newcommand{\glv}{{\mathop{\rm GL}(V)}}

\newcommand{\spec}{\mathop{\rm Spec}}

\newcommand{\df}{\textup{df}}

\begin{document}

\pagestyle{myheadings}
\markboth{\rm T.L.\ G\'omez, A.\ Langer, A.H.W.\ Schmitt, I.\ Sols}
{\rm  Principal Bundles in Arbitrary Characteristic}

\title{Moduli Spaces for Principal Bundles in Arbitrary Characteristic}
\author{Adrian Langer, Alexander H.W.\ Schmitt, Ignacio Sols}
\date{}
%%%
\maketitle
%%%
\begin{abstract}
In this article, we solve the problem of constructing moduli spaces of semista\-ble principal bundles
(and singular versions of them) over smooth projective varieties over algebraically closed ground fields of
positive characteristic.
\end{abstract}
%%%
\tableofcontents
%%%
\section{Introduction}
%%%
In this article, we introduce a formalism for dealing with
principal bundles on projective manifolds defined over an
algebraically closed ground field of arbitrary characteristic which
enables us to construct and compactify the moduli space of
Ramanathan-stable principal bundles. As a major application, we
obtain, under some restrictions on the characteristic of the base
field, the solution of the long-standing problem of constructing
the projective moduli space of semistable principal bundles (with
semisimple structure group) on a smooth projective variety.
\par 
In general, we get different compactifications of the moduli space of
Ramanathan-stable principal $G$-bundles for different representations
$G\hookrightarrow \GL (V)$. In the curve case, all of them are equal,
whereas in higher dimensions we use torsion free sheaves to compactify
the moduli space of semistable principal $G$-bundles, so they naturally
become different.  
\par 
The theory of (semi)stable principal
$G$-bundles starts for the structure group $G=\GL_r(\C)$ as the theory
of (semi)stable vector bundles. Based on his development of Geometric
Invariant Theory, David Mumford proposed the notion of a (semi)stable
vector bundle on a Riemann surface \cite{Mu1}. At about the same time,
Narasimhan and Seshadri made the fundamental discovery that stable
vector bundles on the Riemann surface $X$ are precisely those arising
from irreducible unitary representations of the fundamental group
$\pi_1(X)$ \cite{NS}. (Recall that the relationship between vector
bundles and representations of the fundamental groups was first
investigated by A.\ Weil \cite{Weil}.) Finally, Seshadri gave the GIT
construction of the moduli space of stable vector bundles on a Riemann
surface together with its compactification by S-equivalence classes of
semistable vector bundles \cite{Sesh1}. This construction easily
generalizes to ground fields of arbitrary characteristic.  \par Since
its beginnings, the study of stable $G$-bundles has widely developed
and interacted with other fields. The scope of the theory has been
progressively enlarged by eliminating limitations on the ``three
parameters'' of the theory, i.e., the structure group $G$, the base
manifold $X$, and the ground field $k$. First, in the work of Gieseker
\cite{Gie} and Maruyama \cite{Maru}, the theory of stable vector
bundles was enlarged to a theory of semistable torsion free sheaves on
projective manifolds over fields of characteristic zero. Later,
Simpson brought this theory into its final form \cite{Si}. In the work
\cite{Lang} and \cite{La3}, the barriers of extending Simpson's
results to fields of positive characteristic were finally removed. The
arguments given there improve the formalism even in characteristic
zero.  \par At the time when the results of Gieseker and Maruyama were
published, Rama\-nathan had also treated the theory of principal
$G$-bundles on a compact Riemann surface $X$ for an arbitrary
connected reductive structure group $G$. In the paper
\cite{Ramanathan0}, he introduced the notion of (semi)stability for a
principal $G$-bundle ${\cal P}$ on the Riemann surface $X$ and
generalized the results of the paper \cite{NS}, i.e., linked the
theory of semistable principal bundles on $X$ to the study of
representations of the fundamental group in a compact real form $K$ of
$G$. More important to us is the main result of his PhD thesis,
finished at the Tata Institute in 1976. There, Ramanathan provides an
ingenious GIT construction for the moduli space of semistable
principal $G$-bundles on a compact Riemann surface $X$. Due to the
untimely death of the author, this important result appeared in the
posthumous publication \cite{Ramanathan}. At that time, the subject
had become of general interest to mathematicians and physicists.  \par
In the recent papers \cite{GS}, \cite{GS2}, \cite{Schmitt}, and
\cite{Schmittbesser} two independent---although related---methods for
generalizing Ramanathan's theory to the case of higher dimensional
base manifolds defined over the complex numbers were presented. More
precisely, the moduli space of Ramanathan stable bundles was
constructed and compactified with certain ``generalized" principal
bundles, satisfying a Gieseker type semistability condition.  \par It
thus seemed natural to join the forces of the authors to cope
with the problem of bringing these recent developments to base fields
of arbitrary characteristic. In the present paper, we rewrite the
theory of the paper \cite{Schmitt} from scratch. We will see that the
results of that paper are, in fact, true in positive
characteristic. Furthermore, some of the fundamental discoveries of
the papers \cite{GS}, \cite{GS2}, and \cite{Schmittbesser} also remain
valid over any algebraically closed field. In any case, we manage to construct
our moduli spaces as open subschemes of the moduli spaces of
``$\delta$-semistable pseudo $G$-bundles". In a separate publication
\cite{GLSS2}, we will explain how the approach via the adjoint
representation of a connected reductive group $G$ of \cite{Ramanathan}
and \cite{GS2}, or more generally via faithless representations with
kernel in the center of $G$, may be generalized to positive
characteristic.  \par The main change of philosophy which made the
progress possible is the following: Classically, as suggested by the
work of Ramanan and Ramanathan \cite{RamRam}, one studied
semistability of principal bundles by relating it to the semistability
of associated vector bundles. This works well in characteristic zero
but makes the assumption of sufficiently high characteristic of the
base field necessary while working over fields of positive
characteristic.  In the more recent work quoted above, we linked the
semistability of a principal bundle to the semistability of an
associated \bfseries decorated \rm vector bundle. Unfortunately, the theory
of polynomial representations of the general linear group is more
complicated in positive characteristic (see the book \cite{Green}), so
that the set-up of \cite{Schmitt} and \cite{Schmittbesser} cannot be
directly copied. Nevertheless, the basic idea of that work makes
perfect sense over fields of positive characteristic. Thanks to the
results of \cite{Lang} and \cite{La3}, one may adapt the fundamental
arguments from characteristic zero.  \par Let us introduce a piece of
notation, so that we may state our result in a precise form. In this
paper, we will deal with moduli functors of the form
$$
\begin{array}{rrcl}
\ul{\hbox{\rm M}}^{\rm (s)s}\colon & \ul{\hbox{Sch}}_k &\lra & \ul{\hbox{Set}}\\
& S & \lma & \left\{
\begin{array}{l}
\hbox{Isomorphy classes of families}
\\
\hbox{of (semi)stable objects}
\end{array}\right\}.
\end{array}
$$ In each case, we define S-equivalence on the set of isomorphy
classes of semistable objects (e.g., semistable sheaves or principal
$G$-bundles with fixed numerical data) which, restricted to stable
objects, reduces to isomorphy. Assuming we have the moduli functor
and S-equivalence, we introduce the following convenient
terminology: A \it coarse moduli scheme \rm for the functors
$\ul{\hbox{\rm M}}^{\rm (s)s}$ consists of a scheme ${\rm M}^{\rm
ss}$, an open subscheme ${\rm M}^{\rm s}\subseteq {\rm M}^{\rm ss}$,
and natural transformations of functors
$$
\theta^{\rm(s)s}\colon\ \ul{\hbox{\rm M}}^{\rm (s)s}\lra h_{{\rm M}^{\rm (s)s}}
$$
with the following properties:
\begin{enumerate}
\item[\rm 1.] The space ${\rm M}^{\rm (s)s}$ corepresents
              $\ul{\hbox{\rm M}}^{\rm (s)s}$ with respect to $\theta^{\rm(s)s}$.
              It does so uniformly, if ${\rm Char}(k)>0$, and
              universally, if ${\rm Char}(k)=0$.
              (See \cite{HL}, Definition 2.2.1. Observe that ``uniformly" refers to the base change property for
              \bfseries flat \rm morphisms $\phi$ in that definition.)
\item[\rm 2.] The map $\theta^{\rm s}(k)\colon \ul{\hbox{\rm M}}^{\rm s}(k)\lra
              {\rm M}^{\rm s}(k)$ is a bijection between the set of isomorphy classes
              of stable objects and the closed points of ${\rm M}^{\rm s}$.
\item[\rm 3.] The map $\theta^{\rm ss}(k)\colon \ul{\hbox{\rm M}}^{\rm ss}(k)\lra
              {\rm M}^{\rm ss}(k)$ induces a bijection between the set of S-equivalence clas\-ses
              of semistable objects and the closed points of ${\rm M}^{\rm ss}$.
\end{enumerate}
%%%
The difference between positive and zero characteristic in the above
definition comes from our use of Geometric Invariant Theory, as GIT
quotients in positive characteristic are not necessarily universal
categorical. For $G=\glv$, one can in fact show that, in positive
characteristic, the moduli space of stable sheaves universally
corepresents the moduli functor (see \cite{Lang}, Theorem 0.2). This
follows from the fact that stable sheaves are simple. However, even in
characteristic zero, the sheaves corresponding to stable principal
$G$-bundles on a curve are no longer simple (see \cite{Ramanathan0},
Remark 4.1), so this proof fails in general. We now come to the more
detailed presentation of the contents of our work.
%%%
\subsection{Quasi-projective moduli spaces}
\label{subsec:faithful}
%%%
Let $G$ be a connected semisimple group. Fix a faithful representation
$\rho\colon G\lra \GL(V)$ and note that $\rho(G)\subseteq \SL(V)$. In
characteristic zero, a theory for semistable singular principal
$G$-bundles based on such a representation was developed in
\cite{Schmitt} and \cite{Schmittbesser}. However, some characteristic
zero gadgets such as the Reynolds operator, the
instability flag, and normal forms for homogeneous polynomial
representations were used. In this paper, we will rewrite the theory
from the beginning, such that it becomes independent of the
characteristic of the base field and works without decorated sheaves.
\par We will look at pairs $({\cal A},\tau)$ with a torsion free
${\Oh}_X$-module ${\cal A}$ which has rank $\dim_k(V)$ and trivial
determinant and a homomorphism $\tau\colon {\cal
S}ym^\star(\A\otimes V)^G\lra {\Oh}_X$  of ${\Oh}_X$-algebras which is non-trivial in the
sense that the induced section $\sigma\colon X\lra {\cal S}pec({\cal
S}ym^\star(\A\otimes V)^G)$ is not the zero section. Such a pair is
called a \it pseudo $G$-bundle\rm, and if, furthermore,
$\sigma_U(U)\subset {\cal I}som(V\otimes{\Oh}_U,{\A}^\vee_{|U})/G$,
$U:=U_{\cal A}$ being the maximal open subset where $\A$ is locally
free, we speak of a \it singular principal
$G$-bundle\rm.\footnote{Here, we deviate from the original terminology
in \cite{Schmitt} and \cite{Schmittbesser}.} In the case of a
singular principal $G$-bundle $({\cal A},\tau)$, we get a principal
$G$-bundle $\cal P({\cal A},\tau)$ over $U$, defined by means of base
change:
$$
\xymatrix{
{\cal P}({\cal A},\tau) \ar[r] \ar[d] & {\cal I}som(V\otimes{\Oh}_U,{\A}^\vee_{|U}) \ar[d]
\\
U \ar[r]^-{\sigma_{|U}} & {\cal I}som(V\otimes{\Oh}_U,{\A}^\vee_{|U})/G.
}
$$ We now define the notion of semistability for a singular principal
$G$-bundle $({\cal A},\tau)$.  For this, let $\la\colon {\Bbb
G}_m(k)\lra G$ be a one parameter subgroup of $G$. This yields a
parabolic subgroup $Q_G(\la)$ (see (\ref{parabolic}) below) and a
weighted flag $(V_\bullet(\la),\alpha_\bullet(\la))$ in the vector
space $V$ (see Section \ref{SomeGIT}). A \it reduction of
$({\cal A},\tau)$ to $\la$ \rm is a section $\beta\colon U^\p \lra
{\cal P}({\cal A},\tau)_{|U^\p}/Q_G(\la)$ over an open subset
$U^\p\subseteq U$ with ${\rm codim}_X(X\setminus U^\p)\ge 2$.  This
defines a weighted filtration $({\cal
A}_\bullet(\beta),{\alpha}_{\bullet}(\beta))$ of ${\cal A}$.  Here,
${\alpha}_\bullet(\beta)=(\alpha_t,\dots,\alpha_1)$, if
${\alpha}_\bullet(\lambda)=(\alpha_1,\dots,\alpha_t)$, and the
filtration ${\cal A}_{\bullet}(\beta)\colon\ 0\subsetneq {\cal
A}_{1}\subsetneq \cdots\subsetneq {\cal A}_{t}\subsetneq {\cal A}$ is
obtained as follows: The section
$$
\beta'\colon\ U^\p\stackrel{\beta}{\lra} {\cal P}({\cal A},\tau)_{|U^\p}/Q_G(\la)\hookrightarrow
{\cal I}som(V\otimes{\Oh}_{U'},{\A}^\vee_{|U'})/Q_{\GL(V)}(\la)
$$
yields a filtration
$$
0\subsetneq {\cal A}^\p_{1}\subsetneq \cdots \subsetneq {\cal A}^\p_{t}\subsetneq {\cal A}^\vee_{|U^\p}
$$
of ${\cal A}^\vee_{|U^\p}$ by subbundles with $\rk(\A_{i}^\p)=\dim_k(V_{i})$, $i=1,\dots,t$.
This is because $Q_{\GL(V)}(\la)$ is the $\GL(V)$-stabilizer of the flag $V_{\bullet}(\la)$ and, thus, ${\cal I}som(V\otimes{\Oh}_{U^\p},{\cal A}_{|U^\p}^\vee)/Q_{\GL(V)}(\la)\lra U^\p$ is the bundle of flags in the fibers of
$\A_{|U^\p}^\vee$ having the same dimensions as the members of the flag $V_{\bullet}(\la)$.
We define $\A^{\p\p}_{i}:=\ker(\A_{|U^\p}\lra \A_{t+1-i}^{\p\vee})$, $i=1,\dots,t$,
so that we obtain a filtration
$$
0\subsetneq {\cal A}^{\p\p}_{1}\subsetneq \cdots \subsetneq {\cal A}^{\p\p}_{t}\subsetneq {\cal A}_{|U^\p}
$$
of ${\cal A}_{|U^\p}$ by subbundles. Let $\iota\colon U^\p\lra X$ be the inclusion and define
${\cal A}_{i}$ as the saturation of
${\cal A}\cap \iota_\star({\cal A}^{\p\p}_{i})$, $i=1,\dots,t$.
This is the filtration we denote by ${\cal A}_{\bullet}(\beta)$.
It is worth noting that, if $\la^\p=g\cdot\la\cdot g^{-1}$ for some $g\in G$,
then any reduction to $\la$ may also be interpreted as a reduction to $\la^\p$.
Now, we say that a singular principal $G$-bundle $({\cal A},\tau)$ is \it (semi)stable\footnote{If the word
(semi)stable is used together with the symbol ``$(\le)$", then there are two statements: One for ``semistable"
with ``$\le$" and one for ``stable" with ``$<$".}\rm,
if for every non-trivial one parameter subgroup $\la\colon {\Bbb G}_m(k)\lra G$ and every reduction $\beta$ of $({\cal A},\tau)$ to $\la$, we have
$$
M({\cal A},\tau; \beta):= M\bigl({\cal A}_\bullet(\beta),{\alpha}_\bullet(\beta)\bigr) (\succeq) 0,
$$
where, for every weighted filtration $({\cal A}_\bullet,{\alpha}_\bullet)$ of ${\cal A}$, we set
$$
M({\cal A}_\bullet,{\alpha}_\bullet):=
\sum_{i=1}^t\alpha_{i}\cdot\bigl(\rk \A_{i}\cdot P({\cal A})-\rk\A\cdot P({\cal A}_{i})\bigr).
$$
Finally, there is a notion of S-equivalence which will be explained in Section \ref{subsec:SEqPseudo}
and Remark \ref{rem:SEqSingPrinz}.
%%%
We have the implications
$$
\begin{array}{rcl}
\cal P({\cal A},\tau) \hbox{ is Ramanathan-stable} &\Longrightarrow& ({\cal A},\tau) \hbox{ is stable}
\\
&\Longrightarrow& ({\cal A},\tau) \hbox{ is semistable}
\\
&\Longrightarrow& \cal P({\cal A},\tau) \hbox{ is Ramanathan-semistable}.
\end{array}
$$
More precisely, in our language, Ramanathan's notion of (semi)stability becomes
\begin{equation}
\label{slOpe}
L({\cal A},\tau; \beta):=L\bigl(\A_\bullet(\beta),\alpha_\bullet(\beta)\bigr)
:=
\sum_{i=1}^t\alpha_{i}\cdot\bigl(\rk \A_{i}\cdot \deg({\cal A})-\rk\A\cdot \deg({\cal A}_{i})\bigr)(\ge)0
\end{equation}
for every non-trivial one parameter subgroup $\la\colon {\Bbb G}_m(k)\lra G$ and every reduction $\beta$ of $({\cal A},\tau)$
to $\la$. Here, $\deg$ stands for the degree with respect to the chosen polarization.
%%%
\begin{Rem*}
It is easy to check from the definition that the condition of semistability has to be verified only for the indivisible one parameter subgroups that define maximal parabolic subgroups.
\end{Rem*}
%%%
For a fixed Hilbert polynomial $P$, we define the moduli functors
$$
\begin{array}{rrcl}
\ul{\hbox{\rm M}}_{{P}}^{\rm (s)s}(\rho)\colon & \ul{\hbox{Sch}}_k &\lra & \ul{\hbox{Set}}\\
& S & \lma & \left\{
\begin{array}{l}
\hbox{Isomorphy classes of families of}
\\
\hbox{(semi)stable singular principal $G$-bundles}
\\
\hbox{$({\A}, \tau)$, such that $P(\A)=P$}
\end{array}\right\}.
\end{array}
$$
%%%%
\begin{FThm*}
The coarse moduli space for the functors $\ul{\hbox{\rm M}}_{{P}}^{\rm (s)s}(\rho)$ exists as a quasi-projec\-tive
scheme ${\rm M}_{{P}}^{\rm ss}(\rho)$.
\end{FThm*}
%%%
This moduli space contains the moduli space for Ramanathan-stable principal bundles (whose associated vector
bundle has Hilbert polynomial $P$) as an open subscheme. In particular, we have constructed the moduli spaces
for Ramanathan-stable principal bundles in any characteristic.
\par
The moduli space ${\rm M}_{{P}}^{\rm ss}(\rho)$ will be constructed inside a larger {\bfseries projective} moduli space
${\rm M}_{{P}}^{\widetilde{\delta}\hbox{-}\rm ss}(\rho)$ of $\widetilde{\delta}$-semi\-stable pseudo $G$-bundles,
so that it always comes with a natural compactification.
%%%
\begin{Ex*}  For $G={\rm PGL}_r(k)$, principal $G$-bundles correspond to Azumaya algebras,
so that our construction yields in particular a moduli space for
Azumaya algebras and a compactification by ``Azumaya algebra
sheaves". Such moduli spaces have become of interest recently (see
\cite{HS}, especially Proposition 4.1, and \cite{Yoshi}). More
examples can be found in \cite{GLSS2}.
\end{Ex*}
%%%
\begin{Rem*}[Non-emptiness of moduli spaces]
The above theorem is a mere existence statement. The next step is to investigate the geometry of the moduli spaces. If $X$ is a curve of genus $g \ge 2$, then one can use the moduli stack of principal $G$-bundles. It is a smooth algebraic stack of dimension $(g-1)\cdot \dim(G)$ (see \cite{BehrendThesis}, Corollary 8.1.9). Its connected components are labeled by the elements of the fundamental group $\pi_1(G)$ (see \cite{DriSimp}, Proposition 5, and \cite{Holla}, Proposition 3.15). Estimating the dimension of the locus of unstable principal $G$-bundles, one sees that there are stable principal $G$-bundles for any given topological type $\theta\in\pi_1(G)$ (see \cite{Holla}, Proposition 3.25). (The reader may consult \cite{Ramanathan0} for the topological argument over $k=\C$ and \cite{HeinS}, Proposition 4.2.2, for an existence result on stable principal $G$-bundles with a quasi-parabolic structure.) Our main theorem shows that the moduli space ${\rm M}_{G}^{\rm ss}(\theta)$ of Ramanathan-semistable principal $G$-bundles of ``topological type'' $\theta\in\pi_1(G)$ exists. Using the moduli stack, one checks that it is an irreducible normal variety of dimension $(g-1)\cdot \dim(G)$.
\par 
On higher dimensional base varieties, the geometry of the moduli spaces is completely unknown, even if the base field is $\C$. Note that ${\rm M}_{{P}}^{\rm ss}(\rho)$ always contains the moduli space of slope stable principal $G$-bundles (of the respective numerical invariants) as an open subscheme. To prove non-emptiness of ${\rm M}_{{P}}^{\rm ss}(\rho)$, it hence suffices to construct slope stable principal $G$-bundles. A natural approach is to use stable vector bundles and construct from them principal $G$-bundles by extension of the structure group. On a surface over the field $k=\C$, Balaji used this method to prove interesting existence results for slope stable (and thus stable) principal $G$-bundles (\cite{BUhl}, Theorem 7.10). As he also announces (\cite{BUhl}, Remark 7.2), such existence results are likely to hold in large positive characteristic as well. The details will appear in \cite{BalPar2}.
\end{Rem*}
%%%
\subsection{Projectivity and the semistable reduction theorem}
%%%%
The projectivity of the moduli spaces is not built into our new approach.
The remaining question is thus under which assumptions (on the representation or the characteristic of the base field),
the moduli space ${\rm M}_{{P}}^{\rm ss}(\rho)$ is projective.
Since any connected semisimple group is over $k$ isomorphic to one defined
over the integers, one may assume that $G$ is itself defined over the integers. Then, there is also a faithful
representation $\rho\colon G\lra \GL(V)$ which is defined over the integers. Under this assumption, one may develop
an elegant formalism which provides projective moduli spaces in any dimension, provided that the characteristic of the
base field is either zero or greater than a constant which depends on $\rho$. These results will appear in \cite{GLSS2}.
As remarked before, the moduli spaces will also be projective, if $G$ is one of the classical groups and $\rho$ its standard
representation.
\par
Until very recently, the most general result in that direction was
contained in the work of Balaji and Parameswaran \cite{BP} where
the existence of moduli spaces of semista\-ble principal $G$-bundles
on a smooth projective curve was established under the assumptions
that $G$ is semisimple and the characteristic of the base field is
\bfseries sufficiently large\rm. After the first version of this paper
containing an erroneous proof for semistable reduction appeared, Heinloth managed in
\cite{Heinl} to adapt Langton's algorithm \cite{Langton} to the
setting of principal $G$-bundles. His new approach is to work with
the affine Gra\ss mannian (see \cite{Faltings2} for a discussion
of this object in positive characteristic), so it depends heavily
on the variety being a curve. In our approach, we show a semistable
reduction theorem in all dimensions. 
%%%
\begin{Thm*}[Semistable reduction]
\label{SSRT}
Assume that either $\rho\colon G\lra \GL(V)$ is of low separable
index or $G$ is an adjoint group, $\rho$ is the adjoint
representation and it is of low height. Then, given a semistable
singular principal $G$-bundle ${\cal P}_K$ over $X\times {\rm
Spec}(K)$, where $K$ is the quotient field of the complete discrete
valuation ring $R$, there exists a finite extension $R\subseteq
R'$ such that the pullback ${\cal P}_{K'}$ of ${\cal P}_K$ to
$X\times {\rm Spec}(K')$, $K'$ being the fraction field of $R'$,
extends to a semistable singular principal $G$-bundle ${\cal
P}_{R'}$ over $X\times {\rm Spec}(R')$.
\end{Thm*}
%%%
We also recover the following theorem of Heinloth from \cite{Heinl}:
%%%
\begin{Cor*}[Heinloth] \label{Heinl}
The assertions of the above theorem hold if $X$ is a curve and one
of the following conditions is satisfied:
\begin{itemize}
\item ${\rm Char}(k)=0$.
\item The  simple factors of $G$ \footnote{More precisely, we mean the simple factors of the adjoint form
of $G$.} are of type $A$ and $k$ is arbitrary.
\item The  simple factors of $G$ are of type $A,B,C,D$ and ${\rm Char}(k)\neq 2$.
\item The  simple factors of $G$ are of type $A,B,C,D,G$ and ${\rm Char}(k)>10$.
\item The  simple factors of $G$ are of type $A,B,C,D,G,F,E_6$ and ${\rm Char}(k)>22$.
\item The  simple factors of $G$ are of type $A,B,C,D,G,F,E_6,E_7$ and ${\rm Char}(k)>34$.
\item The  simple factors of $G$ are of type $A,B,C,D,G,F,E$ and ${\rm Char}(k)>58$.
\end{itemize}
%%%
Then, given a semistable principal $G$-bundle ${\cal G}_K$ over $X\times {\rm Spec}(K)$ where $K$ is the spectrum of the complete discrete valuation ring $R$, there exists a finite extension $R\subseteq R'$, such that the pullback
${\cal G}_{K'}$ of ${\cal G}_K$ to $X\times {\rm Spec}(K')$, $K'$ being the fraction field of $R'$, extends to
a semistable principal $G$-bundle ${\cal G}_{R'}$ over $X\times {\rm Spec}(R')$.
\end{Cor*}
%%%
\begin{proof}
Note that semistability is preserved under extension of the structure
group via a central isogeny. So to prove the corollary we can restrict
to a simple group of adjoint type.
\par
In the case of classical groups, the statement can be obtained by familiar methods. We will explain the idea when $G=\mathop{\rm PSO}_n(k)$. Then, we have a short exact sequence of groups
$$
\begin{CD}
\{0\} @>>> \Gm @>>> {\mathop{\rm GO}}_n(k) @>>> {\mathop{\rm PSO}}_n(k)@>>>\{1\},
\end{CD}
$$
where $\mathop{\rm GO} _n(k)$ is the group of matrices, such that
$A^tA=cI$, with $c\in \Gm$. Since $H^2(X_{\rm fl}, \Gm)=0$ ($H^2(X_{\rm fl}, \Gm)=H^2(X, {\cal O}_X)$ by Hilbert's theorem 90, \cite{Milne}, Chapter III, Proposition 4.9), every
principal $\mathop{\rm PSO} _n(k)$-bundle reduces to $\mathop{\rm GO} _n(k)$. Giving a principal $\mathop{\rm GO}_n(k)$-bundle is equivalent to fixing a line bundle $L$ and giving a
pair $(E, \varphi)$, where $E$ is a vector bundle of rank $n$ and
$\varphi\colon E\otimes E\lra L$ is a symmetric non-degenerate bilinear
form. By the theory of decorated vector bundles on curves, developed in \cite{Schmitt0} over $\C$ and extended in this work to positive characteristic (see Remark \ref{rem:DecVb} and \cite{SchmittBook}), there is a moduli space for such objects when dropping the condition on non-degeneracy. Here, the stability concept depends on a parameter $\delta$. The fact that one gets the moduli space for principal ${\rm GO}_n(k)$-bundles for large $\delta$ is exactly the same as for ${\rm SO}_n(k)$-bundles. The latter explained is Example \ref{ex:Classical}.
In the case of other groups the theorem follows directly from the theorem and the remark below.
\end{proof}
%%%
Note that the above proof works only, because in the curve case we need to check semistable reduction only for a single representation. The above statements imply projectivity of the moduli spaces.
%%%
\begin{Cor*}
Under the assumption of the above theorem or corollary, the moduli
space ${\rm M}^{\rm ss}_P(\rho)$ is projective. In particular, in
the curve case, the moduli space ${\rm M}^{\rm ss}(G,t)$ of
semistable principal $G$-bundles (as defined by Ramanathan) of
``topological type'' $t\in \pi_1(G)$ exists as a projective scheme
over $k$.
\end{Cor*}
%%%
%%%
\begin{Rem*} 
i) The low height assumption for the adjoint representation amounts to the following restrictions on the characteristic
of the base field:
\begin{itemize}
\item ${\rm Char}(k)>2n$, if $G$ contains a simple factor  of type $A_n$.
\item ${\rm Char}(k)>4n-2$, if $G$ contains a simple factor  of type $B_n$ or $C_n$.
\item ${\rm Char}(k)>4n-6$, if $G$ contains a simple factor  of type $D_n$.
\item ${\rm Char}(k)>10$, if $G$ contains a simple factor  of type $G_2$.
\item ${\rm Char}(k)>22$, if $G$ contains a simple factor of type $F_4$ or $E_6$.
\item ${\rm Char}(k)>34$, if $G$ contains a simple factor of type  $E_7$.
\item ${\rm Char}(k)>58$, if $G$ contains a simple factor of type $E_8$.
\end{itemize}
\par
ii) Heinloth significantly improved the bounds on the characteristic in his theorem in a recent paper \cite{Hein2}.
\par
iii) In a joint project \cite{HeinS}, Heinloth and the third author have applied the techniques of the current paper to construct moduli spaces for parabolic principal $G$-bundles which are projective under the same hypotheses as stated in the above corollary. In the set-up of parabolic bundles, Heinloth's semistable reduction algorithm could be generalized only to structure groups with classical root systems. For exceptional groups, one had to recur to the approach to semistable reduction which we introduce in the current paper. The work \cite{HeinS} contains an application of moduli spaces of parabolic principal bundles to the cohomology of the moduli stack of principal bundles and may serve as a motivation to study the techniques of the present paper. 
\par 
iv) By the fundamental Theorem \ref{Seshadri's thm} of Seshadri's (which we are going to discuss in the appendix), the existence of ${\rm M}_{{P}}^{\rm ss}(\rho)$ as a {\bfseries projective} scheme implies the semistable reduction theorem for semistable singular principal $G$-bundles. Over higher dimensional base varieties, one may also consider the problem of semistable reduction for {\bfseries slope} semistable singular principal $G$-bundles.
If $k=\C$, this variant of the semistable reduction theorem is established in \cite{BUhl}, Theorem 1.1. Generalizations of that theorem to positive characteristic are announced in loc.\ cit., Remark 7.2. Note however that, over base varieties of dimension at least two, slope semistability is {\bf not equivalent} to semistability (which we are using), so that Balaji's semistable reduction theorem has no implications on the projectivity of ${\rm M}_{{P}}^{\rm ss}(\rho)$.
\end{Rem*}
%%%
\subsection{Notation}
%%%
We work over the algebraically closed field $k$ of characteristic $p\ge 0$.
A \it scheme \rm will be a scheme of finite type over $k$.
For a vector bundle $\E$ over a scheme $X$, we set ${\Pe}(\E):={\rm Proj}(\Sym^\star (\E))$, i.e., ${\Pe}(\E)$ is the projective
bundle of hyperplanes in the fibers of $\E$. An open subset $U\subseteq X$ is said to be \it big\rm, if
${\rm codim}_X(X\setminus U)\ge 2$. The degree $\deg(\E)$ and the Hilbert polynomial $P(\E)$ of a torsion
free coherent ${\Oh}_X$-module $\E$ are taken with respect to the fixed polarization ${\Oh}_X(1)$.
We set $[x]_+:=\max\{\,0,x\,\}$, $x\in\R$.
%%%
\section{Preliminaries}
\label{sec:Prels}
In this section, we collect different results which will be needed throughout the construction
of the moduli space for singular principal $G$-bundles for a semisimple group $G$ via a faithful
representation $G\lra \GL(V)$.
%%%
\subsection{GIT}
\label{SomeGIT}
%%%
We recall some notation and results from Geometric Invariant Theory.
Let $G$ be a reductive group over the field $k$ and $\kappa\colon G\lra \GL(W)$ a representation
on the finite dimensional $k$-vector space $W$. This yields the action
\begin{eqnarray*}
\alpha\colon G\times W &\lra& W
\\
(g,w) &\lma& \kappa(g)(w).
\end{eqnarray*}
Recall that a one parameter subgroup is a homomorphism
$$
\la\colon {\Bbb G}_m(k)\lra G.
$$
Such a one parameter subgroup defines a decomposition
$$
W=\bigoplus_{\gamma\in\Z} W^\gamma
$$
with
$$
W^\gamma=\bigl\{\,w\in W\,|\, \kappa(\la(z))(w)=z^\gamma\cdot w,\ \forall z\in {\Bbb G}_m(k)\,\bigr\},\q
\gamma\in \Z.
$$
Let $\gamma_1<\cdots<\gamma_{t+1}$ be the integers with
$W^\gamma\neq\{0\}$ and
${\gamma}_\bullet(\la):=(\gamma_1,\dots,\gamma_{t+1})$. We define
the flag
$$
W_\bullet(\la):\{0\}\subsetneq W_1:=W^{\gamma_1}\subsetneq W_2:=W^{\gamma_1}\oplus W^{\gamma_2}
\subsetneq\cdots\subsetneq W_t:=W^{\gamma_1}\oplus\cdots\oplus W^{\gamma_t}\subsetneq W
$$
and the tuple ${\alpha}_\bullet(\la):=(\alpha_1,\dots,\alpha_{t})$
of positive rational numbers with
$$
\alpha_i:=\frac{\gamma_{i+1}-\gamma_i}{\dim_k(W)},\q i=1,\dots,t.
$$
If $\la$ is a one parameter subgroup of the special linear group
$\SL(W)$, we will refer to $(W_\bullet(\la),{\alpha}_\bullet(\la))$ as the \it weighted flag
of $\la$\rm. For a point $w\in W\setminus\{0\}$, we define
$$
\mu_\kappa(\la,w):=\max\bigl\{\, \gamma_i\,|\, w\hbox{ has a non-trivial component in $W^{\gamma_i}$,
$i=1,\dots,t+1$}\,\bigr\}.
$$
Note that, for $G=\SL(V)$ and $\la_j\colon {\Bbb G}_m(k)\lra G$, $j=1,2$,
%%%%
\begin{equation}
\label{MuGleich}
\mu_\kappa(\la_1,w)=\mu_\kappa(\la_2,w),
\q\hbox{if }
\bigl(V_{\bullet}(\la_1),{\alpha}_\bullet(\la_1)\bigr)=
\bigl(V_{\bullet}(\la_2),{\alpha}_\bullet(\la_2)\bigr).
\end{equation}
%%%%
(See \cite{GIT}, Proposition 2.7, Chapter 2. Note that we take the weighted flags in $V$ and not in $W$.)
\par
Suppose we are given a projective scheme $X$, a $G$-action $\ol{\sigma}\colon G\times X\lra X$, and
a linearization ${\sigma}\colon G\times \L\lra \L$ of this action in the line bundle $\L$.
For a point $x\in X$ and a one parameter subgroup $\la$, we get the point
$$
x_\infty(\la):=\lim_{z\rightarrow\infty}\ol{\sigma}(\la(z),x).
$$
This point is fixed under the action ${\Bbb G}_m(k)\times X\lra X$,
$(z,x)\lma \ol{\sigma}(\la(z),x)$. Therefore, ${\Bbb G}_m(k)$ acts on the fiber $\L\langle x_\infty(\la)\rangle$.
This action is of the form $l\lma z^\gamma\cdot l$, $z\in {\Bbb G}_m(k)$, $l\in \L\langle x_\infty(\la)\rangle$,
and we set
$$
\mu_\sigma(\la,x):=-\gamma.
$$
\par
For a representation $\kappa$ of $G$ as above, we obtain the action
\begin{eqnarray*}
\ol{\sigma}\colon G\times\P(W^\vee) &\lra& {\Pe}(W^\vee)
\\
\bigl(g,[w]\bigr) &\lma& \bigl[\kappa(g)(w)\bigr]
\end{eqnarray*}
together with an induced linearization $\sigma$ in ${\Oh}_{{\Pe}(W^\vee)}(1)$.
One checks that
\begin{equation}
\label{ComputeMu}
\mu_\kappa(\la,w)=\mu_\sigma\bigl(\la,[w]\bigr),\q \forall w\in W\setminus\{0\},\ \la\colon
{\Bbb G}_m(k)\lra G.
\end{equation}
\par
Finally, we recall that a one parameter subgroup $\la \colon
{\Bbb G}_m(k)\lra G$ gives the parabolic subgroup
\begin{equation}
\label{parabolic}
Q_G(\la):=\bigl\{\,g\in G\,|\, \lim_{z\rightarrow \infty} \la(z)\cdot g\cdot \la(z)^{-1}\hbox{ exists in $G$}
\,\bigr\}.
\end{equation}
The unipotent radical of $Q_G(\la)$ is the subgroup
$$
{\cal R}_u\bigl(Q_G(\la)\bigr):=
\bigl\{\,g\in G\,|\, \lim_{z\rightarrow \infty} \la(z)\cdot g\cdot \la(z)^{-1}=e
\,\bigr\}.
$$
\begin{Rem}
\label{rem:Spri}
In the book \cite{Spri}, one defines the parabolic subgroup
$$
P_G(\la):=\bigl\{\,g\in G\,|\, \lim_{z\rightarrow 0} \la(z)\cdot g\cdot \la(z)^{-1}\hbox{ exists in $G$}\,\bigr\},
$$
i.e.,
$$
P_G(\la)=Q_G(-\la).
$$
Therefore, every parabolic subgroup of $Q$ is of the shape $Q_G(\la)$ for an appropriate one parameter
subgroup $\la$ of $G$. We have chosen a different convention, because it is compatible with our GIT notation.
\end{Rem}
%%%
\subsubsection{Actions on homogeneous spaces}
%%%
Let $H$ be a reductive algebraic group, $G$ a closed reductive subgroup, and $X:=H/G$
the associated affine homogeneous space. Then, the following holds true.
%%%
\begin{Prop}
\label{KraftKuttler}
Suppose that $x\in X$ is a point and $\la\colon {\Bbb G}_m(k)\lra H$ is a one parameter subgroup,
such that $x_0:=\lim_{z\ra\infty} \la(z)\cdot x$ exists in $X$.
Then, $x\in {\cal R}_u(Q_{H}(\la))\cdot x_0$.
\end{Prop}
%%%
\begin{proof}
We may assume $x_0=[e]$, so that $\la$ is a one parameter subgroup of $G$. Define
$$
Y:=\bigl\{\,y\in X\,|\,\lim_{z\ra \infty}\la(z)\cdot y=x_0\,\bigr\}.
$$
This set is closed and invariant under the action of ${\cal R}_u(Q_{H}(\la))$.
Note that viewing $X$ as a variety with ${\Bbb G}_m(k)$-action, $x_0$ is the unique point in $Y$ with a closed
${\Bbb G}_m(k)$-orbit, and by the first lemma in Section III of \cite{Luna} (or Lemma 8.3 in \cite{BR},
or 3.1 in \cite{Hesselink}), there
is a ${\Bbb G}_m(k)$-equivariant morphism $f\colon X\lra T_{x_0}(X)$ which maps $x_0$ to $0$
and is \'etale in $x_0$. Obviously, $f$ maps $Y$ to
%%%
\begin{equation}
\label{TangComp}
\bigl\{\, v\in T_{x_0}X\,|\, \lim_{z\ra \infty} \la(z)\cdot v=0\,\bigr\}
={\frak u}_H(\la)/{\frak u}_G(\la)\subset {\frak h}/{\frak g}.
\end{equation}
%%%
Here, ${\frak u}_H(\la)$ and ${\frak u}_G(\la)$ are the Lie algebras of
${\cal R}_u(Q_H(\la))$ and ${\cal R}_u(Q_G(\la))$, respectively, and ${\frak h}$ and ${\frak g}$
are the Lie algebras of $H$ and $G$, respectively.
Note that ${\frak h}$ and ${\frak g}$ receive their $G$-module structures through
the adjoint representation of $G$, and, moreover, by definition,
$$
{\frak u}_H(\la)=\bigl\{\,v\in {\frak h}\,|\,\lim_{z\ra \infty} \la(z)\cdot v=0\,\bigr\}.
$$
This yields the asserted equality in (\ref{TangComp}).
\par
On the other hand, the dimension of ${\frak u}_H(\la)/{\frak
u}_G(\la)$ equals the one of the ${\cal R}_u(Q_H(\la))$-orbit of
$x_0$ at $X$. By \cite{Hesselink}, Theorem 3.4, $f$ maps $Y$
isomorphically onto ${\frak u}_H(\la)/{\frak u}_G(\la)$.
Therefore, since ${\cal R}_u(Q_H(\la))\cdot x_0\subseteq Y$, the
subset $Y$ must agree with the closed orbit ${\cal
R}_u(Q_H(\la))\cdot x_0$, and we are done.
\end{proof}
%%%
The proof of the above result was communicated to us by Kraft and Kuttler (cf.\ \cite{Schmittbesser}).
Its purpose is to characterize one parameter subgroups of $G$ among the one parameter subgroups of
$\SL(V)$, given a faithful representation $\rho\colon G\lra\GL(V)$.
%%%%%%
\subsubsection{Some specific quotient problems}
\label{subsub:Spec}
%%%
A key of understanding classification problems for vector bundles
together with a section in an associated vector bundle is to study the
representation defining the associated vector bundle. In our case, we have to study a certain GIT problem which we will now
describe.
\par
As before, we fix a representation $\rho\colon G\lra \GL(V)$ on the finite
dimensional $k$-vector space $V$.
We look at the representation
\begin{eqnarray*}
R\colon {\GL}_r(k)\times G&\lra&\GL(k^r\otimes V)
\\
(g,g^\p) &\lma& \Bigl(\, w\otimes v\in k^r\otimes V\lma (g\cdot w)\otimes \rho(g^\p)(v)\,\Bigr).
\end{eqnarray*}
%%%
The representation $R$ provides an action of $G\times \GL_r(k)$ on
$$
(V\otimes k^r)^\vee={\rm Hom}(k^r,V^\vee)\q\hbox{and}\q {\Pe}\bigl({\rm Hom}(k^r,V^\vee)^\vee\bigr)
$$
and induces a $\GL_r(k)$-action on the categorical quotients
$$
{\Bbb H}:={\rm Hom}(k^r,V^\vee)\catqot G\q\hbox{and}\q
\ol{\Bbb H}:={\Pe}\bigl({\rm Hom}(k^r,V^\vee)^\vee\bigr)\catqot G=({\Bbb H}\setminus\{0\})\catqot {\Bbb G}_m(k).
$$
%%%
The coordinate algebra of ${\Bbb H}$ is ${\rm Sym}^{\star}(k^r\otimes_k V)^G$. For $s>0$, we set
$$
{\Bbb W}_s:=\bigoplus_{i=1}^s {\Bbb U}_i,\q {\Bbb U}_i:= \bigl({\rm Sym}^i(k^r\otimes_k V)^G\bigr)^\vee,\q i\ge 0.
$$
If $s$ is so large that $\bigoplus_{i=0}^s {\rm Sym}^i(k^r\otimes_k V)^G$ contains a set of generators for the
algebra ${\rm Sym}^{\star}(k^r\otimes_k V)^G$, then we have a $\GL_r(k)$-equivariant surjection of algebras
$$
{\rm Sym}^{\star}({\Bbb W}_s^\vee)\lra {\rm Sym}^{\star}(k^r\otimes_k V)^G,
$$
and, thus, a $\GL_r(k)$-equivariant embedding
$$
\iota_s\colon {\Bbb H}\hookrightarrow {\Bbb W}_s.
$$
\par
Set ${\Bbb I}:={\rm Isom}(k^r, V^\vee)/G$ ($\cong{\rm GL}_r(k)/G$). This is a dense open subset of ${\Bbb H}$.
The semistability of points $\iota_s(h)$, $h\in {\Bbb H}$, with respect to the action
of the \bfseries special \rm linear group $\SL_r(k)$ is described by the following result.
%%%
\begin{Lem}
\label{lem1:SemStab1001}
{\rm i)} Every point $\iota_s(i)$, $i\in {\Bbb I}$, is $\SL_r(k)$-polystable.
\par
{\rm ii)} A point $\iota_s(h)$, $h\in {\Bbb H}\setminus{\Bbb I}$, is not $\SL_r(k)$-semistable.
\end{Lem}
%%%
\begin{proof}[Proof (compare Lemma {\rm 4.1.1} in {\rm \cite{Schmittbesser}}).]
Ad i). We choose a basis for $V^\vee$. This provides us with the $({\SL}_r(k)\times G)$-invariant
function ${\frak d}\colon \HOM(k^r,V^\vee)\lra k$, $f\lma \det(f)$, which descends to
a (non-constant)
function on ${\Bbb H}$, called again ${\frak d}$. For any $i\in {\Bbb I}$, we clearly have
${\frak d}(\iota_s(i))\neq 0$, so that $\iota_s(i)$ is $\SL_r(k)$-semistable.
Furthermore, for any $f\in {\rm Isom}(k^r,V^\vee)$, the $({\SL}_r(k)\times G)$-orbit of $f$
is just a level set ${\frak d}^{-1}(z)$ for an appropriate $z\in {\Bbb G}_m(k)$. In particular, it is
closed. The image of this orbit is the $\SL_r(k)$-orbit of $i:=[f]$ in ${\Bbb H}$ which is,
therefore, closed. Since $\iota_s$ is a closed, $\SL_r(k)$-equivariant embedding, the
orbit of $\iota_s(i)$ is closed, too.
\par
Ad ii). It is obvious from the construction that the ring of $\SL_r(k)$-invariant functions
on ${\Bbb H}$ is generated by ${\frak d}$. This makes the asserted property evident.
\end{proof}
%%%
A key result is now the following:
%%%
\begin{Prop}
\label{prop1:Recognize1}
Fix a basis for $V$ in order to obtain a faithful representation $\rho\colon G\lra \GL_r(k)$ and a
$\GL_r(k)$-equivariant isomorphism
$$
\phi\colon {\GL}_r(k)/G\lra {\rm Isom}(k^r, V^\vee)/G.
$$
Suppose that $x=\iota_s(i)$ for some $i=\phi(g)\in {\Bbb I}$. Then, for a one parameter subgroup
$\la\colon {\Bbb G}_m(k)
\lra \SL_r(k)$, the following conditions are equivalent:
\begin{itemize}
\item[{\rm i)}] $\mu_{\kappa_s}(\la,x)=0$, $\kappa_s$ being the representation of $\SL_r(k)$ on ${\Bbb W}_s$.
\item[{\rm ii)}] There is a one parameter subgroup $\la^\p\colon {\Bbb G}_m(k)\lra g\cdot G\cdot g^{-1}$
with
$$
\bigl(V_\bullet(\la),\alpha_\bullet(\la)\bigr)=\bigl(V_\bullet(\la^\p),{\alpha}_\bullet(\la^\p)\bigr).
$$
\end{itemize}
\end{Prop}
%%%
\begin{proof}
We may clearly assume $g={\Bbb E}_r$. We first show ``ii)$\Longrightarrow$i)". Since $G$ is the $\GL_r(k)$-stabilizer
of $x$, we have $\mu(\la^\p,x)=0$ for any one parameter subgroup $\la^\p\colon {\Bbb G}_m(k)\lra G$. Now, Formula (\ref{MuGleich}) implies the claim.
\par
We turn to the implication ``i)$\Longrightarrow$ii)". By Lemma \ref{lem1:SemStab1001}, i), there exists an element $g^\p\in \SL_r(k)$, such that
$$
x^\p:=\lim_{z\rightarrow\infty}\la(z)\cdot x=\phi(g^\p).
$$
By Proposition \ref{KraftKuttler}, we may choose $g^\p\in {\cal R}_u(Q_{\SL_r(k)}(\la))$.
In particular, the element $g^\p$ fixes the flag $V_\bullet(\la)$.
Since $\la$ fixes $x^\p$, it lies in $g^\p\cdot G\cdot {g^\p}^{-1}$.
Setting $\la^\p:={g^\p}^{-1}\cdot\la\cdot g^\p$, we obviously have
$
(V_\bullet(\la),{\alpha}_\bullet(\la))=(V_\bullet(\la^\p),{\alpha}_\bullet(\la^\p))$,
and $\la^\p$ is a one parameter subgroup of $G$.
\end{proof}
%%%
Next, we look at the categorical quotient
$$
\ol{\Bbb H}=
\Proj\bigl({\rm Sym}^{\star}(k^r\otimes_k V)^G\bigr).
$$
For any positive integer $d$, we define
$$
{\rm Sym}^{(d)}(k^r\otimes_k V)^G:=\bigoplus_{i=0}^\infty {\rm Sym}^{id}(k^r\otimes_k V)^G.
$$
Then,
by the Veronese embedding,
$$
\Proj\bigl({\rm Sym}^{\star}(k^r\otimes_k V)^G\bigr)\cong\Proj\bigl({\rm Sym}^{(d)}(k^r\otimes_k V)^G\bigr).
$$
We can choose $s$, such that
\begin{itemize}
\item[{\rm a)}] ${\rm Sym}^{\star}(k^r\otimes_k V)^G$ is generated by elements in degree $\le s$.
\item[{\rm b)}] ${\rm Sym}^{(s!)}(k^r\otimes_k V)^G$ is generated by elements in degree $1$, i.e.,
by the elements in the vector space ${\rm Sym}^{s!}(k^r\otimes_k V)^G$.
\end{itemize}
%%%
Set
\begin{equation}
\label{eq:Ves}
{\Bbb V}_s:={\Bbb V}_s(k^r):=\bigoplus_{(d_1,\dots,d_s):\atop d_i\ge 0, \sum i d_i=s!}
\Bigl({\rm Sym}^{d_1}\bigl((k^r\otimes_k V)^G\bigr)\otimes \cdots \otimes
{\rm Sym}^{d_s}\bigl({\rm Sym}^s(k^r\otimes_k V)^G\bigr)\Bigr).
\end{equation}
Obviously, there is a natural surjection ${\Bbb V}_s\lra {\rm Sym}^{s!}(k^r\otimes_k V)^G$ and, thus, a
surjection
$$
{\rm Sym}^{\star}({\Bbb V}_s)\lra {\rm Sym}^{(s!)}(k^r\otimes_k V)^G.
$$
This defines a closed and $\GL_r(k)$-equivariant embedding
$$
\ol{\iota}_s\colon \ol{\Bbb H}\hookrightarrow\P({\Bbb V}_s).
$$
We also define
$$
{\Oh}_{\ol{\Bbb H}}(s!):=\ol{\iota}_s^{\star}\bigl({\Oh}_{{\Pe}({\Bbb V}_s)}(1)\bigr).
$$
Note that
\begin{equation}
\label{eq1:TensorPower}
{\Oh}_{\ol{\Bbb H}}\bigl((s+1)!\bigr)={\Oh}_{\ol{\Bbb H}}(s!)^{\otimes (s+1)}.
\end{equation}
%%%
\begin{Lem}
\label{lem1:SemStab1002}
Let $s$ be a positive integer, such that {\rm a)} and {\rm b)} as above are satisfied, and
$f\in {\HOM}(k^r,\allowbreak V^\vee)$ a $G$-semistable point. Set $h:=\iota_s([f])$ and
$\ol{h}:=\ol{\iota}_s([f])$.
Then, for any one parameter subgroup $\la\colon {\Bbb G}_m(k)\lra \SL_r(k)$, we have
$$
\mu_{\kappa_s}(\la,h)>(=/<)\ 0\q\Longleftrightarrow\q \mu_{\sigma_s}(\la,\ol{h})>(=/<)\ 0.
$$
Here, $\kappa_s$ is the representation of $\SL_r(k)$ on ${\Bbb W}_s$, and $\sigma_s$ is the linearization of the
$\SL_r(k)$-action on $\ol{\Bbb H}$ in ${\Oh}_{\ol{\Bbb H}}(s!)$. In particular, $\ol{h}$ is $\SL_r(k)$-semistable
if and only if $f\in {\rm Isom}(k^r,V^\vee)$.
\end{Lem}
%%%
\begin{proof}
Note that we have the following commutative diagram:
$$
\xymatrix{
{{\HOM}(k^{r}, V^\vee)\catqot G}
\ar@{->}[d]_{{\Bbb G}_m(k){\hbox{-}}}^{{\rm quotient}}
\ar@{^{(}->}[r]^-{{\iota}_s} &
{{\Bbb W}_s\setminus\{0\}} \ar@{->}[d]^-{\alpha}
\\
{\Pe}\bigl({\HOM}(k^{r}, V^\vee)^\vee\bigr)\catqot G
\ar@{^{(}->}[r]^-{\ol{\iota}_s} &{{\Pe}({\Bbb V}_s)}.
}
$$
The morphism $\alpha$ factorizes naturally over the quotient with respect to the ${\Bbb G}_m(k)$-action on
${\Bbb W}_s$ which is given on ${\Bbb U}_i$ by scalar multiplication with $z^{-i}$, $i=1,\dots,s$, $z\in {\Bbb G}_m(k)$.
The explicit description of $\alpha$ is as follows: An element $(l_1,\dots,l_s)\in {\Bbb W}_s$ with
$$
l_i\colon {\rm Sym}^i(k^r\otimes_k V)^G\lra k,\q i=1,\dots,s,
$$
is mapped to the class
$$
\Biggl[\bigoplus_{\ul{d}=(d_1,\dots,d_s):\atop d_i\ge 0, \sum i d_i=s!} l_{\ul{d}}\Biggr]
\colon {\Bbb V}_s \lra k
$$
%%%
$$
l_{\ul{d}} \colon (u^1_1\cdot\cdots\cdot u^{d_1}_1)\cdot\cdots\cdot(u^1_s\cdot\cdots\cdot u^{d_s}_s)\lma \bigl(l_1(u^1_1)\cdot\cdots\cdot l_1(u_1^{d_1})\bigr)\cdot\cdots\cdot \bigl(l_s(u^1_s)\cdot\cdots\cdot l_s(u_s^{d_s})\bigr)
$$
on ${\rm Sym}^{d_1}((W\otimes{\C^r}^\vee)^G)\otimes \cdots \otimes{\rm Sym}^{d_s}({\rm Sym}^s(W\otimes{\C^r}^\vee)^G)$.  With this description, one easily sees
$$
\mu_{\kappa_s}(\la,h)>(=/<)\ 0\q\Longleftrightarrow\q \mu_{\sigma_s}\bigl(\la,\alpha(h)\bigr)>(=/<)\ 0
$$
for all $\la\colon {\Bbb G}_m(k)\lra \SL_r(k)$ and all
$h\in {\Bbb W}_s\setminus \{0\}$. Together with the above diagram, this implies the claim.
\end{proof}
%%%
\par
Let us conclude this section with a formula for the $\mu$-function. Note that
${\Bbb V}_s$ is a $\GL_r(k)$-submodule of
$$
{\Bbb S}_s:=\bigoplus_{(d_1,\dots,d_s):\atop d_i\ge 0, \sum i d_i=s!}
\Bigl({\rm Sym}^{d_1}\bigl(k^r\otimes_k V\bigr)\otimes \cdots \otimes
{\rm Sym}^{d_s}\bigl({\rm Sym}^s(k^r\otimes_k V)\bigr)\Bigr).
$$
Since ${\Bbb G}_m(k)$ is a linearly reductive group, the weight spaces inside ${\Bbb V}_s$ with respect
to any one parameter subgroup $\la\colon {\Bbb G}_m(k)\lra \SL_r(k)$ are the intersection of the weight
spaces for $\la$ inside ${\Bbb S}_s$ with the subspace ${\Bbb V}_s$.
\par
The module ${\Bbb S}_s$ is a quotient module of $(W^{\otimes s!})^{\oplus N}$, $W:=k^r$. Therefore, the weight spaces
inside ${\Bbb S}_s$ with respect
to a one parameter subgroup $\la\colon {\Bbb G}_m(k)\lra \SL_r(k)$ are the projections of
the corresponding weight spaces in $(W^{\otimes s!})^{\oplus N}$. The latter may be easily described.
Given a one parameter subgroup $\la\colon {\Bbb G}_m(k)\lra \SL_r(k)$, we obtain the decomposition
%%%
$$
W=\bigoplus_{i=1}^{t+1} W^i
$$
into eigenspaces and the corresponding weights $\gamma_1<\cdots<\gamma_{t+1}$.
Set $I:=\{\, 1,\dots,t+1\,\}^{\times s!}$, and for $(i_1,\dots,i_{s!})\in I$ define
$$
W_{i_1,\dots,i_{s!}}:= W^{i_1}\otimes\cdots\otimes W^{i_s}.
$$
Then, all the weight spaces inside $(W^{\otimes s!})^{\oplus N}$ may be written as direct sums
of some subspaces of the form $(W_{i_1,\dots,i_{s!}})^{\oplus N}$.
\par
If we let $W^{\bullet}_{i_1,\dots,i_{s!}}$ be the image of $(W_{i_1,\dots,i_{s!}})^{\oplus N}$ in ${\Bbb S}_s$
and $W^\star_{i_1,\dots,i_{s!}}$ the intersection of $W^\bullet_{i_1,\dots,i_{s!}}$ with ${\Bbb V}_s$,
$(i_1,\dots,i_{s!})\in I$,
we understand that all the weight spaces inside ${\Bbb V}_s$ are direct sums of subspaces of the form
$W^\star_{i_1,\dots,i_{s!}}$, $(i_1,\dots,i_{s!})\in I$.
\par
This enables us to compute the weights in ${\Bbb V}_s$ in terms of the weighted flag
$(W_\bullet(\la),\allowbreak \alpha_\bullet(\la))$. Indeed, we define
$$
W_{i_1}\cdot\cdots\cdot W_{i_{s!}}
$$
as the image of $(W_{i_1}\otimes\cdots\otimes W_{i_{s!}})^{\oplus N}$ under the projection map $(W^{\otimes s!})^{\oplus N}
\lra {\Bbb S}_s$
and
$$
W_{i_1}\star\cdots\star W_{i_{s!}}:=(W_{i_1}\cdot\cdots\cdot W_{i_{s!}})\cap {\Bbb V}_s,\q (i_1,\dots,i_{s!})\in I.
$$
Altogether, we compute for $[l]\in {\Pe}({\Bbb V}_s)$ and $\la\colon {\Bbb G}_m(k)
\lra\SL_r(k)$ with weighted flag $(W_\bullet(\la),\alpha_\bullet(\la))$ as before
%%%%
\begin{equation}
\label{eq1:NochEinKomp0}
\mu_{\sigma_s}\bigl(\la, [l]\bigr) =
-\min\bigl\{\,\gamma_{i_1}+\cdots+\gamma_{i_{s!}}\,|\,
 (i_1,\dots,i_{s!})\in I:
l_{|W_{i_1}\star\cdots\star W_{i_{s!}}}\not\equiv 0\,\bigr\}.
\end{equation}
%%%
\subsubsection{Good quotients}
%%%
Suppose the algebraic group $G$ acts on the scheme $X$. In the framework of his GIT, Mumford defined the notion of a
\emph{good quotient} \cite{GIT}. Moreover, a \it universal (uniform) categorical quotient \rm is a categorical
quotient $(Y,\phi)$ for $X$ with respect to the action of $G$, such that, for every (every flat) base change
$Y^\p\lra Y$, $Y^\p$ is the categorical quotient for $Y^\p\times_Y X$ with respect to the induced $G$-action.
In particular, $Y\times Z$ is the categorical quotient for $X\times Z$ with respect to the given $G$-action on the
first factor.
%%%
\begin{Ex}
\label{uniform} i) Mumford's GIT produces good, uniform
(universal, if ${\rm Char}(k)=0$)
categorical quotients (see \cite{GIT}, Theorem 1.10, page 38).
\par
ii) If $G$ and $H$ are algebraic groups and we are given an action of $G\times H$ on the scheme $X$,
such that the good, universal, or uniform categorical quotients $X\catqot G$ and $(X\catqot G)\catqot H$ exist, then
$$
X\catqot (G\times H)=(X\catqot G)\catqot H.
$$
This follows from playing around with the universal property of a categorical quotient.
For good quotients, one might also use the argument from \cite{OST}.
\end{Ex}
%%%
The following lemma is well known (see \cite{Gie}, Lemma 4.6, and \cite{Ramanathan}, Lemma 5.1 (both in characteristic zero),
\cite{SeshGeomRed}, Theorem 2, (ii)). We recall the proof for the convenience of the reader.
%%%
\begin{Lem}
\label{lem:Gieseker/Ramanathan}
Let $G$ be a reductive linear algebraic group acting on the schemes $X_1$ and $X_2$, and let $\psi\colon X_1\lra X_2$ be an affine $G$-equivariant morphism. Suppose that there exists a good quotient $X_2\lra X_2\catqot G$. Then, there also exists a good quotient $X_1\lra X_1\catqot G$, and the induced morphism $\ol{\psi}\colon X_1\catqot G\lra X_2\catqot G$ is affine. Moreover, the following holds:
%%%
\begin{enumerate}
\item[\rm 1.] If $\psi$ is finite, then $\ol{\psi}$ is also finite.
\item[\rm 2.] If $\psi$ is finite and $X_2\catqot G$ is a geometric
quotient, then $X_1\catqot G$ is also a geometric quotient.
\end{enumerate}
\end{Lem}
%%%
\begin{proof}
If $X_2\catqot G$ is affine, then $X_1$ and $X_2$ are also affine,
and the existence of $X_1\catqot G$ is well known (see \cite{GIT}, Theorem
A.1.1). In general, the existence of $X_1\catqot G$ affine over
$X_2\catqot G$ is an easy exercise on gluing affine quotients (see
\cite{Ramanathan}, proof of Lemma 5.1).
The only non-trivial statement in the lemma is 1.
It follows from the last part of \cite{GIT}, Theorem
A.1.1. The point is that, if $\psi$ is finite, then $X_1$ is the
spectrum of the sheaf $\psi_\star({\Oh}_{X_1})$ of ${\Oh}_{X_2}$-algebras
which is coherent as an ${\Oh}_{X_2}$-module. Hence, by the theorem
cited above, $(\psi_\star {\Oh} _{X_1}) ^G$ is a coherent ${\Oh}
_{X_2/\hskip -2pt/ G}$-module, which is also an ${\Oh} _{X_2/\hskip -2pt/
G}$-algebra whose spectrum is $X_1\catqot G$. Hence, $\ol{\psi}$ is
a finite morphism.
\end{proof}
%%%
\subsection{Destabilizing one parameter subgroups}
%%%
We have seen in Lemma \ref{lem1:SemStab1001} that a point
$\iota_s(h)\in {\Bbb W}_s$ is not semistable for the $\SL_r(k)$-action, if
$h=[f]\in {\rm Hom}(k^r, V^\vee)\catqot G$ is the image of a homomorphism
$f\colon k^r\lra V^\vee$ which is not an isomorphism. This conclusion still
holds, if we replace $k$ by a non-algebraically closed ground field $K$.
What is, unfortunately, not automatic in positive characteristic is the fact
that there exists a one parameter subgroup $\la\colon {\Bbb G}_m(K)\lra
\SL_r(K)$ with $\mu_{\kappa_s}(\la, h)<0$ in this case.
This property is, however, needed in our approach to the semistable reduction
theorem for semistable singular principal $G$-bundles. Therefore,
we will now explain under which assumptions on the characteristic of the
ground field we will be able to verify the existence of a one parameter
subgroup $\la$ with $\mu_{\kappa_s}(\la, h)<0$.
%%%
\subsubsection{Preliminaries}
%%%
Let us collect two basic results. The first one is a generalization
of Kempf's results on the instability flag.
%%%
\begin{Thm}[Hesselink]
\label{thm:Hessel1}
Let $K$ be a not necessarily algebraically closed field. Suppose we are
given a representation $\kappa\colon \SL_r(K)\lra \GL(W)$ on the finite
dimensional $K$-vector space $W$, a point $h\in W$, a {\bfseries separable}
extension $\widetilde{K}/K$, and a one parameter subgroup
$\widetilde{\la}\colon {\Bbb G}_m(\widetilde{K})\lra \SL_r(\widetilde{K})$ with
$$
\mu(\widetilde{\la},h)<0.
$$
Then, there also exists a one parameter subgroup $\la\colon {\Bbb G}_m({K})\lra \SL_r({K})$ with $\mu({\la},h)<0$.
\end{Thm}
%%%
\begin{proof}
This follows from Theorem 5.5 in \cite{Hesselink2}.
\end{proof}
%%%
\begin{Prop}
\label{prop:Emb1}
Let $G$ be a reductive linear algebraic group and $\rho\colon G\lra \GL(W)$ a rational
representation. If $W'$ is a $G$-invariant subspace which possesses a direct complement as $G$-module,
then the categorical quotient $W'\catqot G$ embeds into $W\catqot G$.
\end{Prop}
%%%
\begin{proof}
We have to show that the surjection
$$
\iota^\#\colon {\rm Sym}^\star(W^\vee)=k[W]\lra k[W']={\rm Sym}^\star({W'}^\vee)
$$
of locally finite $G$-modules also induces a surjection on the algebras of invariant elements.
Our assumption says that $W=W'\oplus W''$ splits as a $G$-module. This shows that
${\rm Sym}^\star({W'}^\vee)$ embeds as a $G$-submodule into ${\rm Sym}^\star(W^\vee)$,
such that the restriction of $\iota^\#$ onto it is simply the identity. This easily yields the
claim.
\end{proof}
%%%
%\begin{Cor}
%\label{cor:Emb2}
%Let $G$ and $\rho$ be as before and $U\subset k^r$, $r:=\dim(V)$, a subspace. Then, the quotient ${\rm Hom}(V,U)\catqot G$ embeds into ${\rm Hom}(k^r, V^\vee)\catqot G$.
%\end{Cor}
%%%
%\begin{proof}
%Choose a direct complement $U'$, so that $k^r=U\oplus U'$ as $k$-vector space. The resulting
%decomposition
%$$
%{\rm Hom}(k^r, V^\vee)={\rm Hom}(V,U)\oplus {\rm Hom}(V,U')
%$$
%obviously is a decomposition of ${\rm Hom}(k^r, V^\vee)$ as a $G$-module, and we may apply
%Proposition \ref{prop:Emb1}.
%\end{proof}
%%%
\subsubsection{Digression on low height representations}
%%%
The general references for the following assertions are \cite{BP} and
\cite{SerreMour}.
Let $\rho\colon G\lra \GL(V)$ be a representation of the reductive linear
algebraic group $G$.
Then one attaches to $\rho$ its \it height \rm ${\rm ht}_G(\rho)$
(\cite{SerreMour}, p.\ 20; \cite{BP}, Definition 1) and its
\it separable index \rm $\psi_G(\rho)$ (\cite{BP}, Definition 6).
%%%
\begin{Rem}
By \cite{BP}, Remark 10, one has the estimate
$$
\psi_G(\rho)\le {\rm rank}(G)!\cdot {\rm ht}_G(\rho)^{{\rm rank}(G)}.
$$
\end{Rem}
%%%
We say that $\rho$ is a \it representation of low height \rm(\it low separable
index\rm), if ${\rm ht}_G(\rho)<p$ ($\psi_G(\rho)<p$).
(Of course, $p$ is the characteristic of the base field $k$.)
Here is a list of properties that representations of low height and low
separable index do enjoy.
%%%
\begin{Fact}
\label{Fact:LowHeight}
Assume that $\rho$ is a representation of low height.
\par
{\rm i)} The representation $\rho$ is semisimple.
\par
{\rm ii)} The stabilizer of any point $v\in V$ is a {\bfseries saturated}
subgroup of $G$. (See {\rm\cite{SerreMour}}, p.\ 22, for the definition of a saturated
subgroup.)
\par
{\rm iii)} Suppose that $\rho$ is also {\bfseries non-degenerate}, i.e.,
the connected component of the kernel is a torus.
Let $v\in V$ be a {\bfseries polystable} point. Then,
$\rho$ is also semisimple as a $G_v$-module.
\par
{\rm iv)} If $\rho$ is of low separable index, then the action of $G$
on $V$ that is induced by $\rho$ is separable.
\end{Fact}
%%%
\begin{proof}
Ad {\rm i)}. This is Theorem 6 in \cite{SerreMour}. Ad ii). The asserted
property is evident from the definition of a saturated subgroup given in
\cite{SerreMour}.
Ad {\rm iii)}. This property results from ii) and Theorems 8 and 9 in \cite{SerreMour}.
(Note that ${\rm ht}_G(\rho)<p$ also implies that the
Coxeter number $h_G$ of $G$ is at most $p$, by \cite{SerreMour}, p.\ 20.)
Ad {\rm iv)}. This is Theorem 7 in \cite{BP}.
\end{proof}
%%%
\subsubsection{Digression on the slice theorem}
%%%
The references for this section are the papers \cite{BR} and \cite{DZ}.
We assume that $\rho\colon G\lra \GL(V)$ is such that $\rho^\vee$ is
of low separable index and non-degenerate and look at the resulting action
of $G$ on ${\rm Hom}(k^r, V^\vee)$, $r:=\dim(V)$. (Note that the height and
separable index of ${\rm Hom}(k^r, V^\vee)\cong (V^\vee)^{\oplus r}$ agrees with the one
of $V^\vee$.) The results collected in Facts
\ref{Fact:LowHeight} imply that the slice theorem of Bardsley and Richardson
\cite{BR} may always be applied. Let
us review the formalism. Suppose that $f\in {\rm Hom}(k^r, V^\vee)$ is a
polystable point. Then, its stabilizer $G_f$ is a reduced, reductive,
and saturated closed subgroup of $G$, by Facts \ref{Fact:LowHeight}.
The tangent space $T_f(G\cdot f)$ to the orbit of $f$ at $f$ is
a $G_f$-submodule of ${\rm Hom}(k^r, V^\vee)$. By \ref{Fact:LowHeight}, iii),
we may find a $G_f$-submodule $N$ of ${\rm Hom}(k^r, V^\vee)$, such that
$$
{\rm Hom}(k^r, V^\vee)=T_f(G\cdot f)\oplus N
$$
as $G_f$-module. Then, Proposition 7.4 of \cite{BR} asserts that
there is a function $h\in k[N]^{G_f}$, such that $S:=N_h$ is an
{\bfseries \'etale slice at $f$}, i.e., we have the cartesian diagram
$$
\xymatrix{
G\times^{G_f} S \ar[r]^-\psi \ar[d]& {\rm Hom}(k^r, V^\vee)\ar[d]
\\
S\catqot G_f \ar[r]^-{\psi/\hskip -2pt/  G} & {\rm Hom}(k^r, V^\vee)\catqot G
}
$$
in which $\psi$ and $\psi\catqot G$ are \'etale morphisms.
\par
Next, we discuss the stratification given in \cite{DZ}, \S 2.
To this end, let ${\cal T}$ be the set of conjugacy classes
of stabilizers of closed points in ${\rm Hom}(k^r, V^\vee)$.
We say that a point $f\in {\rm Hom}(k^r, V^\vee)$ is \it of type
$\tau\in {\cal T}$\rm, if the stabilizer
$G_{f'}$ of a point $f'\in \ol{G\cdot f}$ with closed orbit belongs to $\tau$. If $\phi\in {\rm Hom}(k^r, V^\vee)\catqot G$,
then the \it type of $\phi$ \rm is the type of a point $f$ with
closed orbit that maps to $\phi$ under the quotient morphism.
For $\tau\in {\cal T}$, we set
$$
{\rm Hom}(k^r, V^\vee)_\tau:=\bigl\{\,f\in {\rm Hom}(k^r, V^\vee)\,|\,
f \hbox{ is of type }\tau\,\bigr\}
$$
and
$$
\bigl({\rm Hom}(k^r, V^\vee)\catqot G\bigr)_\tau:=\bigl\{\,\phi\in
{\rm Hom}(k^r, V^\vee)\catqot G\,|\,
\phi \hbox{ is of type }\tau\,\bigr\}.
$$
Similarly, we define $S_\tau$ and $(S\catqot G_f)_\tau$, if $S$ is an
\'etale slice as above.
Finally, for $\nu,\tau\in {\cal T}$, we write $\nu\preceq \tau$,
if there are points $f$ and $l$ in ${\rm Hom}(k^r, V^\vee)$ of type
$\nu$ and $\tau$, respectively, such that $G_f\supseteq G_l$.
By \cite{DZ}, Proposition 2.4 and 2.5, we have the following result.
%%%
\begin{Prop}
\label{prop:Strat}
For any $\tau\in {\cal T}$, the set $({\rm Hom}(k^r, V^\vee)\catqot G)_\tau$
is an irreducible locally closed subset of ${\rm Hom}(k^r, V^\vee)\catqot G$
with
$$
\ol{\bigl({\rm Hom}(k^r, V^\vee)\catqot G\bigr)_\tau}
=
\bigcup_{\nu \preceq\tau} \bigl({\rm Hom}(k^r, V^\vee)\catqot G\bigr)_\nu.
$$
\end{Prop}
%%%
The last statement which we are going to need is the following:
%%%
\begin{Prop}
\label{prop:Etale}
Let $f\in {\rm Hom}(k^r, V^\vee)$ be a point with closed orbit of type $\tau\in {\cal T}$.
\par
{\rm i)} The morphism $\psi\catqot G$ induces an \'etale morphism
$$
(S\catqot G_f)_\tau\lra (X\catqot G)_\tau.
$$
\par
{\rm ii)} The natural map
$$
\eta\colon S^{G_f}\lra (S\catqot G_f)_\tau
$$
is an isomorphism.
\end{Prop}
%%%
\begin{proof}
Part i) is Proposition 2.6, i), of \cite{DZ}. Moreover, that proposition also shows that $\eta$ is
a bijection.
Now, $N^{G_f}$ is a $G_f$-submodule of $N$. Since ${\rm Hom}(k^r, V^\vee)$ and hence $N$ is a semisimple
$G$-module, Proposition \ref{prop:Emb1} proves that the natural map
$$
N^{G_f}= N^{G_f}\catqot G_f\lra N\catqot G_f
$$
is a closed embedding. Finally, by construction of the \'etale slice, we
have the cartesian diagram
$$
\xymatrix{
{S=N_h} \ar@{>->}[r]\ar[d] & N \ar[d]
\\
S\catqot G_f=(N\catqot G_f)_h \ar@{>->}[r] & N\catqot G_f
}
$$
in which the horizontal maps are open embeddings. Using this diagram, one easily infers our claim.
\end{proof}
%%%
\subsubsection{Finding the destabilizing one parameter subgroup}
%%%
%We put ourselves in the context of Corollary \ref{cor:Emb2} and
Let $Y$ be any (irreducible)
quasi-projective variety, and $E$ a vector bundle of rank $r$ on $Y$.
We define as usual
$$
m\colon {\cal H}:={\cal H}{\it om}(E, V^\vee\otimes{\Oh}_Y)\lra Y
$$
and let
$$
\ol{m}\colon \ol{\cal H}:={\cal H}\catqot G\lra Y
$$
be its good whence categorical quotient. Set $\pi\colon {\cal H}\lra \ol{\cal H}$ to be the quotient
map.
%%%
\begin{Thm}
Let
$$
\sigma\colon Y \lra {\cal H}{\it om}(E, V^\vee\otimes {\Oh}_Y)\catqot G
$$
be a section. Then, there are a finite separable extension $K$ of the function field $k(Y)$ of $Y$ and
a $K$-valued point $\widetilde{\eta}\in {\cal H}$, such that
$$
\pi(\widetilde{\eta})=\sigma(\eta),\q \eta\hbox{ being the generic point of $Y$}.
$$
\end{Thm}
%%%
\begin{proof}
By shrinking $Y$, we may assume that the vector bundle $E$ is trivial, so that we define
$$
\sigma'\colon Y\stackrel{\sigma}{\lra} \ol{\cal H}\cong({\rm Hom}(k^r, V^\vee)\catqot G)\times Y \lra {\rm Hom}(k^r, V^\vee)\catqot G.
$$
Let $\tau\in {\cal T}$ be minimal with respect to ``$\preceq$", such that
$$
\sigma'(\eta)\in ({\rm Hom}(k^r, V^\vee)\catqot G)_\tau.
$$
Then, we may choose a point $f\in {\rm Hom}(k^r, V^\vee)$ with closed orbit, such that
$\pi(f)\in \sigma'(Y)\cap {\rm Hom}(k^r, V^\vee)_\tau$. Let $S$ be an \'etale slice at $f$. Then, by Proposition \ref{prop:Etale}, i),
we have the \'etale map
$$
{\rm et}\colon (S\catqot G_f)_\tau\lra ({\rm Hom}(k^r, V^\vee)\catqot G)_\tau,\q \tau\hbox{ the type of } f.
$$
By construction, $\sigma'(Y)$ meets ${\rm et}((S\catqot G_f)_\tau)$, so that, in particular,
$\sigma'(\eta)$ lies in the image of ${\rm et}$. Hence, we find a finite separable extension $K$ of
$k(Y)$ and a $K$-valued point of $(S\catqot G_f)_\tau$ that maps to $\sigma'(\eta)$ under ${\rm et}$.
We now conclude with Proposition \ref{prop:Etale}, ii).
\end{proof}
%%%
\begin{Cor}
\label{cor:InstabFlag}
In the above situation, look at the embedding
$$
\iota_s\colon {\rm Hom}\bigl(k(Y)^r,V^\vee\otimes_k k(Y)\bigr)\catqot \bigl(G\times_{{\rm Spec}(k)}
{\rm Spec}(k(Y))\bigr)\hookrightarrow {\Bbb W}_s\otimes_k k(Y).
$$
If $\sigma(\eta)$ is not an element of ${\rm Isom}(k(Y)^r,V^\vee\otimes_k k(Y))\catqot (G\times_{{\rm Spec}(k)}
{\rm Spec}(k(Y))$, then there is a one parameter subgroup which defined over $k(Y)$ and which destabilizes $(\iota_s\circ \sigma)(\eta)$.
\end{Cor}
%%%
\begin{proof}
The point $\pi(\widetilde{\eta})$ gives a homomorphism
$$
h\colon K^r\lra V^\vee\otimes_k K
$$
whose kernel $B$ is non-trivial. There is a one parameter subgroup
$\widetilde{\la}\colon {\Bbb G}_m(K)\lra \SL_r(K)$ with weighted flag
$(0\subsetneq B\subsetneq K^r, (1))$. It satisfies $\mu(\widetilde{\la}, h)<\dim(B)-r<0$.
One easily sees that also
$$
\mu_{\kappa_s}\bigl(\widetilde{\la}, (\iota_s\circ \sigma)(\eta)\bigr)<0.
$$
The result therefore follows from Theorem \ref{thm:Hessel1}.
\end{proof}
%%%
\subsubsection{An improvement for adjoint groups}
%%%
Here, we assume that $G$ is an adjoint simple group of exceptional type.
Let $p$ be such that the adjoint representation of $G$ is of low height.
This implies that $p$ is also a good prime for $G$ (i.e., $p\neq 2,3$ for type $E_6$, $E_7$,
$F_4$, and $G_2$, and $p\neq 2,3,5$ for type $E_8$).
\par
By our previous discussion, it suffices to show that the action of $G$ on
${\rm Hom}(k^r,{\frak g}^\vee)$ is separable, ${\frak g}$ the Lie algebra of $G$. We recall
\begin{Thm}
\label{thm:VerGoodChar}
Under the assumption that the characteristic of $k$ is a good prime for $G$, the
Killing form on ${\frak g}$ is non-degenerate.
\end{Thm}
%%%
\begin{proof}
This follows from the computation of the determinant of the Killing form
in \cite{Sel}.
\end{proof}
%%%
Thus, we may split the $G$-module ${\rm End}(\frak g)$ as
${\frak m}\oplus {\frak g}$ and derive the ${\rm Ad}$-equivariant map
$$
\psi\colon G\stackrel{{\rm Ad}}{\lra} \GL({\frak g})\subset {\rm End}({\frak g})\lra {\frak g}
$$
with $(d\psi)_e={\id}_{\frak g}$. Note that ${\rm ad}$ is a left inverse to the last map.
\par
Now, let $Y\in {\frak g}$ be an element. Then, the Lie algebra of the (scheme theoretic) stabilizer $G_Y$
of $Y$ under the adjoint representation is
$$
{\frak g}_Y:=\bigl\{\, X\in {\frak g}\,|\, [X,Y]=0\,\bigr\}.
$$
On the other hand, we have the commutative diagram
\begin{equation}
\label{eq:Commute}
\xymatrix{
G \ar[r]^-{\rm Ad}\ar[d]_-{\psi} & \GL({\frak g})\ar[d]
\\
{\frak g} \ar[r]^-{{\rm ad}} & {\rm End}({\frak g}).
}
\end{equation}
We know
$$
{\rm ad}\bigl({\rm Ad}(g)(Y)\bigr)= {\rm Ad}(g)\cdot {\rm ad}(Y)\cdot {\rm Ad}(g)^{-1},\q
g\in G,\ Y\in {\frak g}.
$$
Thus, $g\in G_Y$ if and only if
$$
[{\rm Ad}(g),{\rm ad}(Y)]=0.
$$
By diagram (\ref{eq:Commute}), this is equivalent to
$$
[\psi(g), Y]=0.
$$
Therefore, under the \'etale morphism $\psi\colon G\lra {\frak g}$, the stabilizer $G_Y$ is the
preimage of the Lie algebra ${\frak g}_Y$. In particular, $G_Y$ is a reduced group
scheme. The same argument shows that the action of $G$ on ${{\frak g}}^r$ is separable.
By Theorem \ref{thm:VerGoodChar}, we see
%%%
\begin{Cor}
If the characteristic of the field $k$ is a good prime for $G$, then the action of $G$ on
${\rm Hom}(k^r,{\frak g}^\vee)\cong {\rm Hom}(k^r,{\frak g})$ is separable.
\end{Cor}
%%%
\begin{Rem}
i) Under the hypothesis that the characteristic is good, the isogeny $G\lra {\rm Ad}(G)$
is separable for all simple exceptional groups (see \cite{KW}, Chapter VI, Remark 1.7).
Thus, over curves we may use this result to deal with arbitrary simple groups.
\par
ii) For simple groups of type $A$, $B$, $C$, and $D$ and $k$ of very good characteristic (i.e., $\neq 2$ for
$B$, $C$, and $D$, and $n \not\equiv -1{\rm mod}\ p$ for $A_n$),
there also exist invariant scalar products on ${\rm End}({\frak g})$ which induce
non-degenerate forms on ${\frak g}$. These come from the trace form for the standard
representation of the respective classical group. Since an adjoint group is the product of its
simple factors, we get the result for all adjoint groups.
\end{Rem}
%%%
\subsection{Some $G$-linearized sheaves}
\label{sec:GLinearized}
%%%
Assume that $\rho\colon G\lra\GL(V)$ is any representation.
Let $B$ be a scheme and $\A$ a coherent ${\Oh}_B$-module. Equip $B$ with the trivial $G$-action.
We obtain the $G$-linearized
sheaf $\A\otimes V$. It follows easily from the universal property of the symmetric
algebra (\cite{EGA}, Section (9.4.1)) that ${\cal S}ym^{\star}(\A\otimes V)$ inherits a $G$-linearization.
Note that the algebra ${\cal S}ym^{\star}(\A\otimes V)$ is naturally graded and that the $G$-linearization
preserves this grading.
Let ${\cal S}ym^{\star}(\A\otimes V)^G$ be the sub-algebra of $G$-invariant elements in ${\cal S}ym^{\star}(\A\otimes V)$.
The $G$-linearization provides a $\pi$-invariant action of $G$ on
$$
\Hom\bigl(\A,V^\vee\otimes{\Oh}_B\bigr):=\Spec\bigl({{\cal S}ym}^{\star}(\A\otimes V)\bigr),
$$
$\pi\colon {\cal S}ym^{\star}(\A\otimes V)\lra B$ being the natural projection.
Then, the categorical quotient of the scheme
$\Spec(\allowbreak {\cal S}ym^{\star}(\A\otimes V))$ by the $G$-action is given
through
$$
\Spec\bigl({{\cal S}ym}^{\star}(\A\otimes V)\bigr)\catqot G=
\Spec\bigl({{\cal S}ym}^{\star}(\A\otimes V)^G\bigr)\stackrel{\ol{\pi}}{\lra} B.
$$
In characteristic zero, the construction commutes with base change. In positive characteristic,
we have to be more careful. Let $f\colon A\lra B$ be a morphism
of schemes. The natural isomorphism
$$
f^\star\bigl({\Sym}^\star(\A\otimes V)\bigr)\lra {\Sym}^\star\bigl(f^\star(\A)\otimes V\bigr)
$$
of $G$-linearized sheaves gives rise to the homomorphism
$$
{\rm bc}(f)\colon f^\star\bigl({\Sym}^\star(\A\otimes V)^G\bigr)\lra
{\Sym}^\star\bigl(f^\star(\A)\otimes V\bigr)^G.
$$
%%%
\begin{Lem}
\label{UnivCat}
If $\A$ is locally free, then ${\rm bc}(f)$ is an isomorphism.
\end{Lem}
%%%
\begin{proof}
If $\A$ is locally free of rank $r$, then $\Sym^\star(\A\otimes V)^G$ is the algebra
that is associated to $\A$ and the $\GL_r(k)$-module ${\rm Sym}^\star(k^r\otimes V)^G$, and the assertion
is clear.
\end{proof}
%%%
We also note the following property.
%%%
\begin{Lem}
\label{SurJect}
Let $\psi\colon \A^\p\lra \A$ be a surjective map of ${\Oh}_{B}$-modules. Assume that $\A'$ and $\A$
are locally free. Then, the induced homomorphism
$$
{{\cal S}ym}^{\star}(\A^\p\otimes V)^G\lra {{\cal S}ym}^{\star}(\A\otimes V)^G
$$
of ${\Oh}_B$-algebras is surjective as well.
\end{Lem}
%%%
\begin{proof}
This follows like Lemma \ref{UnivCat}, taking into account Proposition \ref{prop:Emb1}.
\end{proof}
%%%
\subsection{Polynomial representations}
\label{subsec:Poly}
%%%
A representation $\kappa\colon \GL_r(k)\lra \GL(U)$ is called a \it polynomial representation\rm, if it extends to a (multiplicative)
map $\ol\kappa\colon M_r(k)\lra {\rm End}(U)$. We say that $\kappa$ is \it homogeneous of degree $u\in\Z$\rm, if
$$
\kappa(z\cdot {\Bbb E}_r)=z^u\cdot{\id}_U,\q\forall z\in {\Bbb G}_m(k).
$$
Let $P(r,u)$ be the abelian category of homogeneous polynomial representations of $\GL_r(k)$ of degree $u$.
It comes with the duality functor 
\begin{eqnarray*}
^\star\colon P(r,u) &\lra& P(r,u)
\\
\kappa &\lma& (\kappa\circ {\id}_{\GL_r(k)}^\vee)^\vee.
\end{eqnarray*}
%%%
Here, $.^\vee$ stands for the corresponding dual representation. An example for a representation
in $P(r,u)$ is the \it $u$th divided power \rm $({\rm Sym}^u({\id}_{\GL_r(k)}))^\star$, i.e., 
the representation of $\GL_r(k)$ on
$$
D^u(W):=\bigl({\rm Sym}^u(W^\vee)\bigr)^\vee,\q W:=k^r.
$$
More generally, we look, for $u,v>0$, at the $\GL_r(k)$-module
\begin{equation}
\label{eq:DivPower}
{\Bbb D}^{u,v}(W):=\bigoplus_{(u_1,\dots,u_v):\atop u_i\ge 0, \sum_{i=1}^v u_i=u} \bigl(D^{u_1}(W)\otimes\cdots
\otimes D^{u_v}(W)\bigr).
\end{equation}
%%%
\begin{Lem}
\label{lem:Normalf}
Let $\kappa\colon \GL_r(k)\lra \GL(U)$ be a homogeneous polynomial representation of degree $u$.
Then, there exists an integer $v>0$, such that $U$ is a quotient of the $\GL(U)$-module
${\Bbb D}^{u,v}(W)$.
\end{Lem}
%%%
\begin{proof}
For the proof, we refer to \cite{HeinS}.
\end{proof}
%%%
\subsection{Extension of the structure group}
\label{subsec:boundedness}
%%%
We will need Theorem 8.4 of \cite{Lang2}:
%%%
\begin{Thm}
Let $G$ be a connected reductive group and $\rho\colon G\lra \GL(V)$ a representation which maps the radical of $G$ to the center of $\GL(V)$.  Then, there are the following cases:
\par
{\rm i)} Assume either ${\rm Char}(k)=0$ or
$\mu_{\max}(\Omega_X)\le 0$. If $(U,{\cal P})$ is a Ramanathan semistable
rational principal $G$-bundle on $X$ and $(U,\E)$ is the rational
vector bundle with fiber $V$ associated to $(U,{\cal P})$, then
$(U,\E)$ is (strongly) slope semistable.
\par
{\rm ii)} If  ${\rm Char}(k)=p>0$ and $\mu_{\max}(\Omega_X)> 0$,
there is a constant $C(\rho)$, depending only on $\rho$, such that
for any Ramanathan semistable rational principal $G$-bundle $(U,{\cal P})$ on
$X$ with associated rational vector bundle $(U,\E)$, one finds
$$
0\le \mu_{\max}(\E)-\mu_{\min}(\E)\le L_{\max}(\E)-L_{\min}(\E)\le C(\rho)\cdot
\frac{[L_{\max}(\Omega_X)]_+}{p}.
$$
\end{Thm}
%%%
\begin{Cor}
 \label{cor:ExtofSG}
Let $G$ be a connected reductive group and $\rho\colon G\lra \GL(V)$ a representation which maps the radical of $G$ to the center of $\GL(V)$. There is a constant $D(\rho)$ which depends only on $\rho$, such that 
$$
\mu(\E)-D(\rho)\le \mu_{\min}(\E)
$$
for any semistable rational principal $G$-bundle $(U,{\cal P})$ on
$X$ with associated rational vector bundle $(U,\E)$.
\end{Cor}
%%%
\begin{proof}
One has
$$
\mu(\E)\le \mu_{\max}(\E)=\bigl(\mu_{\max}(\E)-\mu_{\min}(\E)\bigr)+\mu_{\min}(\E).
$$
The assertion follows directly from this and the previous theorem. 
\end{proof}
%%%
\subsection{An extension property}
%%%
\begin{Prop}
\label{Extension}
Let $S$ be a scheme and $\F_S$ a vector bundle
on $S\times X$. Let $Z\subset S\times X$ be a closed subset, such
that ${\rm codim}_X(Z\cap (\{s\}\times X))\ge 2$ for every point
$s\in S$. Denote by $\iota\colon U:= (S\times X)\setminus
Z\subseteq S\times X$ the inclusion. Then, the natural map
$$
\F_S\lra \iota_{\star}(\F_{S|U})
$$
is an isomorphism.
\end{Prop}
%%%
\begin{proof}[Proof (after Maruyama {\rm \cite{Maru}}, page 112)]
Since this is a local question, we may clearly assume
$\F_S={\Oh}_{S\times X}$. Note that $Z$ is ``stable under
specialization" in the sense of \cite{EGAIV2}, (5.9.1), page 109.
By \cite{EGAIV2}, Theorem (5.10.5), page 115, one has to show that
$\inf_{x\in Z}{\rm depth}({\Oh}_{S\times X, x})\ge 2$. Since $X$ is
smooth, the morphism $\pi_S\colon S\times X\lra X$ is smooth.
Thus, by \cite{EGAIV4}, Proposition (17.5.8), page 70,
$$
\dim({\Oh}_{S\times X, x})-{\rm depth}({\Oh}_{S\times X, x})= \dim({\Oh}_{S,s})-{\rm depth}({\Oh}_{S,s})
$$
for every point $x\in S\times X$ and $s:=\pi_S(x)$.
This implies
\begin{equation}
\label{depth}
{\rm depth}({\Oh}_{S\times X, x})\ge \dim({\Oh}_{S\times X, x})-\dim({\Oh}_{S,s})=\dim({\Oh}_{\pi_S^{-1}(s),x}).
\end{equation}
Since for any point
$x\in \pi_S^{-1}(s)$, one has $\dim({\Oh}_{\pi_S^{-1}(s),x})={\rm codim}_{\pi_S^{-1}(s)}(\ol{\{x\}})$,
we derive the desired estimate ${\rm depth}({\Oh}_{S\times X, x})\ge 2$ for every point $x\in Z$ from the fact that
${\rm codim}_{\pi_S^{-1}(s)}(Z\cap \pi_S^{-1}(s))\ge 2$ and (\ref{depth}).
\end{proof}
%%%
\begin{Cor}
\label{NatTrans}
Suppose $S$ is a scheme, $\E_S$ is a coherent ${\Oh}_{S\times X}$-module,
and $\F_S$ is a locally free sheaf on $S\times X$. Let
$U\subseteq S\times X$ be an open subset whose complement $Z$ satisfies
 ${\rm codim}_X(Z\cap (\{s\}\times X))\ge 2$ for every point
$s\in S$. Then, for any homomorphism $\widetilde{\phi}_S\colon
\E_{S|U}\lra \F_{S|U}$, there is a unique extension
$$
\phi_S\colon \E_S\lra \F_S
$$
to $S\times X$. In particular, for a base change morphism $f\colon T\lra S$, we have
$$
\phi_T=(f\times {\id}_X)^{\star}(\phi_S).
$$
Here, $\phi_T$ is the extension of ${(f\times {\id}_X)_{|(f\times\id_X)^{-1}(U)}}^{\star}(\widetilde{\phi}_S)$.
\end{Cor}
%%%
\begin{proof}
An extension is given by
$$
\phi_S\colon \E_S\lra
\iota_{\star}(\E_{S|U})\stackrel{\iota_{\star}(\widetilde{\phi}_S)}{\lra}
\iota_{\star}(\F_S)\stackrel{\rm Proposition\
\ref{Extension}}{=} \F_S.
$$
Since $\E_S$ can be written as the quotient of a locally free sheaf, the uniqueness also follows
from Proposition \ref{Extension}. The final statement is clearly a consequence of the uniqueness property.
\end{proof}
%%%
\section{Fundamental results on semistable singular principal bundles}
%%%
After reviewing several elementary properties, we show that a singular principal $G$-bundle $(\A,\tau)$ is
slope semistable in the sense which has been defined in the introduction if and only if the associated
rational principal $G$-bundle $(U,{\cal P}(\A,\tau))$ is semistable in the sense of Ramanathan.
%%%
\subsection{The basic formalism of singular principal bundles}
\label{subsec:Pseudo}
%%%
Since $G$ is a semisimple group, the basic formalism of pseudo
$G$-bundles in positive characteristic is exactly the same as in
characteristic zero. Therefore, we may refer the reader to
\cite{Schmitt}, Section 3.1, for more details (be aware that in
this reference the term ``singular principal $G$-bundle" is used
for our ``pseudo $G$-bundle"). We fix a faithful representation
$\rho\colon G\lra \GL(V)$. Then, a \it pseudo $G$-bundle \rm
$(\A,\tau)$ consists of a torsion free coherent ${\Oh}_X$-module $\A$
of rank $\dim_k(V)$ with trivial determinant and a homomorphism
$\tau\colon {\cal S}ym^{\star}(\A\otimes V)^G\lra {\Oh}_X$ which is
non-trivial in the sense that it is not just the projection onto
the degree zero component. Let $U\subseteq X$ be the maximal open
subset where $\A$ is locally free. Since $\rho(G)\subseteq
\SL(V)$, we have the open immersion
$$
{\cal I}som(\A_{|U}, V^\vee\otimes{\Oh}_U)/G
\subset
\Hom(\A, V^\vee\otimes {\Oh}_X)\catqot G.
$$
Recall the following alternatives.
%%%
\begin{Lem}
\label{Alternative}
Let $(\A,\tau)$ be a pseudo $G$-bundle and
$$
\sigma\colon X\lra \Hom(\A, V^\vee\otimes {\Oh}_X)\catqot G
$$
the section defined by $\tau$. Then, either
$$
\sigma(U)\subset {\cal I}som(\A_{|U}, V^\vee\otimes{\Oh}_U)/G
$$
or
$$
\sigma(U)\subset \bigl(\Hom(\A_{|U}, V^\vee\otimes {\Oh}_U)\catqot G\bigr)
\big\backslash \bigl({\cal I}som(\A_{|U}, V^\vee\otimes{\Oh}_U)/G\bigr).
$$
\end{Lem}
%%%
\begin{proof}
See \cite{Schmitt}, Corollary 3.4.
\end{proof}
%%%
In the former case, we call $(\A,\tau)$ a \it singular principal $G$-bundle\rm.
We may form the base change diagram
$$
\xymatrix{
{\cal P}(\A,\tau) \ar[r]\ar[d] & \Hom(\A_{|U}, V^\vee\otimes {\Oh}_U)\ar[d]
\\
U \ar[r]^-{\sigma_{|U}} & \Hom(\A_{|U}, V^\vee\otimes {\Oh}_U)\catqot G.
}
$$
If $(\A,\tau)$ is a singular principal $G$-bundle, then ${\cal P}(\A,\tau)$ is a principal $G$-bundle
over $U$ in the usual sense, i.e., a rational principal $G$-bundle on $X$ in the sense of Ramanathan
(see Section \ref{ClassicalSemStab}).
\par
A {\it family of pseudo $G$-bundles parameterized by the scheme
$S$} is a pair $(\A_S,\tau_S)$ which consists of an $S$-flat
family $\A_S$ of torsion free sheaves on $S\times X$ and a
homomorphism $\tau_S\colon {\Sym}^\star(\A_S\otimes V)^G\lra
{\Oh}_{S\times X}$. We say that the family $(\A_S^1,\tau_S^1)$ is
{\it isomorphic} to the family $(\A_S^2,\tau_S^2)$, if there is an
isomorphism $\psi_S\colon \A_S^1\lra \A_S^2$, such that the
induced isomorphism ${\Sym}^\star(\A^1_S\otimes V)^G\lra
{\Sym}^\star(\A^2_S\otimes V)^G$ carries $\tau_S^2$ into $\tau_S^1$.
\par
We also need a more general looking concept.
A {\it pre-family of pseudo $G$-bundles parameterized by the scheme
$S$} is a pair $(\A_S,\tau'_S)$ which is composed of an $S$-flat
family $\A_S$ of torsion free sheaves on $S\times X$ and a
homomorphism $\tau'_S\colon {\Sym}^\star(\A_{S|U}\otimes V)^G\lra
{\Oh}_{U}$. Here, $U\subseteq S\times X$ is the maximal open subset where $\A_S$ is locally free.
The pre-family $(\A_S^1,{\tau'}_S^1)$ is
{\it isomorphic} to the pre-family $(\A_S^2,{\tau'}_S^2)$, if there is an
isomorphism $\psi_S\colon \A_S^1\lra \A_S^2$, such that the
induced isomorphism ${\Sym}^\star(\A^1_{S|U}\otimes V)^G\lra
{\Sym}^\star(\A^2_{S|U}\otimes V)^G$ transforms ${\tau'}_S^2$ into ${\tau'}_S^1$.
%%%
\begin{Lem}
\label{lem:AllowsBaseChange}
Let $S$ be a scheme of finite type over $k$.
Then, the assignment
%%%
\begin{eqnarray*}
\left\{\,\begin{array}{r}
                    \hbox{Isomorphy classes of}
                    \\
                    \hbox{families of pseudo $G$-bundles}
                    \\
                    \hbox{parameterized by $S$}
                    \end{array}
\,\right\}  &\lra& \left\{\,\begin{array}{l}
                    \hbox{Isomorphy classes of}
                    \\
                    \hbox{pre-families of pseudo $G$-bundles}
                    \\
                    \hbox{parameterized by $S$}
                    \end{array}
                    \,\right\}
\\
(\A_S,\tau_S) &\lma& (\A_S, \tau_{S|U})
\end{eqnarray*}
%%%
is a bijection.
\end{Lem}
%%%
\begin{proof}
If $(\A_S,\tau'_S)$ is a pre-family, denote by $U$ the maximal open subset where $\A_S$ is locally free
and by $\iota\colon U\lra S\times X$ the inclusion. Set
$$
\tau_S\colon {\Sym}^\star(\A_{S}\otimes V)^G\lra \iota_\star\bigl({\Sym}^\star(\A_{S|U}\otimes V)^G\bigr)
\stackrel{\iota_\star(\tau'_S)}{\lra} \iota_\star({\Oh}_U)\stackrel{\rm\ref{Extension}}{=}{\Oh}_{S\times X}.
$$
Then, $(\A_S,\tau_S)$ is a real family of pseudo $G$-bundles which maps to $(\A_S,\tau'_S)$ under the above map.
\par
It remains to verify injectivity. Let $\A_S$ be a flat family of torsion free coherent ${\Oh}_X$-modules
on $S\times X$. Over every affine open subset $W\subset S\times X$, the algebra $\Sym^\star(\A_{S|W}\otimes V)^G$
is finitely generated. Since $S\times X$ is according to our assumption quasi-compact, we may find an
$s>0$, such that $\Sym^\star(\A_S\otimes V)^G$ is generated by the coherent ${\Oh}_{S\times X}$-module
$$
{\Bbb W}_s(\A):=\bigoplus_{i=1}^s {\Sym}^i(\A_S\otimes V)^G.
$$
Since $\tau_S$ is determined by $\tau'_S:=\tau_{S|{\Bbb W}_s(\A)}$, it remains to show that
$\tau'_S$ is determined by its restriction to $U$. But this is an immediate consequence of
Corollary \ref{NatTrans}.
\end{proof}
%%%
By Lemma \ref{UnivCat}, we have a pullback operation on pre-families of pseudo $G$-bundles with respect
to base change morphisms $T\lra S$. Using the above arguments, we also obtain a pullback operation
for families of pseudo $G$-bundles.
%%%
\subsection{Semistable rational principal $G$-bundles}
\label{ClassicalSemStab}
%%%
We now review the formalism introduced by Ramanathan and compare it with our setup.
A \it rational principal $G$-bundle on $X$ \rm is a pair $(U,{\cal P})$ which consists of a big open subset
$U\subseteq X$ and a principal $G$-bundle ${\cal P}$ on $U$.
Such a rational principal $G$-bundle is said to be \it (semi)stable\rm, if for every open subset $U^\p\subseteq U$
which is big in $X$, every parabolic subgroup $P$ of $G$, every reduction $\beta\colon U^\p\lra {\cal P}_{|U^\p}/P$
of the structure group of ${\cal P}$ to $P$ over $U^\p$, and every antidominant character (see below) $\chi$ on $P$,
we have
$$
{\rm deg}\bigl({\cal L}(\beta,\chi)\bigr) (\ge)0.
$$
Note that the antidominant character $\chi$ and the principal $P$-bundle ${\cal P}_{|U^\p}\lra {\cal P}_{|U^\p}/P$
define a line bundle on ${\cal P}_{|U^\p}/P$. Its pullback to $U^\p$ via $\beta$
is the line bundle ${\cal L}(\beta,\chi)$. Since $U^\p$ is big in $X$, it makes sense to speak about the
degree of ${\cal L}(\beta,\chi)$.
\par
We fix a pair $(B,T)$ which consists of a Borel subgroup $B\subset G$ and a maximal torus $T\subset B$.
If $P$ and $P^\p$ are conjugate in $G$, a reduction $\beta$ of a principal $G$-bundle to $P$ may equally
be interpreted as a reduction to $P^\p$. Thus, it suffices to consider parabolic subgroups
of the type $P_G(\la)$, $\la\colon {\Bbb G}_m(k)\lra T$ a one parameter subgroup, which contain $B$.
Here, we use the convention (compare Remark \ref{rem:Spri})
$$
P_G(\la):=Q_G(-\la)=\bigl\{\, g\in G\,|\,
\lim_{z\rightarrow 0} \la(z)\cdot g\cdot \la(z)^{-1}\hbox{ exists in $G$}\,\bigr\}.
$$
Let $X_\star(T)$ and $X^\star(T)$ be the free $\Z$-modules of one parameter subgroups and characters of $T$,
respectively. We have the canonical pairing $\langle .,.\rangle\colon X_\star(T)\times X^\star(T)\lra \Z$.
Set $X_{\star,\K}(T):=X_\star(T)\otimes_\Z \K$ and $X^\star_\K(T):= X^\star(T)\otimes_\Z\K$, and let
$\langle .,.\rangle_\K\colon X_{\star,\K}(T)\times X^\star_\K(T)\lra \K$ be the $\K$-bilinear extension of
$\langle.,.\rangle$, $\K=\Q,\R$. Finally, suppose $(.,.)^\star\colon X^\star_\R(T)\times X^\star_\R(T)\lra \R$ is a scalar product which is invariant under the Weyl group $W(T)={\cal N}(T)/T$. This also yields the product
$(.,.)_\star\colon X_{\star,\R}(T)\times X_{\star,\R}(T)\lra \R$. We assume that $(.,.)^\star$ is defined over $\Q$.
\par
The datum $(B,T)$ defines the set of positive roots $R$ and the set $R^\vee$ of coroots (see \cite{Spri}).
Let ${\cal C}\subset X_{\star,\R}(T)$ be the cone spanned by the elements of $R^\vee$ and ${\cal D}\subset X^\star_\R(T)$ the dual cone of ${\cal C}$ with respect to $\langle.,.\rangle_\R$. Equivalently, the cone ${\cal D}$ may be characterized as being the dual cone of the cone spanned by the elements in $R$ with respect to $(.,.)^\star$. Indeed, one has
$\langle\, \cdot\,, \alpha^\vee\rangle_\R= 2 (\,\cdot\,,\alpha)^{\star}/(\alpha,\alpha)^\star$, $\alpha\in R$.
Now, a character $\chi\in X^\star(T)$ is called \it dominant\rm, if it lies in ${\cal D}$, and \it antidominant\rm, if $-\chi$ lies in ${\cal D}$. A character $\chi$ of a parabolic subgroup containing $B$ is called \it (anti)dominant\rm, if its restriction to $T$ is (anti)dominant.
\par
In the definition of semistability, we may clearly assume that $(\chi,\alpha)^\star>0$ for every $\alpha\in R$ with $\langle\la,\alpha\rangle>0$, if $P=P_G(\la)$. Otherwise, we may choose $\la^\p$, such that $(\chi,\alpha)^\star>0$ if and only if $\langle\la^\p,\alpha\rangle>0$. Then, $P_G(\la^\p)$ is a parabolic subgroup which contains $P_G(\la)$,
$\chi$ is induced by a character $\chi^\p$ on $P_G(\la^\p)$, the reduction $\beta$ defines the reduction $\beta^\p\colon {\cal P}_{|U^\p}/P_G(\la)\lra {\cal P}_{|U^\p}/P_G(\la^\p)$, and $\L(\beta,\chi)=\L(\beta^\p,\chi^\p)$.
Every one parameter subgroup $\la$ of $T$ defines a character $\chi_\la$ of $T$, such that $\langle \la,\chi\rangle=(\chi_\la,\chi)^\star$ for all $\chi\in X^\star(T)$. Finally, observe that the cone ${\cal C}^\p$ of one parameter subgroups $\la$ of $T$, such that $B\subseteq P_G(\la)$, is dual to the cone spanned by the roots. Thus,
$$
\forall \la\in X_\star(T):\q B\subseteq P_G(\la)\q\Longleftrightarrow\q \chi_\la\in {\cal D}.
$$
If one of those conditions is verified, then $(P_G(\la),\chi_\la)$ consists of a parabolic subgroup containing $B$ and a dominant character on it. Similarly, if $Q_G(\la)$ contains $B$, then $\chi_\la$ is an antidominant character on $Q_G(\la)$. Our discussion shows:
%%%
\begin{Lem}
A rational principal $G$-bundle $(U,{\cal P})$ is (semi)stable if and only if for every open subset $U^\p\subseteq U$
which is big in $X$, every non-trivial one parameter subgroup $\la\in {\cal C}^\p$, and every reduction $\beta\colon U^\p\lra {\cal P}_{|U^\p}/Q_G(\la)$ of the structure group of ${\cal P}$ to $Q_G(\la)$ over $U^\p$, we have
$$
{\rm deg}\bigl({\cal L}(\beta,\chi_\la)\bigr) (\ge)0.
$$
\end{Lem}
%%%
For any one parameter subgroup $\la$ of $G$, we may find a pair $(B^\p,T^\p)$ consisting of a Borel subgroup
$B^\p$ of $G$ and a maximal torus $T^\p\subset B^\p$, such that $\la\in X_\star(T^\p)$ and
$B^\p\subseteq Q_G(\la)$. Then, there exists an element $g\in G$, such that $(g\cdot B\cdot g^{-1},
g\cdot T\cdot g^{-1})= (B^\p, T^\p)$, and we obtain
\begin{eqnarray*}
{(.,.)^\star}^\p\colon X^\star_\R(T^\p)\times X^\star_\R(T^\p) &\lra& \R
\\
(\chi,\chi^\p) &\lma& \bigl(\chi(g^{-1}\cdot .\cdot g), \chi^\p(g^{-1}\cdot .\cdot g)\bigr)^\star.
\end{eqnarray*}
%%%
Since $(.,.)^\star$ is invariant under $W(T)$, the product ${(.,.)^\star}^\p$ does not depend on $g$.
In particular, we obtain a character $\chi_\la$ on $Q_G(\la)$. This character does not depend on
$(B^\p,T^\p)$ (see \cite{Schmittbesser}, (2.31)). We conclude
%%%
\begin{Lem}
A rational principal $G$-bundle $(U,{\cal P})$ is (semi)stable if and only if for every open subset $U^\p\subseteq U$ which is big in $X$, every non-trvial one parameter subgroup $\la\colon {\Bbb G}_m(k)\lra G$, and every reduction $\beta\colon U^\p\lra {\cal P}_{|U^\p}/Q_G(\la)$ of the structure group of ${\cal P}$ to $Q_G(\la)$ over $U^\p$, we have
$$
{\rm deg}\bigl({\cal L}(\beta,\chi_\la)\bigr) (\ge)0.
$$
\end{Lem}
%%%
Suppose that $\rho\colon G\lra \GL(V)$ is a faithful representation. We may assume that $T$ maps to the maximal torus $\widetilde{T}\subset \GL(V)$, consisting of the diagonal matrices. The character group $X^\star(\widetilde{T})$ is freely generated by the characters $e_i\colon {\rm diag}(\la_1,\dots,\la_n)\lma \la_i$, $i=1,\dots,n$. We define
%%%
\begin{eqnarray*}
(.,.)^\star_{\widetilde{T}}\colon  X^\star_\R(\widetilde{T})\times X^\star_\R(\widetilde{T})&\lra& \R
\\
\Bigl(\sum_{i=1}^n x_i\cdot e_i, \sum_{i=1}^n y_i\cdot e_i\Bigr)&\lma & \sum_{i=1}^n x_i\cdot y_i.
\end{eqnarray*}
%%%
The scalar product $(.,.)_{\widetilde{T}}^\star$ is clearly defined over $\Q$ and invariant under the
Weyl group $W(\widetilde{T})$. The product $(.,.)_{\widetilde{T}}^\star$ therefore restricts to a scalar product $(.,.)^\star$ on $X^\star_\R(T)$ with the properties we have asked for. We find a
nice formula for ${\rm deg}({\cal L}(\beta,\chi_\la))$.
Indeed, if $(U,{\cal P})$ is a rational principal $G$-bundle, and if ${\E}$ is the vector bundle on $U$ associated to ${\cal P}$ by means of $\rho$, then we have, for every one parameter subgroup $\la\colon {\Bbb G}_m(k)\lra G$,
the embedding
$$
\iota\colon {\cal P}/Q_G(\la) \hookrightarrow {\cal I}som(V\otimes{\Oh}_X,\E)/Q_{\GL(V)}(\la).
$$
As usual, we obtain a weighted filtration $(V_\bullet(\lambda),\alpha_\bullet(\lambda))$ of $V$, and, for every reduction $\beta\colon U^\p\lra {\cal P}_{|U^\p}/\allowbreak Q_G(\la)$ over a big open subset $U^\p\subseteq U$, the reduction $\iota\circ \beta$ corresponds to a filtration
$$
\E_{\bullet}(\beta): 0\subsetneq \E_{1}\subsetneq \cdots\subsetneq \E_{t}\subsetneq \E_{|U^\p}
$$
by subbundles with $\rk(\E_{i})=\dim_k (V_{i})$, $i=1,\dots,t$. With the weighted filtration $(\E_{\bullet}(\beta),\allowbreak {\alpha}_\bullet(\la))$ of $\E_{|U^\p}$, we find
$$
{\rm deg}\bigl({\cal L}(\beta,\chi_\la)\bigr)= L\bigl(\E_{\bullet}(\beta),{\alpha}_\bullet(\la)\bigr)
=\sum_{i=1}^{t} \alpha_{i}\cdot \bigl(\rk(\E_{i})\cdot \deg(\E)-\rk(\E)\cdot \deg(\E_{i})\bigr).
$$
To see this, observe that the character $\chi_\la$ is, by construction, the restriction of a character $\chi$ of $\widetilde{T}$, so that
${\cal L}(\beta,\chi_\la)={\cal L}(\iota\circ\beta,\chi)$. The degree of the latter line bundle is computed in
Example 2.15 of \cite{SchmittGlobal} and gives the result stated above.
Thus, we conclude
%%%
\begin{Lem}
A rational principal $G$-bundle $(U,{\cal P})$ is (semi)stable if and only if, for every open subset $U^\p\subseteq U$ which is big in $X$, every non-trivial one parameter subgroup $\la\colon {\Bbb G}_m(k)\lra G$, and every reduction $\beta\colon U^\p\lra {\cal P}_{|U^\p}/Q_G(\la)$ of the structure group of ${\cal P}$ to $Q_G(\la)$ over $U^\p$, we have
$$
L\bigl(\E_{\bullet}(\beta),{\alpha}_\bullet(\la)\bigr)(\ge)0.
$$
\end{Lem}
%%%
\begin{Rem}
If $(U,{\cal P})$ is given as a singular principal $G$-bundle $({\A},\tau)$, then, in the notation of the introduction, we have
$$
\A_{|U}=\E^\vee,\ \A_{i|U}=\ker(\E^\vee\lra \E_{t+1-i}^\vee)\,\q\hbox{and}\q\alpha_{i;\beta}=\alpha_{t+1-i},\q i=1,\dots,t.
$$
Then, one readily verifies
$$
L({\A},\tau;\beta)=L\bigl(\E_\bullet(\beta),\alpha_\bullet(\la)\bigr).
$$
This proves that the definition of slope semistability (\ref{slOpe}) given in the introduction is the
original definition of Ramanathan. We have arrived at our notion of semistability, by replacing degrees
by Hilbert polynomials. Thus, our semistability concept is a ``Gieseker version" of Ramanathan semistability.
\end{Rem}
%%%
\section{Dispo sheaves}
\label{sec:DisShea}
%%%
In the papers \cite{Schmitt} and \cite{Schmittbesser}, the theory of decorated sheaves was used to construct projective moduli spaces
for singular principal $G$-bundles in characteristic zero. Due to the more difficult representation
theory of general linear groups in positive characteristic, this
approach is not available in (low) positive characteristic. Nevertheless, one may still associate to
any singular principal bundle a more specific object than a decorated sheaf, namely a so-called ``dispo
sheaf". The moduli theory of these dispo sheaves may be developed along the lines of the theory of decorated sheaves
in \cite{Schmitt0} and \cite{GS}, making several non-trivial modifications.
%%%
\subsection{The basic definitions}
%%%
For this section, we fix the representation $\rho\colon G\lra \SL(V)\subseteq \GL(V)$ and a positive integer $s$ as in Section
\ref{subsub:Spec}.
\par
Suppose that $\A$ is a coherent ${\Oh}_X$-module of rank $r:=\dim_k(V)$. Then, the \it sheaf of invariants in the
symmetric powers of $\A$ \rm is defined as
$$
{\Bbb V}_s(\A):=\bigoplus_{(d_1,\dots,d_s):\atop d_i\ge 0, \sum i d_i=s!}
\Bigl({\Sym}^{d_1}\bigl((\A\otimes_k V)^G\bigr)\otimes \cdots \otimes
{\Sym}^{d_s}\bigl({\Sym}^s(\A\otimes_k V)^G\bigr)\Bigr).
$$
A \it dispo\footnote{decorated with invariants in symmetric powers} sheaf \rm(\it of type $(\rho,s)$\rm) is a pair
$(\A,\phi)$ which consists of a torsion free sheaf $\A$ of rank $r$ with $\det(\A)\cong{\Oh}_X$ on $X$ and a non-trivial
homomorphism
$$
\phi\colon {\Bbb V}_s(\A)\lra {\Oh}_X.
$$
Two dispo sheaves $(\A_1,\phi_1)$ and $(\A_2,\phi_2)$ are said to
be \it isomorphic\rm, if there exists an isomorphism $\psi\colon \A_1\lra \A_2$,
such that, with the induced isomorphism ${\Bbb V}_s(\psi)\colon {\Bbb V}_s(\A_1)\lra {\Bbb V}_s(\A_2)$, one obtains
$$
\phi_1=\phi_2\circ {\Bbb V}_s(\psi).
$$
A \it weighted filtration $(\A_\bullet,{\alpha}_\bullet)$ of the torsion free sheaf $\A$ \rm consists of a filtration
$$
0\subsetneq \A_1\subsetneq\cdots\subsetneq \A_t\subsetneq \A_{t+1}= \A
$$
of $\A$ by saturated subsheaves and a tuple ${\alpha}_\bullet=(\alpha_1,\dots,\alpha_t)$ of positive rational numbers.
Given such a weighted filtration $(\A_\bullet,\alpha_\bullet)$, we introduce the quantities
%%%
\begin{eqnarray*}
M(\A_\bullet,{\alpha}_\bullet)&:=&\sum_{j=1}^t \alpha_j\cdot
\bigl(\rk(\A_j)\cdot P(\A)-\rk(\A)\cdot P(\A_j)\bigr),
\\
L(\A_\bullet,{\alpha}_\bullet)&:=&\sum_{j=1}^t \alpha_j\cdot
\bigl(\rk(\A_j)\cdot \deg(\A)-\rk(\A)\cdot \deg(\A_j)\bigr).
\end{eqnarray*}
%%%
Next, let $(\A,\phi)$ be a dispo sheaf and $(\A_\bullet,\alpha_\bullet)$ a weighted filtration of $\A$.
Fix a flag
$$
W_\bullet\colon 0\subsetneq W_1\subsetneq\cdots
\subsetneq W_t\subsetneq W:=k^r\q\hbox{with } \dim_k(W_i)=\rk(\A_i),\ i=1,\dots,t.
$$
We may find a small open subset $U$, such that
\begin{itemize}
\item $\phi_{|U}$ is a surjection onto ${\Oh}_U$;
\item there is a trivialization $\psi\colon \A_{|U}\lra W\otimes{\Oh}_U$ with $\psi(\A_{i|U})=W_i\otimes{\Oh}_U$, $i=1,\dots,t$.
\end{itemize}
%%%
In presence of the trivialization $\psi$, the homomorphism $\phi_{|U}$ provides us with the morphism
$$
\beta\colon U\lra {\Pe}\bigl({\Bbb V}_s(\A)\bigr)\stackrel{{\Bbb V_s}(\psi)}{\lra}{\Pe}({\Bbb V}_s)\times U\lra {\Pe}({\Bbb V}_s).
$$
(Consult Section \ref{SomeGIT} for the notation ``${\Bbb V}_s$".) Finally, let $\la\colon {\Bbb G}_m(k)\lra \SL(W)$
be a one parameter subgroup with $(W_\bullet,\alpha_\bullet)$ as its weighted flag. With these choices made, we
set
$$
\mu(\A_\bullet,\alpha_\bullet;\phi):=\max\bigl\{\,\mu_{\sigma_s}\bigl(\la,\beta(x)\bigr)\,|\, x\in U\,
\bigr\}.
$$
(The linearization $\sigma_s$ has been introduced in Lemma \ref{lem1:SemStab1002}.)
As in \cite{Schmitt0}, p.\ 176, one checks that the quantity $\mu(\A_\bullet,\alpha_\bullet;\phi)$
depends only on the data $(\A_\bullet,\alpha_\bullet)$ and $\phi$.
%%%
\begin{Rem}
\label{rem:DescSemStabAnotherWay}
i) Let us outline another, intrinsic definition of the number
$\mu(\A_\bullet,\alpha_\bullet;\phi)$.
First, observe that ${\Bbb V}_s(\A)$ is a submodule of
$$
{\Bbb S}_s(\A):=\bigoplus_{(d_1,\dots,d_s):\atop d_i\ge 0, \sum i d_i=s!}
\Bigl({\Sym}^{d_1}\bigl(\A\otimes_k V\bigr)\otimes \cdots \otimes
{\Sym}^{d_s}\bigl({\Sym}^s(\A\otimes_k V)\bigr)\Bigr)
$$
and that ${\Bbb S}_s(\A)$ is a quotient of $(\A^{\otimes s!})^{\oplus N}$.
Let $(\A_\bullet,\alpha_\bullet)$ be a weighted filtration of $\A$. Set $I:=\{\,1,\dots,t+1\,\}^{\times s!}$
and $\A_{t+1}:=\A$. For $(i_1,\dots,i_{s!})\in I$, define
$$
\A_{i_1}\cdot\cdots\cdot \A_{i_{s!}}
$$
as the image of the subsheaf $(\A_{i_1}\otimes\cdots\otimes \A_{i_{s!}})^{\oplus N}$ of $(\A^{\otimes s!})^{\oplus N}$
in ${\Bbb S}_s$ and
$$
\A_{i_1}\star\cdots\star \A_{i_{s!}}:=(\A_{i_1}\cdot\cdots\cdot \A_{i_{s!}})\cap {\Bbb V}_s(\A).
$$
\par
The \it standard weight vectors \rm are
$$
\gamma_r^{(i)}:=\bigl(\underbrace{i-r,\dots,i-r}_{i\times},\underbrace{i,\dots,i}_{(r-i)\times}\bigr),\q i=1,\dots,r-1.
$$
Given a weighted filtration $(\A_\bullet,{\alpha}_\bullet)$ of the torsion free sheaf $\A$, we obtain the
\it associated weight vector \rm
$$
\bigl(\underbrace{\gamma_{1},\ldots,\gamma_{1}}_{(\rk \A_1)\times}, \underbrace{\gamma_{2},\ldots,\gamma_{2}}_{(\rk
\A_2-\rk \A_1)\times},\ldots,\underbrace{\gamma_{t+1},\ldots,\gamma_{t+1}}_{(\rk \A-\rk \A_{t})\times}\bigr) :=\sum_{j=1}^t \alpha_j\cdot \gamma_r^{(\rk\A_j)}.
$$
(We recover $\alpha_j=(\gamma_{j+1}-\gamma_j)/r$, $j=1,\dots,t$.)
For a dispo sheaf $(\A,\phi)$ and a weighted filtration $(\A_\bullet,{\alpha}_\bullet)$ of $\A$, we finally find
with (\ref{eq1:NochEinKomp0})
\begin{equation}
\label{eq:AnotherWayToComputeMu}
\mu(\A_\bullet,{\alpha}_\bullet;\phi)=-\min\bigl\{\,\gamma_{i_1}+\cdots+\gamma_{i_{s!}}\,|\,
(i_1,\dots,i_{s!})\in I: \phi_{|\A_{i_1}\star \cdots\star\A_{i_{s!}}} \not\equiv 0\,\bigr\}.
\end{equation}
\par 
ii) We need a variant of the former definition. Let $(\A,\phi)$ be a dispo sheaf. We look at the representation $\kappa$ of $\GL_{r}(k)$ on ${\Bbb V}_s(k^r)$. Let $U$ be the maximal open subset on which $\A$ is locally free. Then, ${\Bbb V}_s(\A)_{|U}={\Bbb V}_s(\A_{|U})$ is the vector bundle that is associated to the vector bundle $\A_{|U}$ via the representation $\kappa$. Since $\kappa$ is clearly a polynomial representation, we can write ${\Bbb V}_s(k^r)$ as a quotient of ${\Bbb D}^{s!,v}(k^r)$ for an appropriate integer $v>0$. We let ${\Bbb D}^{s!,v}(\A_{|U})$ be the vector bundle with fiber ${\Bbb D}^{s!,v}(k^r)$ that is associated to $\A_{|U}$. By construction, we have a surjection ${\Bbb D}^{s!,v}(\A_{|U})\lra {\Bbb V}_s(\A_{|U})$, so that $\phi_{|U}$ induces a homomorphism
$$
\widetilde{\phi}\colon {\Bbb D}^{s!,v}(\A_{|U})\lra \Oh_U.
$$
Note that ${\Bbb D}^{s!,v}(k^r)$ is a subrepresentation of $({k^r}^{\otimes s!})^{\oplus N}$ for a suitable integer $N>0$. Hence, ${\Bbb D}^{s!,v}(\A_{|U})$ is a subbundle of $(\A_{|U}^{\otimes s!})^{\oplus N}$.
\par
Let $(\A_\bullet,\alpha_\bullet)$ be a weighted filtration of $\A$. As before, $I=\{\,1,\dots,t+1\,\}^{\times s!}$
and $\A_{t+1}=\A$. This time, we set
$$
\A_{i_1}\star\cdots\star \A_{i_{s!}}:=(\A_{i_1|U}\otimes \cdots\otimes \A_{i_{s!}|U})^{\oplus N}\cap {\Bbb D}^{s!,v}(\A_{|U}),\q (i_1,\dots,i_{s!})\in I.
$$
Then,
\begin{equation}
\label{eq:AnotherWayToComputeMu2}
\mu(\A_\bullet,{\alpha}_\bullet;\phi)=-\min\bigl\{\,\gamma_{i_1}+\cdots+\gamma_{i_{s!}}\,|\,
(i_1,\dots,i_{s!})\in I: \phi_{|\A_{i_1}\star \cdots\star\A_{i_{s!}}} \not\equiv 0\,\bigr\}.
\end{equation}
\end{Rem}
%%%%
Fix a positive polynomial $\delta\in\Q[x]$ of degree at most $\dim(X)-1$.
Now, we say that a dispo sheaf $(\A,\phi)$ is \it $\delta$-(semi)stable\rm, if the inequality
$$
M(\A_\bullet,{\alpha}_\bullet)+\delta\cdot \mu(\A_\bullet,{\alpha}_\bullet;\phi)(\succeq)0
$$
holds for every weighted filtration $(\A_\bullet,{\alpha}_\bullet)$ of $\A$.
\par
Let $\ol{\delta}$ be a non-negative rational number. We call a dispo sheaf $(\A,\phi)$
\it $\ol{\delta}$-slope (semi)stable\rm, if the inequality
$$
L(\A_\bullet,{\alpha}_\bullet)+\ol{\delta}\cdot \mu(\A_\bullet,{\alpha}_\bullet;\phi)(\ge)0
$$
holds for every weighted filtration $(\A_\bullet,{\alpha}_\bullet)$ of $\A$. Note that, for
$\delta=\ol{\delta}/(\dim(X)-1)!\cdot x^{\dim(X)-1}+\cdots$ (where $n=\dim X$), we have
%%%
\begin{equation}
\label{eq1:slopy}
(\A,\phi)\hbox{ is $\delta$-semistable}\q\Longrightarrow\q(\A,\phi)\hbox{ is $\ol{\delta}$-slope semistable}.
\end{equation}
%%%
\subsection{Global boundedness}
%%%
\begin{Thm}
\label{thm:DispShBound}
Fix a Hilbert polynomial $P$, a representation $\rho$, and an integer $s$ as above.
Then, the set of isomorphy classes of torsion free sheaves $\A$ on $X$ with Hilbert polynomial $P$ for which
there do exist a positive rational number $\ol\delta$ and a $\ol\delta$-slope semistable dispo sheaf $(\A,\phi)$ of type $(\rho,s)$ is bounded.
\end{Thm}
%%%
\begin{proof}
This is a slight modification of the proof of the corresponding result in \cite{GLSS2}. We will use the notation in Remark \ref{rem:DescSemStabAnotherWay}, ii).
\par 
Suppose $(\A,\phi)$ is a dispo sheaf which is $\ol{\delta}$-slope semistable for some 
$\ol{\delta}>0$. Assume $\A$ is not slope semistable as a sheaf and consider
its slope  Harder--Narasimhan filtration
$$
\A_\bullet: 0 = \A_0 \subsetneq \A_1 \subsetneq \A_2 \subsetneq \cdots \subsetneq \A_t \subsetneq 
\A_{t+1}=\A.
$$
We use the notation $\A^i=\A_i/\A_{i-1}$, $r_i:=\rk(\A_i)$, $r^i:=\rk(\A^i)$, and $\mu^i:=\mu(\A^i)$, 
$i=1,\dots,t+1$. Define
$$
C(\A_{\bullet})=\bigl\{\,{\gamma}=(\gamma _1, \dots ,
\gamma_{t+1})\in \R^{t+1}\,|\, \gamma_1\le \gamma_2\le \cdots
\le \gamma_{t+1}, \sum_{i=1}^{t+1}\gamma_i\cdot r^i=0\,\bigr\}.
$$
We equip $\R^{t+1}$ with the maximum norm $\|.\|$. For all ${\gamma}\in C(\A_{\bullet})\setminus \{0\}$, we have
$$
\sum_{i=1}^t \frac{\gamma_{i+1}-\gamma_i}{r}\cdot \bigl( r\cdot \deg(\A_i)-r_i\cdot \deg(\A) \bigr)< 0,
$$
so that the ${\delta}$-semistability of $(\A,\phi)$ implies
$$
f({\gamma}):=\mu\bigl(\A_{\bullet},{\alpha}_\bullet({\gamma});\phi\bigr)>0,\q
{\alpha}_\bullet({\gamma}):=\left(\,\frac{\gamma_2-\gamma_1}{r},\dots,\frac{\gamma_{t+1}-\gamma_t}{r}\,\right).
$$
Consider the set
$$
K:=C(\A_{\bullet})\cap \bigl\{\,{\gamma} \in \R^{t+1}\,|\,\|{\gamma}\|=1\bigr\}.
$$
Obviously $K$ is a compact set and $f$ is piecewise linear whence continuous, so that $f$ attains its infimum on $K$. It is easy to
see that there are only finitely many possibilities for the function $f$, so that we may bound this infimum from below by a
constant $C_0>0$ which depends only on the input data. 
\par
As usual, we let $U$ be the maximal open subset where $\A$ is locally free. We have the induced homomorphism $\widetilde{\phi}\colon {\Bbb D}^{s!,v}(\A_{|U})\lra \Oh_U$. Take a tuple $(i_1,\dots ,i_{s!})$ with $\widetilde{\phi}_{|\A_{i_1}\star\cdots\star \A_{i_{s!}}}\not\equiv 0$ which is minimal with respect to the lexicographic ordering of the index set $I$.
Define 
$$
\ol\A_{i_1,\dots,i_{s!}}
$$
as the quotient of $\A_{i_1}\star\cdots\star \A_{i_{s!}}$ by the subbundle that is generated by the $\A_{i'_1}\star\cdots\star \A_{i'_{s!}}$ for the index tuples $(i'_1,\dots,i'_{s!})$ which are strictly smaller than $(i_1,\dots ,i_{s!})$ in the lexicographic ordering.
By construction, $\widetilde{\phi}$ factorizes over a non-zero homomorphism
$$
\ol\phi\colon \ol\A_{i_1,\dots,i_{s!}}\lra \Oh_U,
$$
whence
$$
\mu_{\min}(\ol\A_{i_1,\dots,i_{s!}})\le 0.
$$
\par
In order to compute $\mu_{\min}(\ol\A_{i_1,\dots,i_{s!}})$, we observe that $\ol\A_{i_1,\dots,i_{s!}}$ is a subbundle of
$$
\bigl(\A^{i_1}_{|U}\otimes\cdots\otimes\A^{i_{s!}}_{|U}\bigr)^{\oplus N}.
$$
In fact, $\ol\A_{i_1,\dots,i_{s!}}$ is the vector bundle that is associated to the vector bundle
$$
\A_{|U}^{i_1}\oplus\cdots\oplus\A_{|U}^{i_t}
$$
by means of a representation $\rho_{i_1,\dots,i_{s!}}$ of $\GL_{r^1}(k)\times\cdots\times\GL_{r^t}(k)$ which is a subrepresentation of the representation on
$$
\bigl(k^{r^{i_1}}\otimes\cdots\otimes k^{r^{i_{s!}}}\bigr)^{\oplus N}.
$$
This already shows 
$$
\mu(\ol\A_{i_1,\dots,i_{s!}})=\mu^{i_1}+\cdots+\mu^{i_{s!}}.
$$
Now, we can apply Corollary \ref{cor:ExtofSG} to see that
$$ 
\mu^{i_1}+\cdots+\mu^{i_{s!}}\le \mu(\ol\A_{i_1,\dots,i_{s!}}) +D(\rho_{i_1,\dots,i_{s!}}).
$$
There are only finitely many possibilities for the representation $\rho_{i_1,\dots,i_{s!}}$, so that we may replace the constant $D(\rho_{i_1,\dots,i_{s!}})$ in the above inequality by a constant $C$ which depends only on $\rho$ and $s$.
\par
Altogether, we have demonstrated
%%%
\begin{equation}
\label{eq1:estimate}
\mu^{i_1}+\dots +\mu^{i_{s!}}\le C.
\end{equation}
%%%
Take the point
$$
{\gamma}:=\bigl(\mu(\A)-\mu^1, \dots , \mu(\A)-\mu^{t+1}\bigr)=\bigl(-\mu^1, \dots ,-\mu^{t+1}\bigr)
\in \R^{t+1}.
$$
By construction,
${\gamma} \in C(\A_{\bullet})\setminus\{0\}$ and
$$
f({\gamma})=\mu(\A_{\bullet}, \alpha_\bullet({\gamma});\phi)\le C.
$$
But $f$ is linear on each ray, so
$$
f({\gamma})=\|{\gamma}\| \cdot f\Bigl(\frac{\gamma}{\|{\gamma}\|}\Bigr)\ge C_0\cdot \|{\gamma}\|.$$
Now, this shows that either $\mu^1 =\|{\gamma}\|\le C^{\p}:={C}/{C_0}$
or $-\mu^{t+1}= \|{\gamma}\|\le C^{\p}$, i.e., 
$$
\hbox{either}\q\mu_{\max}(\A)\le \mu(\A)+C^{\p}
\q\hbox{or}\q\mu_{\min}(\A)\ge \mu(\A)-C^{\p}.
$$
The theorem finally follows from the boundedness theorem of Maruyama--Langer \cite{Lang}.
\end{proof}
%%%
\begin{Cor}
 \label{cor:AsymBound}
Fix the background data as in the theorem. There is a polynomial $\delta_\infty$, such that for every polynomial $\delta\succ \delta_\infty$ and every dispo sheaf $(\A,\phi)$ of type $(\rho,s)$ in which $\A$ has Hilbert polynomial $P$, the following conditions are equivalent:
\begin{itemize}
\item[{\rm i)}] $(\A,\phi)$ is $\delta$-(semi)stable.
\item[{\rm ii)}] For every weighted filtration $(\A_\bullet,\alpha_\bullet)$ of $\A$, one has
$$
\mu(\A_\bullet,\alpha_\bullet;\phi)\ge 0,
$$ 
 and
$$
M(\A_\bullet,\alpha_\bullet)(\succeq)0,
$$
for every weighted filtration $(\A_\bullet,\alpha_\bullet)$ with $\mu(\A_\bullet,\alpha_\bullet;\phi)=0$.
\end{itemize}
\end{Cor}
%%%
\begin{proof}
Let us call a dispo sheaf $(\A,\phi)$ which satisfies i) \it asymptotically (semi)stable\rm. Using \cite{GLSS2}, one can find a polynomial $\delta_0$, such that for every dispo sheaf $(\A,\phi)$ of type $(\rho,s)$ in which $\A$ has Hilbert polynomial $P$, the following holds true:
\begin{itemize}
 \item If $(\A,\phi)$ is asymptotically (semi)stable, then it is $\delta$-(semi)stable for every polynomial $\delta\succ \delta_0$.
\item Assume $\delta_0\prec \delta_1\prec \delta_2$. If $(\A,\phi)$ is $\delta_2$-(semi)stable, it is also $\delta_1$-(semi)stable.
\end{itemize}
%%%
Note also: If $(\A,\phi)$ is not asymptotically semistable, then it will not be $\delta$-semistable for any polynomial $\delta\succ\hskip -2pt\succ 0$. 
\par
What remains to show is that we can find $\delta_\infty$, such that, for every dispo sheaf $(\A,\phi)$ of type $(\rho,s)$ in which $\A$ has Hilbert polynomial $P$ and for every two polynomials $\delta_\infty\prec \delta_1\prec \delta_2$, the implication
$$
\hbox{$(\A,\phi)$ is $\delta_1$-semistable}\q\Longrightarrow\q\hbox{$(\A,\phi)$ $\delta_2$-semistable}
$$
is also correct. In \cite{GLSS2}, we referred to the instability flag for this. This is only adequate, if the characteristic of the base field is very large. We cannot assume this here. Therefore, we will give a different argument which relies only on general properties of semistability.
\par
As before, we will use the finite set ${\cal T}$ which depends only on $\rho$ and $s$, such that the condition of semistability of a dispo sheaf $(\A,\phi)$ has to be tested only for weighted filtrations $(\A_\bullet,\alpha_\bullet)$ of $\A$ with $((\rk\A_1,\dots,\rk\A_t),\alpha_\bullet)\in {\cal T}$. If $\delta_0\prec \delta_1\prec \delta_2$ and $(\A,\phi)$ is a dispo sheaf of type $(\rho,s)$ and Hilbert polynomial $P$ which is $\delta_1$-semistable but not $\delta_2$-semistable, there are a weighted filtration $(\A_\bullet,\alpha_\bullet)$ of $\A$ with $((\rk\A_1,\dots,\rk\A_t),\alpha_\bullet)\in {\cal T}$ and a polynomial $\delta_1\preceq \delta_\star\prec \delta_2$, such that
\begin{itemize}
 \item $M(\A_\bullet,\alpha_\bullet)+\delta_1\cdot \mu(\A_\bullet,\alpha_\bullet;\phi)\succeq 0$,  $M(\A_\bullet,\alpha_\bullet)+\delta_2\cdot \mu(\A_\bullet,\alpha_\bullet;\phi)\prec 0$, and
$M(\A_\bullet,\alpha_\bullet)+\delta_\star\cdot \mu(\A_\bullet,\alpha_\bullet;\phi) =0$. (Note that this implies $M(\A_\bullet,\alpha_\bullet)\succ 0$ and $\mu(\A_\bullet,\alpha_\bullet;\phi)\le 0$.)
\item $(\A,\phi)$ is $\delta_\star$-semistable.
\end{itemize}
There is the admissible deformation ${\df}_{(\A_\bullet,\alpha_\bullet)}(\A,\phi)=(\A_{\rm gr},\phi_{\rm gr})$.
It is performed with respect to the stability parameter $\delta_\star$, and $(\A_{\rm gr},\phi_{\rm gr})$ is still $\delta_\star$-semistable. We have 
$$
\A_{\rm gr}=\bigoplus_{i=1}^{t+1} \A^i\q\hbox{with}\q \A^i=\A_i/\A_{i-1},\q i=1,\dots,t+1.
$$
If we define $\A_{\rm gr,\bullet}$ via $\A_{\rm gr, i}:=\bigoplus_{j=1}^i \A^i$, $i=1,\dots,t$, it is clear that
$M(\A_{\rm gr,\bullet},\alpha_\bullet)=M(\A_\bullet,\alpha_\bullet)$ and $\mu(\A_{\rm gr,\bullet},\alpha_\bullet;\phi_{\rm gr})=\mu(\A_\bullet,\alpha_\bullet;\phi)$. 
Since $(\A_{\rm gr},\phi_{\rm gr})$ is $\delta_\star$-semistable, $\A_{\rm gr}$ belongs to a bounded family of torsion free sheaves with Hilbert polynomial $P$, by Theorem \ref{thm:DispShBound}. Moreover, ${\cal T}$ is finite, so that 
there are only finitely many possibilities for the polynomial $M(\A_{\bullet},\alpha_\bullet)$. There are evidently only finitely many choices for $\mu(\A_{\bullet},\alpha_\bullet;\phi_{\rm gr})$. The equation
$$
M(\A_{\bullet},\alpha_\bullet)+\delta_\star\cdot\mu(\A_{\bullet},\alpha_\bullet;\phi_{\rm gr})=0
$$
leaves therefore only finitely many options for $\delta_\star$. If we choose $\delta_\infty$ larger than the maximal possible value for $\delta_\star$, the assertion of the corollary will hold.
\end{proof}
%%%
\subsection{S-equivalence}
\label{subsec:SEquivDispo}
%%%
An important issue is the correct definition of S-equivalence of properly semistable dispo
sheaves. For this, suppose we are given a $\delta$-semistable dispo sheaf $(\A,\phi)$
and a weighted filtration $(\A_\bullet,{\alpha}_\bullet)$ of $\A$ with
$$
M(\A_\bullet,\alpha_\bullet)+\delta\cdot \mu(\A_\bullet,\alpha_\bullet;\phi)\equiv 0.
$$
We want to define the \it associated admissible deformation
${\rm df}_{(\A_\bullet,\alpha_\bullet)}(\A,\phi)=(\A_{\rm df},
\phi_{\rm df})$\rm. Of course, we set $\A_{\rm df}=\bigoplus_{i=0}^t\A_{i+1}/\A_i$.
Let $U$ be the maximal (big!) open subset where $\A_{\rm df}$ is locally free.
We may choose a one parameter subgroup
$\la\colon {\Bbb G}_m(k)\lra \SL_r(k)$ whose weighted flag $(W_\bullet(\la), {\alpha}_\bullet(\la))$ in $k^r$
satisfies:
%%%%
\begin{itemize}
\item $\dim_k(W_i)=\rk\A_i$, $i=1,\dots,t$, in $W_\bullet(\la): 0\subsetneq W_1\subsetneq\cdots\subsetneq W_t\subsetneq k^r$;
\item ${\alpha}_\bullet(\la)={\alpha}_\bullet$.
\end{itemize}
%%%%
Then, the given filtration $\A_\bullet$ corresponds to a reduction of the structure group of
${\cal I}som({\Oh}_U^{\oplus r},\A_{|U})$ to $Q(\la)$. On the other hand, $\la$ defines a decomposition
$$
{\Bbb V}_s= U^{\gamma_1}\oplus\cdots\oplus U^{\gamma_{u+1}},\q \gamma_1<\cdots<\gamma_{u+1}.
$$
Now, observe that $Q(\la)$ fixes the flag
\begin{equation}
\label{eq1:AFiltration}
0\subsetneq U_1:=U^{\gamma^1}\subsetneq U_2:=(U^{\gamma^1}\oplus U^{\gamma^2})\subsetneq\cdots\subsetneq
U_u:=(U^{\gamma_1}\oplus\cdots\oplus U^{\gamma_u})\subsetneq {\Bbb V}_s.
\end{equation}
Thus, we obtain a $Q(\la)$-module structure on
\begin{equation}
\label{eq1:L-ModuleIso}
\bigoplus_{i=0}^u U_{i+1}/U_i\cong {\Bbb V}_s.
\end{equation}
Next, we write $Q(\la)={\cal R}_u(Q(\la))\rtimes L(\la)$ where
$L(\la)\cong \GL(W_1/W_0)\times\cdots\times\GL(k^r/W_t)$
is the centralizer of $\la$. Note that (\ref{eq1:L-ModuleIso}) is an isomorphism of $L(\la)$-modules.
The process of passing from $\A$ to $\A_{\rm df}$ corresponds to first reducing the structure group
to $Q(\la)$, then extending it to $L(\la)$ via $Q(\la)\lra Q(\la)/{\cal R}_u(Q(\la))\cong L(\la)$, and
then extending it to $\GL_r(k)$ via the inclusion $L(\la)\subset \GL_r(k)$.
By (\ref{eq1:AFiltration}), ${\Bbb V}_s(\A_{|U})$ has a filtration
$$
0\subsetneq {\cal U}_1\subsetneq {\cal U}_2\subsetneq\cdots\subsetneq
{\cal U}_u\subsetneq {\Bbb V}_s(\A_{|U}),
$$
and, by (\ref{eq1:L-ModuleIso}), we have a canonical isomorphism
$$
{\Bbb V}_s(\A_{|U}) \cong \bigoplus_{i=1}^{u+1} {\cal U}_{i}/{\cal U}_{i-1}.
$$
Now, for $i_0$ with $\gamma_{i_0}=-\mu(\A_\bullet,\alpha_\bullet;\phi)$,
the restriction $\phi_{i_0}$ of $\phi_{|U}$ to ${\cal U}_{i_0}$ is non-trivial, and thus we may define
$\widetilde{\phi}_{\rm df}$ as the map induced by $\phi_{i_0}$ on ${\cal U}_{i_0}/{\cal U}_{i_0-1}$ and as zero on the other components.
Then, we finally obtain
$$
\phi_{\rm df}\colon {\Bbb V}_s(\A)\lra \iota_\star\bigl({\Bbb V}_s(\A_{|U})\bigr)
\stackrel{\iota_{\star}(\widetilde{\phi}_{\rm df})}{\lra}
\iota_{\star}({\Oh}_U)={\Oh}_X,
$$
$\iota\colon U\lra X$ being the inclusion.
A dispo sheaf $(\A,\phi)$ is said to be \it $\delta$-polystable\rm,
if it is $\delta$-semistable and isomorphic to every admissible deformation
${\rm df}_{(\A_\bullet,\alpha_\bullet)}(\A,\phi)=(\A_{\rm df},
\phi_{\rm df})$ associated to a filtration $(\A_\bullet,\alpha_\bullet)$ of $\A$ with
$$
M(\A_\bullet,\alpha_\bullet)+\delta\cdot \mu(\A_\bullet,\alpha_\bullet;\phi)\equiv 0.
$$
By the GIT construction of the moduli space which will be given in Section \ref{subsec:ConstModDispSh}, one has the following:
%%%
\begin{Lem}
Let $(\A,\phi)$ be a $\delta$-semistable dispo sheaf.
Then, there is a $\delta$-polystable admissible deformation ${\rm gr}(\A,\phi)$
of $(\A,\phi)$. The dispo sheaf ${\rm gr}(\A,\phi)$ is unique up to isomorphy.
\end{Lem}
%%%
In general, not every admissible deformation will immediately lead to a polystable dispo sheaf,
but any iteration of admissible deformations (leading to non-isomorphic dispo sheaves) will do so after finitely
many steps. We call two $\delta$-semistable dispo sheaves $(\A,\phi)$ and $(\A',\phi')$
\it S-equivalent\rm, if ${\rm gr}(\A,\phi)$ and ${\rm gr}(\A',\phi')$ are isomorphic.
%%%
\begin{Rem}
\label{rem1:AnotherWay}
Another way of looking at S-equivalence is the following: With the notation as above, we may choose
an open subset $U\subseteq X$ (no longer big), such that $\phi$ is surjective over $U$ and we have an isomorphism
$\psi\colon \A_{|U}\cong k^r\otimes {\Oh}_U$ with
$\psi(\A_i)= W_i\otimes {\Oh}_U$ for $i=1,\dots,t$. For such a trivialization, we obtain, from $\phi_{|U}$, the morphism
$$
\beta\colon U\lra {\Pe}\bigl({\Bbb V}_s(\A_{|U})\bigr)\stackrel{{\Bbb V}_s(\psi)}{\cong}
{\Pe}({\Bbb V}_s)\times U\lra {\Pe}({\Bbb V}_s).
$$
For the morphism $\beta_{\rm df}\colon U\lra {\Pe}({\Bbb V}_s)$ associated to $\phi_{{\rm df}|U}$,
we discover the relationship
\begin{equation}
\label{eq1:OrbitClosures}
\beta_{\rm df}(x)=\lim_{z\rightarrow\infty} \la(z)\cdot \beta(x),\q x\in U.
\end{equation}
\end{Rem}
%%%
\subsection{The main theorem on dispo sheaves}
%%%
With the definitions which we have encountered so far, we may introduce the moduli functors
\begin{eqnarray*}
\ul{\rm M}_P^{\delta\rm\hbox{-}(s)s}(\rho,s)\colon \ul{\rm Sch}_k &\lra& \ul{\rm Sets}
\\
S &\lma&\left\{\,
\begin{array}{l}
\hbox{Isomorphy classes of families of}
\\
\hbox{$\delta$-(semi)stable dispo sheaves $\A_S$ of}
\\
\hbox{type $(\rho,s)$ with Hilbert polynomial $P$}
\\
\hbox{parameterized by the scheme $S$}
\end{array}
\,\right\}
\end{eqnarray*}
%%%
for $\delta$-(semi)stable dispo sheaves $(\A,\phi)$ of type $(\rho,s)$ on $X$ with Hilbert polynomial $P(\A)=P.$
%%%
\begin{Thm}
\label{thm:ModDispShea}
Given the input data $P$, $\rho$, $s$, and $\delta$ as above, then the moduli space
${\rm M}_P^{\delta\rm\hbox{-}ss}(\rho,s)$ for $\delta$-semistable dispo sheaves $(\A,\tau)$ of type
$(\rho,s)$ with $P(\A)=P$ exists as a {\bfseries projective} scheme.
\end{Thm}
%%%
\subsection{The proof of the main theorem on dispo sheaves}
\label{subsec:ConstModDispSh}
%%%
In this section, we will outline how a GIT construction may be used for proving the main
auxiliary result Theorem \ref{thm:ModDispShea}. Once one has the correct set-up, the details
become mere applications of the techniques of the papers \cite{Schmitt0} or \cite{GS} and \cite{La3}.
%%%
\subsubsection{Construction of the parameter space}
%%%
As we have seen in Theorem \ref{thm:DispShBound}, there is a constant $C$, such that
$\mu_{\max}(\A)\le C$ for every ${\delta}$-semistable dispo sheaf $(\A,\phi)$ with
$P(\A)=P$, i.e., $\A$ lives in a bounded family. Thus, we may
choose an $n_0\gg 0$ with the following properties: For every sheaf $\A$ with Hilbert polynomial $P$ and
$\mu_{\max}(\A)\le C$ and every $n\ge n_0$, one has
%%%
\begin{itemize}
\item $H^i(\A(n))=\{0\}$ for $i>0$;
\item $\A(n)$ is globally generated.
\end{itemize}
%%%
We also fix a $k$-vector space $U$ of dimension $P(n)$. Let ${\frak Q}$ be the quasi-projective scheme which
parameterizes quotients $q\colon U\otimes{\Oh}_X(-n)\lra \A$ where $\A$ is a torsion free sheaf with Hilbert polynomial
$P$ and $H^0(q(n))$ an isomorphism. Let
$$
{\frak q}_{{\frak Q}}\colon U\otimes\pi_X^\star\bigl({\Oh}_X(-n)\bigr)\lra {\cal A}_{{\frak Q}}
$$
be the universal quotient. Setting
$$
{\Bbb V}_s(U):=\bigoplus_{(d_1,\dots,d_s):\atop d_i\ge 0, \sum i d_i=s!}
\Bigl({\rm Sym}^{d_1}\bigl((U\otimes_k V)^G\bigr)\otimes \cdots \otimes
{\rm Sym}^{d_s}\bigl({\rm Sym}^s(U\otimes_k V)^G\bigr)\Bigr),
$$
there is a homomorphism
$$
{\Bbb V}_s(U)\otimes  \pi_X^\star\bigl({\Oh}_{\frak Q\times X}(-s!\cdot n)\bigr)
\lra {\Bbb V}_s(\A_{{\frak Q}}),
$$
which is surjective over the open subset where $\A_{\frak Q}$ is locally free (see Lemma \ref{SurJect}).
For a point $[q\colon U\otimes{\Oh}_X(-n)\lra \A]\in {\frak Q}$, any
homomorphism $\phi\colon {\Bbb V}_s(\A) \lra{\Oh}_X$ is determined by the induced homomorphism
$$
{\Bbb V}_s(U)\lra H^0\bigl({\Oh}_X(s!\cdot n)\bigr)
$$
of vector spaces. Hence, our parameter space should be a subscheme of
$$
{\frak D}^0:={\frak Q} \times \underbrace{{\Pe}\Bigl({\rm Hom}\bigl({\Bbb V}_s(U), H^0({\Oh}_X(s!\cdot n))\bigr)^\vee\Bigr)}_{=:{\Pe}}.
$$
Note that, over ${\frak D}^0\times X$, there is the universal homomorphism
$$
\varphi''': {\Bbb V}_s(U)\otimes {\Oh}_{{\frak D}^0\times X}\lra H^0\bigl({\Oh}_X(s!\cdot n)\bigr)\otimes
\pi_\P^\star\bigl({\Oh}_{{\Pe}}(1)\bigr).
$$
Let $\varphi'' = {\rm ev} \circ {\varphi}'''$ be the composition of ${\phi}'''$ with the evaluation map
$$
{\rm ev}\colon H^0\bigl({\Oh}_X(s!\cdot n)\bigr)\otimes{\Oh}_{{\frak D}^0\times X}\lra
\pi_X^\star\bigl({\cal O}_X (s!\cdot n)\bigr).
$$
We twist $\varphi''$ by $\id_{\pi_X^\star({\Oh}_X(-s!\cdot n))}$ in order to obtain
$$
\varphi'\colon {\Bbb V}_s(U)\otimes \pi_X^\star\bigl({\Oh}_X(-s!\cdot n)\bigr) \lra \pi_\P^\star\bigl({\Oh}_{{\Pe}}(1)\bigr).
$$
Set $\A_{\frak D^0}:=\pi_{\frak Q\times X}^\star(\A_{\frak Q})$. We have the homomorphism
$S\colon{\Bbb V}_s(U)\otimes\pi_X^\star({\Oh}_X(-s!\cdot n))\lra {\Bbb V}_s({\cal A}_{\frak D^0})$.
Therefore, we may define a closed
subscheme ${\frak D}$ of ${\frak D}^0$ by the condition that $\phi'$ vanishes on $\ker(S)$.
Declaring $ {\cal A}_{{\frak D}}:=
({\cal A}_{{\frak D^0}})_{|{\frak D}\times X} $, there is thus the homomorphism
$$
{\phi}_{{\frak D}} \colon {\Bbb V}_s(\A_{\frak D}) \lra \pi_\P^\star\bigl({\Oh}_{{\Pe}}(1)\bigr)
$$
with $\phi_{|\frak D\times X}=\phi_{\frak D}\circ S$.
(To be precise, we first get $\phi_{\frak D}$ on the maximal open subset $V\subset {\frak D}\times X$
where $\A_{\frak D}$ is locally free and then extend it to ${\frak D}\times X$, using Corollary \ref{NatTrans}.
By the same token, $\phi_{|\frak D\times X}=\phi_{\frak D}\circ S$ is true, because it holds over $V$.)
The family
$(\A_{{\frak D}}, \phi_{{\frak D}})$ is the \it universal family of dispo sheaves parameterized by ${\frak D}$\rm.
By its construction, it has the features listed below.
%%%
\begin{Prop}[Local universal property]
\label{LUP0} Let $S$ be a scheme and $(\A_S,\phi_S)$ a family of
${\delta}$-semi\-stable dispo sheaves with Hilbert
polynomial $P$ parameterized by $S$. Then, there exist a covering
of $S$ by open subschemes $S_i$, $i\in I$, and morphisms
$\beta_i\colon S_i\lra {\frak D}$, $i\in I$, such that the family
$(\A_{S|S_i},\phi_{S|S_i})$ is isomorphic to the pullback of the
universal family on ${\frak D}\times X$ by $\beta_i\times\id_X$
for all  $i\in I$.
\end{Prop}
%%%
\subsubsection{The group action}
%%%
There is a natural action of $\GL(U)$ on the quot scheme ${\frak
Q}$ and on ${\frak D}^0$. This action leaves the closed subscheme
${\frak D}$ invariant, and therefore yields an action
$$
\Gamma\colon \GL(U)\times {\frak D}\lra {\frak D}.
$$
\begin{Prop}[Gluing property]
\label{Glue0} Let $S$ be a scheme and $\beta_i\colon S\lra {\frak
D}$, $i=1,2$, two morphisms, such that the pullback of the
universal family via $\beta_1\times \id_X$ is isomorphic to its
pullback via $\beta_2\times \id_X$. Then, there is a morphism
$\Xi\colon S\lra \GL(U)$, such that $\beta_2$ equals the morphism
$$
S\stackrel{\Xi\times\beta_1}{\lra}\GL(U)\times {\frak D}\stackrel{\Gamma}{\lra} \frak D.
$$
\end{Prop}
%%%
\subsubsection{Good quotients of the parameter space}
%%%
For a point $z\in {\frak D}$, we let $(\A_z,\phi_z)$ be the dispo sheaf obtained from the universal
family by restriction to $\{z\}\times X$. It will be our task to show that the set ${\frak D}^{\delta\hbox{-}\rm(s)s}$
parameterizing those points $z\in {\frak D}$ for which $(\A_z,\phi_z)$ is $\delta$-(semi)stable are open subsets of
${\frak D}$ which possess a good or geometric quotient. This can be most conveniently done by applying GIT.
To this end, we first have to exhibit suitable linearizations of the group action. We will use here the
approach by Gieseker in order to facilitate the computations. The experienced reader should have no problem
in rewriting the proof in Simpson's language.
\par
There is a projective subscheme ${\frak A}\hookrightarrow {\rm Pic}(X)$, such that the morphism
$\det\colon {\frak Q}\lra {\rm Pic}(X)$, $[q\colon U\otimes{\Oh}_X(-n)\lra \A]\lma[\det(\A)]$ factorizes over ${\frak A}$.
We choose a Poincar\'e sheaf ${\cal P}_{\frak A}$ on ${\frak A}\times X$.
Then, there is an integer $n_1$, such that for every integer $n\ge n_1$ and every line bundle $\L$ on $X$
with $[\L]\in {\frak A}$, the bundle $\L(rn)$ is globally generated and satisfies $h^i(\L(rn))=0$ for all $i>0$.
For such an $n$, the sheaf
$$
{\cal G}:=\pi_{\frak A\star}\Bigl({\cal P}_{\frak A}\otimes\pi_X^\star\bigl({\Oh}_X(rn)\bigr)\Bigr)
$$
is locally free. We then form the projective bundle
$$
{\Bbb G}_1:={\Pe}\Bigl({\cal H}om\bigl(\bigwedge^r U\otimes {\Oh}_{\frak A}, {\cal G}\bigr)^\vee\Bigr)
$$
over the scheme $\frak A$. For our purposes, we may always replace the Poincar\'e sheaf ${\cal P}_{\frak A}$
by its tensor product with the pullback of the dual of a
sufficiently ample line bundle on ${\frak A}$, so that we can achieve that ${\Oh}_{{\Bbb G}_1}(1)$ is ample.
The homomorphism
$$
\bigwedge^r \bigl(q_{\frak Q}\otimes {\id}_{\pi_X^\star({\Oh}_X(n))}\bigr)\colon \bigwedge^r U\otimes{\Oh}_{\frak Q\times X}
\lra \det(\A_{\frak Q})\otimes \pi_X^\star\bigl({\Oh}_X(rn)\bigr)
$$
defines a $\GL(U)$-equivariant and injective morphism
$$
{\frak Q}\lra {\Bbb G}_1.
$$
We declare
$$
{\Bbb G}_2:={\Pe}\Bigl({\rm Hom}\bigl({\Bbb V}_s(U), H^0({\Oh}_X(s!\cdot n))\bigr)^\vee\Bigr)\q\hbox{and}\q {\Bbb G}:={\Bbb G}_1\times {\Bbb G}_2.
$$
Then, we obtain the injective and $\SL(U)$-equivariant morphism
$$
{\rm Gies}\colon {\frak D}\lra {\Bbb G}.
$$
The ample line bundles ${\Oh}_{\Bbb G}(\nu_1,\nu_2)$, $\nu_1,\nu_2\in \Z_{>0}$, are naturally $\SL(U)$-linearized,
and we choose $\nu_1$ and $\nu_2$ in such a way that
\begin{equation}
\label{eq:GoodLins}
\frac{\nu_1}{\nu_2}=\frac{p-s!\cdot \delta(n)}{r\cdot\delta(n)}.
\end{equation}
%%%
\begin{Thm}
\label{thm:CharSemStab}
There exists $n_7\in \Z_{>0}$, such that for all $n\ge n_7$ the following property is verified:
For a point $z\in {\frak D}$, the Gieseker point ${\rm Gies}(z)\in {\Bbb G}$ is (semi)stable
with respect to the above linearization if and only if $(\A_z,\phi_z)$ is a $\delta$-(semi)stable
dispo sheaf of type $(\rho,s)$.
\end{Thm}
%%%
In the following, we will prove the theorem in several stages. As the first step,
we establish the following result.
%%%
\begin{Prop}
\label{prop:Bounded007}
There is an $n_2>0$, such that the following holds true: The set ${\frak S}$ of isomorphy classes
of torsion free sheaves $\A$ with Hilbert polynomial $P$ for which there exist an $n\ge n_2$ and
a point $z=([q\colon U\otimes{\Oh}_X(-n)\lra\A],\phi)\in {\frak D}$, such that
${\rm Gies}(z)$ is semistable with respect
to the above linearization, is bounded.
\end{Prop}
%%%
\begin{proof}
We would like to find a lower bound for $\mu_{\min}(\A)$ for a sheaf $\A$ as in the proposition.
Then, we may conclude with Theorem 4.2 of \cite{Lang}.
\par
Let ${\cal Q}=\A/{\cal B}$ be a torsion free quotient sheaf of $\A$. We have the exact sequence
$$
\begin{CD}
0 @>>> H^0\bigl({\cal B}(n)\bigr) @>>> H^0\bigl(\A(n)\bigr) @>>> H^0\bigl({\cal Q}(n)\bigr).
\end{CD}
$$
Let $\la\colon {\Bbb G}_m(k)\lra \SL(U)$ be a one parameter subgroup with weighted flag
$$
\bigl(U_\bullet(\la): 0\subsetneq U_1:=H^0\bigl(q(n)\bigr)^{-1}\bigl(H^0({\cal B}(n))\bigr) \subsetneq  U,\q
\alpha_\bullet(\la)=(1)\bigr).
$$
Define ${\cal B}':=q(U_1\otimes {\Oh}_X(-n))$. If ${\rm Gies}(z)=([M],[L])$, then
$$
\mu(\la,[M])= P(n)\cdot\rk({\cal B}')- h^0({\cal B}(n))\cdot r
\le P(n)\cdot\rk({\cal B})- h^0({\cal B}(n))\cdot r.
$$
Similarly to the proof of Theorem \ref{thm:DispShBound}, one
finds
$$
\mu(\la,[L])\le s!\cdot \bigl(P(n)-h^0({\cal B}(n))\bigr).
$$
The assumption that ${\rm Gies}(z)$ is semistable thus gives
\begin{eqnarray*}
0 &\le& \frac{\nu_1}{\nu_2}\cdot\mu(\la,[M])+\mu(\la,[L])
\\
&\le& \frac{P(n)-s!\cdot \delta(n)}{r\cdot\delta(n)}\cdot \bigl(P(n)\cdot\rk({\cal B})- h^0({\cal B}(n))\cdot r
\bigr)+ s!\cdot \bigl(P(n)-h^0({\cal B}(n))\bigr)
\\
&=&\frac{P(n)^2\cdot\rk({\cal B})}{r\cdot \delta(n)}-\frac{P(n)\cdot h^0({\cal B}(n))}{\delta(n)}-\frac{s!\cdot P(n)\cdot
\rk({\cal B})}{r}+s!\cdot P(n).
\end{eqnarray*}
%%%
We multiply this by $r\cdot \delta(n)/P(n)$ and find
$$
P(n) \rk({\cal B})-r h^0({\cal B}(n))+\delta(n) s!(r-1)
\ge
P(n) \rk({\cal B})-r h^0({\cal B}(n))+\delta(n) s!\bigl(r-\rk({\cal B})\bigr)\ge 0.
$$
The first exact sequence implies $h^0({\cal B}(n))\ge P(n)-h^0({\cal Q}(n))$. This enables us to transform
the above inequality into
\begin{equation}
\label{eq:LowerEsti}
\frac{h^0({\cal Q}(n))}{r}\ge \frac{P(n)}{r}-\frac{\delta(n)\cdot s!\cdot (r-1)}{\rk({\cal Q})\cdot r}
\ge \frac{P(n)}{r}-\frac{\delta(n)\cdot s!\cdot (r-1)}{r}.
\end{equation}
For a semistable sheaf $\E$ with $\mu(\E)\ge 0$, \cite{La3} provides the estimate
\begin{equation}
\label{eq:LangerEsti}
\frac{h^0(\E)}{\rk(\E)}\le \deg(X)\cdot {\frac{\mu(\E)}{\deg(X)}+f(r)+\dim(X)\choose \dim(X)}
\le \frac{\deg(X)}{\dim(X)!}\cdot \Bigl(\frac{\mu(\E)}{\deg(X)}+f(r)+\dim(X)\Bigr)^{\dim(X)}.
\end{equation}
%%%
If $\mu(\E)<0$, we have of course $h^0(\E)=0$.
The right-hand side $R(n)$ of (\ref{eq:LowerEsti}) is a positive polynomial of degree $\dim(X)$ with leading coefficient
$\deg(X)/\dim(X)!$. We can bound it from below by a polynomial of the form
$$
\frac{\deg(X)}{\dim(X)!}\cdot\Bigl(\frac{C}{\deg(X)}+f(r)+\dim(X)+n\Bigr)^{\dim(X)}.
$$
Assume that $n_2$ is so large that the value of this polynomial is positive and smaller than $R(n)$ for all
$n\ge n_2$. Then, (\ref{eq:LowerEsti}),
applied to the minimal destabilizing quotient ${\cal Q}$ of ${\cal A}$,
together with (\ref{eq:LangerEsti}) yields
$$
\mu_{\min}(\A)\ge C,
$$
and we are done.
\end{proof}
%%%
\begin{Thm}
\label{thm:Part1}
There is an $n_3$, such that for every $n\ge n_3$ and every point $z\in {\frak D}$ with
(semi)stable Gieseker point ${\rm Gies}(z)\in {\Bbb G}$, the dispo sheaf $(\A_z,\phi_z)$ is
$\delta$-(semi)stable.
\end{Thm}
%%%
\begin{proof}
As in \cite{SchmittGlobal}, Proposition 2.14, one may show that there is a finite
set
\begin{eqnarray*}
{\cal T}=\Bigl\{\, ({r}_\bullet^j,{\alpha}_\bullet^j)&\big|& {r}_\bullet^j=(r^j_1,\dots,r^j_{t_j}): 0<r^j_1<\cdots<r^j_{t_j}< r;
\\
&&{\alpha}^j_\bullet=(\alpha^j_1,\dots,\alpha^j_{t_j}): \alpha^j_i\in \Q_{>0},\ i=1,\dots,t_j,\ j=1,\dots,t\,\Bigr\},
\end{eqnarray*}
%%%
depending only on the $\GL_r(k)$-module ${\Bbb V}_s$,
such that the condition of $\delta$-(semi)stability of a dispo sheaf $(\A,\phi)$ of type $(\rho,s)$ with $P(\A)=P$
has to be verified only for weighted filtrations $(\A_\bullet,\alpha_\bullet)$ with
$$
\bigl((\rk(\A_1),\dots,\rk(\A_t)), \alpha_\bullet\bigr)\in {\cal T}.
$$
We may prescribe a constant $C'$. Then, there exists a constant $C''$, such that for every dispo sheaf $(\A,\phi)$
of type $(\rho,s)$ with $P(\A)=P$ and $[\A]\in {\frak S}$ and every weighted filtration $(\A_\bullet,\alpha_\bullet)$,
such that $((\rk(\A_1),\dots,\rk(\A_t)), \alpha_\bullet)\in {\cal T}$ and
\begin{equation}
\label{eq:FirstCase}
\mu(\A_i)\le C'',\q \hbox{for \bfseries one \rm index } i\in\{\,1,\dots,t\,\},
\end{equation}
one has
$$
L(\A_\bullet,\alpha_\bullet)>C'.
$$
It is easy to determine a constant $C'''$ which depends only on ${\Bbb V}_s$
with
$$
\mu(\A_\bullet,\alpha_\bullet;\phi)\ge -C'''
$$
for every weighted filtration $(\A_\bullet,\alpha_\bullet)$ of a sheaf $\A$ as above
with $((\rk(\A_1),\dots,\rk(\A_t)),\allowbreak \alpha_\bullet)\in {\cal T}$.
We choose $C'\ge \ol{\delta}\cdot C'''$. Then, for a dispo sheaf $(\A,\phi)$ of type $(\rho,s)$ with $[\A]\in {\frak S}$
and a weighted filtration $(\A_\bullet,\alpha_\bullet)$, such that $((\rk(\A_1),\dots,\rk(\A_t)), \alpha_\bullet)\in {\cal T}$
and (\ref{eq:FirstCase}) holds,
one has
$$
L(\A_\bullet,\alpha_\bullet)+\ol\delta\cdot \mu(\A_\bullet,\alpha_\bullet;\phi)> C'-\ol\delta\cdot C'''\ge 0,
$$
so that also
$$
M(\A_\bullet,\alpha_\bullet)+\delta\cdot \mu(\A_\bullet,\alpha_\bullet;\phi)\succ 0.
$$
Thus, the condition of $\delta$-(semi)stability has to be verified only for weighted filtrations
$(\A_\bullet,\alpha_\bullet)$ with $((\rk(\A_1),\dots,\rk(\A_t)), \alpha_\bullet)\in {\cal T}$ for which
(\ref{eq:FirstCase}) fails. But these live in bounded families. We conclude
%%%
\begin{Cor}
\label{cor:SimplifyStability}
There is a positive integer $n_4\ge n_3$, such that any $n\ge n_4$ has the following property: For every dispo sheaf
$(\A,\phi)$ of type $(\rho,s)$ for which $[\A]$ belongs to the bounded family ${\frak S}$, the conditions stated below
are equivalent.
\begin{itemize}
\item[{\rm 1.}] $(\A,\phi)$ is $\delta$-(semi)stable.
\item[{\rm 2.}] For every weighted filtration $(\A_\bullet,\alpha_\bullet)$ with
$((\rk(\A_1),\dots,\rk(\A_t)), \alpha_\bullet)\in {\cal T}$, such that $\A_j(n)$ is globally generated and
$h^i(\A_j(n))=0$ for all $i>0$, $j=1,\dots,t$, one has
$$
\sum_{j=1}^t\alpha_j\cdot\bigl(h^0(\A(n))\cdot\rk(\A_j)-h^0(\A_j(n))\cdot \rk(\A)\bigr)+
\delta(n)\cdot \mu(\A_\bullet,\alpha_\bullet;\phi)(\ge)0.
$$
\end{itemize}
\end{Cor}
%%%
We assume that $n\ge n_4$.
Now, let $z=([q\colon U\otimes{\Oh}_X(-n)\lra \A],\phi)\in {\frak D}$ be a point with (semi)stable Gieseker point
${\rm Gies}(z)$.
Then, $[\A]$ belongs to the bounded family ${\frak S}$. Therefore, it suffices to check Criterion 2.\ in Corollary
\ref{cor:SimplifyStability} for establishing the $\delta$-(semi)stability of $(\A,\phi)$.
\par
Let, more generally, $(\A_\bullet,\alpha_\bullet)$ be a weighted filtration of $\A$, such that
%%%
\begin{itemize}
\item $\A_j(n)$ is globally generated, $j=1,\dots,t$;
\item $h^i(\A_j(n))=0$, $i>0$, $j=1,\dots,t$.
\end{itemize}
%%%
Since $H^0(q(n))$ is an isomorphism, we define the subspaces
$$
U_j:=H^0\bigl(q(n)\bigr)^{-1}\bigl(H^0(\A_j(n))\bigr)\subsetneq U,\q j=1,\dots,t.
$$
Define the \it standard weight vectors \rm
$$
\gamma_p^{(i)}:=\bigl(\underbrace{i-p,\dots,i-p}_{i\times},\underbrace{i,\dots,i}_{(p-i)\times}\bigr),\q i=1,\dots,p-1,
$$
and choose a basis $\ul u=(u_1,\dots,u_p)$ of $U$, such that
$$
\langle\, u_1,\dots,u_{l_j}\,\rangle = U_j,\q l_j=\dim_k(U_j)=h^0(\A_j(n)),\ j=1,\dots,t.
$$
These data yield the weight vector
$$
\ul\gamma=(\gamma_1,\dots,\gamma_p):=\sum_{j=1}^t \alpha_j\cdot \gamma_p^{(l_j)}
$$
and the one parameter subgroup $\la:=\la(\ul u,\ul \gamma)\colon {\Bbb G}_m(k)\lra\SL(U)$
with
$$
\la(z)\cdot \sum_{i=1}^p c_i\cdot u_i=\sum_{i=1}^p z^{\gamma_i}\cdot c_i\cdot u_i,\q z\in {\Bbb G}_m(k).
$$
Similarly, we define the one parameter subgroups $\la^j:=\la(\ul u, \gamma_p^{(l_j)})$, $j=1,\dots,t$.
Let
$$
L\colon {\Bbb V}_s(U)\lra H^0\bigl({\Oh}_X(s!\cdot n)\bigr)
$$
be a linear map that represents the second component of ${\rm Gies}(z)$. We wish to compute
$\mu(\la,[L])$. First, we note that the choice of the basis $\ul u$ provides
an identification
$$
{\rm gr}(U):=\bigoplus_{j=1}^{t+1} U^j\cong U,\q U^j:=U_{j}/U_{j-1},\ j=1,\dots,t+1,
$$
which we will use without further mentioning in the following. Define $I:=\{\,1,\dots,t+1\,\}^{\times s!}$.
In analogy to the considerations at the very end of Section \ref{subsub:Spec},
we introduce the subspaces
$$
U^\star_{i_1,\dots,i_{s!}}\subset {\Bbb V}_s(U),\q (i_1,\dots,i_{s!})\in I.
$$
As before, we check that all weight spaces with respect to the one parameter subgroup $\la$ inside
${\Bbb V}_s(U)$ are direct sums of some of these subspaces. In addition, the subspaces
$U^{\star}_{i_1,\dots,i_{s!}}$ are eigenspaces for the one parameter subgroups $\la^1,\dots,\la^t$.
More precisely, $\la^j$ acts on $U^\star_{i_1,\dots,i_{s!}}$ with the weight
$$
s!\cdot l_j-\nu_j(i_1,\dots,i_{s!})\cdot p,\q (i_1,\dots,i_{s!})\in I,\ j=1,\dots,t.
$$
In that formula, we have used
$$
\nu_j(i_1,\dots,i_{s!}):=\#\bigl\{\, i_k\le j\,|\,k=1,\dots,s!\,\bigr\}.
$$
Thus, we find
%%%
\begin{equation}
\label{eq:LaterCompare}
\mu(\la,[L])=-\min\Bigl\{\,\sum_{j=1}^t\alpha_j\cdot\bigl(s!\cdot l_j-\nu_j(i_1,\dots,i_{s!})\cdot p\bigr)\,\big|\,
U_{i_1,\dots,i_{s!}}^\star\not\subseteq \ker(L),\ (i_1,\dots,i_{s!})\in I\,\Bigr\}.
\end{equation}
%%%
Fix an index tuple $(i_1^0,\dots,i_{s!}^0)\in I$ for which the minimum is achieved.
\par
Let
$$
M\colon \bigwedge^r U\lra H^0\bigl(\det(\A)(rn)\bigr)
$$
represent the first component of ${\rm Gies}(z)$. It is well known that
\begin{eqnarray*}
\mu(\la,[M])&=&\sum_{j=1}^t\alpha_j\cdot\bigl(h^0(\A(n))\cdot \rk(\A_j)-h^0(\A_j(n))\cdot \rk(\A)\bigr)
\\
&=&\sum_{j=1}^t\alpha_j\cdot\bigl(p\cdot \rk(\A_j)-h^0(\A_j(n))\cdot r\bigr).
\end{eqnarray*}
%%%
Since we assume ${\rm Gies}(z)$ to be (semi)stable, we have
%%%
\begin{eqnarray*}
0
&(\le)& \frac{\nu_1}{\nu_2}\cdot\mu(\la,[M]) + \mu(\la,[L])
\\
\\
&=& \frac{\nu_1}{\nu_2} \sum_{j=1}^t\alpha_j\bigl(p \rk(\A_j)-h^0(\A_j(n)) r\bigr)
- \sum_{j=1}^t\alpha_j\bigl(s! l_j-\nu_j(i^0_1,\dots,i^0_{s!}) p\bigr)
\\
\\
&=& \frac{p-s! \delta(n)}{r\delta(n)} \sum_{j=1}^t\alpha_j\bigl(p \rk(\A_j)-h^0(\A_j(n)) r\bigr)
-\sum_{j=1}^t\alpha_j\bigl(s! l_j-\nu_j(i^0_1,\dots,i^0_{s!}) p\bigr)
\\
\\
&{=}&\sum_{j=1}^t \alpha_j\Bigl(\frac{p^2\rk(\A_j)}{r \delta(n)}-\frac{p s! \rk(\A_j)}{r}
-\frac{p h^0(\A_j(n))}{\delta(n)}\Bigr) + \sum_{j=1}^t\alpha_j \nu_j(i_1^0,\dots,i_{s!}^0) p.
\end{eqnarray*}
%%%
For the last equation, we have used $l_j=h^0(\A_j(n))$.
We multiply this inequality by $r\cdot\delta(n)/p$. This leads to the inequality
$$
\sum_{j=1}^t \alpha_j\cdot\bigl(p\cdot \rk(\A_j)- h^0(\A_j(n))\cdot r\bigr)+\delta(n)\cdot\Bigl(- \sum_{j=1}^t
\alpha_j\cdot\bigl(s!\cdot \rk(\A_j)-\nu_j(i_1^0,\dots,i_{s!}^0)\cdot r\bigr)\Bigr)(\ge)0.
$$
To conclude, we have to verify
\begin{equation}
\label{eq:LastStep}
\mu(\A_\bullet,\alpha_\bullet;\phi)\ge -\sum_{j=1}^t
\alpha_j\cdot\bigl(s!\cdot \rk(\A_j)-\nu_j(i_1^0,\dots,i_{s!}^0)\cdot r\bigr).
\end{equation}
%%%
If $\ul\gamma$ is the weight vector with the distinct weights $\gamma_1<\cdots<\gamma_{t+1}$
which is associated to $(\A_\bullet,\alpha_\bullet)$ as in
Remark \ref{rem:DescSemStabAnotherWay}, then one easily checks that
$$
\gamma_{i^0_1}+\cdots+\gamma_{i^0_{s!}}=\sum_{j=1}^t
\alpha_j\cdot\bigl(s!\cdot \rk(\A_j)-\nu_j(i_1^0,\dots,i_{s!}^0)\cdot r\bigr).
$$
In view of (\ref{eq:AnotherWayToComputeMu}), it remains to show that
\begin{equation}
\label{eq:FinishArgument}
\phi_{|\A_{i_1^0}\star\cdots\star \A_{i_{s!}^0}}\not\equiv 0.
\end{equation}
To this end, note that, up to a scalar, $L$ is given as
$$
\Bigl(\iota\colon{\Bbb V}_s(U)\lra H^0\bigl({\Bbb V}_s(\A)(s!\cdot n)\bigr)\Bigr)\stackrel{
H^0(\phi\otimes{\id}_{{\Oh}_X(s!\cdot n)})}{\lra}
H^0\bigl({\Oh}_X(s!\cdot n)\bigr).
$$
The image of $U^\star_{i_1^0,\dots,i^0_{s!}}$ lies in the subspace
$H^0((\A_{i_1^0}\star\cdots\star \A_{i_{s!}^0})(s!\cdot n))$ and that shows (\ref{eq:FinishArgument}).
\end{proof}
%%%
We now turn to the converse direction in the proof of Theorem \ref{thm:CharSemStab}. Again, we need a preparatory
result.
%%%
\begin{Prop}
\label{prop:SectSemStab}
There is a positive integer $n_5$, such that any $\delta$-(semi)stable dispo sheaf $(\A,\phi)$ of type $(\rho,s)$
with $P(\A)=P$ satisfies
$$
\sum_{j=1}^t \alpha_j\cdot\bigl(P(n)\cdot \rk(\A_j)-h^0(\A_j(n))\cdot r\bigr)+\delta(n)\cdot \mu(\A_\bullet,\alpha_\bullet;
\phi)(\ge) 0
$$
for every weighted filtration $(\A_\bullet,\alpha_\bullet)$ of $\A$ and every $n\ge n_5$.
\end{Prop}
%%%
\begin{proof}
By Theorem \ref{thm:DispShBound}, we know that there is a bounded family ${\frak S}'$ of torsion
free sheaves with Hilbert polynomial $P$, such that $[\A]\in {\frak S}'$ for any $\delta$-semistable dispo
sheaf $(\A,\phi)$ of type $(\rho,s)$ with $P(\A)=P$. We choose a constant $C$, such that $\mu_{\max}(\A)\le C$
for every torsion free sheaf $\A$ on $X$ with $[\A]\in {\frak S}'$. Given an additional positive constant
$C'$, we subdivide the class of torsion free sheaves ${\cal B}$ which might occur as saturated subsheaves of
${\Oh}_X$-modules $\A$ with $[\A]\in {\frak S}$ into two classes:
\begin{itemize}
\item[A.] $\mu({\cal B})\ge -C'$.
\item[B.] $\mu({\cal B})<-C'$.
\end{itemize}
By a lemma of Grothendieck's (\cite{HL}, Lemma 1.7.9), the sheaves ${\cal B}$ falling into class A live again
in bounded families, so that we may always assume that our $n$ is large enough, such that any such sheaf
${\cal B}$ satisfies $h^i({\cal B}(n))=0$, $i>0$.
\par
If $\E$ is a torsion free sheaf on $X$ with Harder--Narasimhan
filtration $0=:\E_0\subsetneq \E_1\subsetneq\cdots\subsetneq \E_t\subsetneq \E_{t+1}:=\E$, then
$$
h^0(\E)\le \sum_{i=1}^{t+1}h^0(\E_i/\E_{i-1}),
$$
so that (\ref{eq:LangerEsti}) gives, with $F(r):=\max\{\, f(i)\,|\, i=1,\dots,r\,\}$ and $\rk(\E)\le r$,
%%%
\begin{eqnarray*}
h^0(\E) &\le&\bigl(\rk(\E)-1\bigr)\cdot \frac{\deg(X)}{\dim(X)!}\cdot \Bigl(\frac{\mu_{\max}(\E)}{\deg(X)}+F(r)+\dim(X)\Bigr)^{\dim(X)}+
\\
&&+ \frac{\deg(X)}{\dim(X)!}\cdot \Bigl(\frac{\mu(\E)}{\deg(X)}+F(r)+\dim(X)\Bigr)^{\dim(X)}.
\end{eqnarray*}
%%%
For any sheaf $\A$ with $[\A]\in {\frak S}'$ and any saturated subsheaf ${\cal B}$ of $\A$ which belongs to the
class B, we thus find
%%%
\begin{eqnarray}
\nonumber
h^0({\cal B}(n))&\le&  \bigl(\rk({\cal B})-1\bigr)\cdot\frac{\deg(X)}{\dim(X)!}\cdot \Bigl(\frac{C}{\deg(X)}+F(r-1)+\dim(X)+n\Bigr)^{\dim(X)}
\\
\label{eq:SectEsti}
&&+ \frac{\deg(X)}{\dim(X)!}\cdot \Bigl(\frac{-C'}{\deg(X)}+F(r-1)+\dim(X)+n\Bigr)^{\dim(X)}
\\
\nonumber
&=:&R\bigl(\rk({\cal B}),C'\bigr)(n).
\end{eqnarray}
%%%
We choose $C'$ so large that
\begin{equation}
\label{eq:PolySuffLarge}
P\cdot \rk({\cal B})- h^0({\cal B}(n))\cdot r\ge P\cdot \rk({\cal B})- R\bigl(\rk({\cal B}),C'\bigr)\cdot r
= K\cdot x^{\dim(X)-1}+\cdots \succ \delta\cdot s!\cdot (r-1).
\end{equation}
%%%
For all $n\gg 0$, (\ref{eq:PolySuffLarge}) remains true when evaluated at $n$.
\par
Now, let $(\A_\bullet,\alpha_\bullet)$ be any weighted filtration of $\A$.
Write $\{\,1,\dots,t\,\}=I_{\rm A}\sqcup I_{\rm B}$ with $i\in I_{\rm A/B}$ if and only if $\A_i$ belongs to the
class A/B. Let $1\le i^{\rm A/B}_1< \cdots<i^{\rm A/B}_{t_{\rm A/B}}\le t$ be the indices in $I_{\rm A/B}$ and define the weighted
filtrations $(\A^{\rm A/B}_\bullet,\alpha^{\rm A/B}_\bullet)$ with
\begin{eqnarray*}
\A^{\rm A/B}_\bullet&:&0\subsetneq \A_1^{\rm A/B}:=\A_{i^{\rm A/B}_1}\subsetneq \cdots\subsetneq
\A_{t_{\rm A/B}}^{\rm A/B}:=\A_{i^{\rm A/B}_{t_{\rm A/B}}}\subsetneq \A,
\\
\alpha^{\rm A/B}_\bullet&=&(\alpha^{\rm A/B}_1,\dots,\alpha^{\rm A/B}_{t_{\rm A/B}})
:=(\alpha_{i^{\rm A/B}_1},\dots,\alpha_{i^{\rm A/B}_{t_{\rm A/B}}}).
\end{eqnarray*}
It is easy to see that
\begin{equation}
\label{eq:Split}
\mu(\A_\bullet,\alpha_\bullet;\phi)\ge\mu(\A^{\rm A}_\bullet,\alpha^{\rm A}_\bullet;\phi)
-s!\cdot(r-1)\cdot \sum_{j=1}^{t_{\rm B}} \alpha^{\rm B}_{j}.
\end{equation}
Now, we compute
\begin{eqnarray*}
&&\sum_{j=1}^t \alpha_j\cdot\bigl(P(n)\cdot \rk(\A_j)-h^0(\A_j(n))\cdot r\bigr)+\delta(n)\cdot \mu(\A_\bullet,\alpha_\bullet;
\phi)
\\
&\stackrel{\rm(\ref{eq:Split})}{\ge}&
\sum_{j=1}^{t_{\rm A}} \alpha_j^{\rm A}\cdot \bigl(P(n)\cdot \rk(\A^{\rm A}_j)-h^0(\A^{\rm A}_j(n))\cdot r\bigr)
+\delta(n)\cdot \mu(\A^{\rm A}_\bullet,\alpha^{\rm A}_\bullet;\phi)+
\\
&&+\sum_{j=1}^{t_{\rm B}} \alpha_j^{\rm B}\cdot \bigl(P(n)\cdot \rk(\A^{\rm B}_j)-h^0(\A^{\rm B}_j(n))\cdot r\bigr)
-\delta(n)\cdot s!\cdot(r-1)\cdot \sum_{j=1}^{t_{\rm B}} \alpha^{\rm B}_{j}
\\
&\stackrel{n\gg 0}{=}&
M(\A^{\rm A}_\bullet,\alpha^{\rm A}_\bullet)(n)+\delta(n)\cdot
\mu(\A^{\rm A}_\bullet,\alpha^{\rm A}_\bullet;\phi)+
\\
&&+\sum_{j=1}^{t_{\rm B}} \alpha_j^{\rm B}\cdot \Bigl(\bigl(P(n)\cdot \rk(\A^{\rm B}_j)-h^0(\A^{\rm B}_j(n))\cdot r\bigr)
-\delta(n)\cdot s!\cdot(r-1)\Bigr)
\\
&\stackrel{{\rm (\ref{eq:PolySuffLarge})\&} n\gg 0}{\ge}&
M(\A^{\rm A}_\bullet,\alpha^{\rm A}_\bullet)(n)+\delta(n)\cdot
\mu(\A^{\rm A}_\bullet,\alpha^{\rm A}_\bullet;\phi)
\q\stackrel{n\gg 0}{(\ge)}\q 0.
\end{eqnarray*}
%%%
The last estimate results from the condition of $\delta$-(semi)stability, applied to the weighted filtration
$(\A^{\rm A}_\bullet,\alpha^{\rm A}_\bullet)$. We still have to justify that, in this last estimate, $n$ can be uniformly chosen
for all polynomials of the form $M(\A_\bullet,\alpha_\bullet)+\delta\cdot \mu(\A_\bullet,\alpha_\bullet;\phi)$
where all the members of the filtration $\A_\bullet$ belong to the class A. We use again the set ${\cal T}$ which has been introduced
in the proof of Theorem \ref{thm:Part1}. Then, any polynomial of the above form can be written as a positive
rational linear combination of polynomials
$M(\A^i_\bullet,\alpha^i_\bullet)+\delta\cdot \mu(\A^i_\bullet,\alpha^i_\bullet;\phi)$ where
$((\rk(\A_1^i),\dots,\rk(\A^i_{t_i}),\alpha_\bullet^i)\in {\cal T}$ and all the members of the filtration $\A^i_\bullet$
belong to the class A, $i=1,\dots,u$.
By the boundedness of the set of isomorphy classes of sheaves in the class A, it now follows that there
are only finitely many polynomials of the form
$M(\A_\bullet,\alpha_\bullet)+\delta\cdot \mu(\A_\bullet,\alpha_\bullet;\phi)$ where all the members of $\A_\bullet$
belong to the class A and $((\rk(\A_1),\dots,\rk(\A_{t}),\alpha_\bullet)\in {\cal T}$. This proves our last claim and the proposition.
\end{proof}
%%%
\begin{Thm}
\label{thm:Part2}
There exists a positive integer $n_6$, enjoying the following property: If $n\ge n_6$ and
$(\A,\phi)$ is a $\delta$-(semi)stable dispo sheaf of type $(\rho,s)$ with $P(\A)=P$, then, for a point
$z\in {\frak D}$ of the form $z=([q\colon U\otimes{\Oh}_X(-n)\lra\A], \phi)$, the associated Gieseker
point ${\rm Gies}(z)$ is (semi)stable for the given linearization.
\end{Thm}
%%%%
\begin{proof}
Let $\la\colon {\Bbb G}_m(k)\lra \SL(U)$ be a one parameter subgroup and suppose ${\rm Gies}(z)=([M],[L])$. Then, we
have to verify that
$$
\frac{\nu_1}{\nu_2}\cdot \mu(\la,[M])+\mu(\la,[L]) (\ge)0.
$$
The one parameter subgroup $\la$ provides the weighted flag $(U_\bullet(\la),\beta_\bullet(\la))$ with
$$
U_\bullet(\la): 0=:U_0\subsetneq U_1\subsetneq\cdots\subsetneq U_\tau\subsetneq U_{\tau+1}:=U;\q
\beta_\bullet(\la)=(\beta_1,\dots,\beta_\tau).
$$
For each $h\in\{\,1,\dots,\tau\,\}$, we let $\widetilde{\A}_h$ be the saturated subsheaf that is generically generated
by $q(U_h\otimes{\Oh}_X(-n))$. There may be improper inclusions among the $\widetilde{\A}_h$'s. After clearing
these, we obtain the filtration
$$
\A_\bullet: 0=:\A_0\subsetneq \A_1\subsetneq \cdots\subsetneq \A_t\subsetneq \A_{t+1}:=\A.
$$
For $j=1,\dots,t$, we define
$$
T(j):=\bigl\{\,h\in\{\,1,\dots,\tau\,\}\,|\, \widetilde{\A}_{h}=\A_j\,\bigr\}
$$
and
$$
\alpha_j:=\sum_{h\in T(j)}\beta_h.
$$
This gives the weighted filtration $(\A_\bullet,\alpha_\bullet)$ of $\A$. By Proposition \ref{prop:SectSemStab},
\begin{equation}
\label{eq:Sect1}
\sum_{j=1}^t \alpha_j\cdot\bigl(P(n)\cdot \rk(\A_j)-h^0(\A_j(n))\cdot r\bigr)+\delta(n)\cdot \mu(\A_\bullet,\alpha_\bullet;
\phi)(\ge) 0.
\end{equation}
Recall from (\ref{eq:AnotherWayToComputeMu}) that
\begin{equation}
\label{eq:Remindi}
\mu(\A_\bullet,{\alpha}_\bullet;\phi)=-\min\bigl\{\,\gamma_{i_1}+\cdots+\gamma_{i_{s!}}\,|\,
(i_1,\dots,i_{s!})\in I: \phi_{|\A_{i_1}\star \cdots\star\A_{i_{s!}}} \not\equiv 0\,\bigr\}.
\end{equation}
Let $(i_1^0,\dots,i_{s!}^0)\in I=\{\, 1,\dots,t+1\,\}^{\times s!}$ be an index tuple which computes the minimum.
With
$$
\nu_j(i_1,\dots,i_{s!}):=\#\bigl\{\, i_k\le j\,|\,k=1,\dots,s!\,\bigr\},
$$
one calculates
$$
\gamma_{i^0_1}+\cdots+\gamma_{i^0_{s!}}
=\sum_{j=1}^t\alpha_j\cdot\bigl(s!\cdot \rk(\A_j)-\nu_j(i^0_1,\dots,i^0_{s!})\cdot r\bigr).
$$
Thus, (\ref{eq:Sect1}) transforms into
$$
\sum_{j=1}^t \alpha_j\cdot\bigl(P(n)\cdot \rk(\A_j)-h^0(\A_j(n))\cdot r\bigr)+\delta(n)\cdot
\Bigl(-\sum_{j=1}^t\alpha_j\cdot\bigl(s!\cdot \rk(\A_j)-\nu_j(i^0_1,\dots,i^0_{s!})\cdot r\bigr)\Bigr)
(\ge) 0.
$$
A computation as in the proof of Theorem \ref{thm:Part1}, but performed backwards, shows that this implies
\begin{equation}
\label{eq:Sect2}
\frac{\nu_1}{\nu_2}\cdot \Bigl(\sum_{j=1}^t \alpha_j\cdot\bigl(P(n)\cdot \rk(\A_j)-h^0(\A_j(n))\cdot r\bigr)
\Bigr)-\sum_{j=1}^t\alpha_j\cdot\bigl(s!\cdot h^0(\A_j(n))-\nu_j(i^0_1,\dots,i^0_{s!})\cdot p\bigr)
(\ge) 0.
\end{equation}
First, we see that
\begin{equation}
\label{eq:Estimate0}
\mu(\la,[M])=\sum_{h=1}^\tau \beta_h\cdot\bigl(P(n)\cdot \rk(\widetilde{\A}_h)-\dim_k(U_h)\cdot r\bigr)
\ge \sum_{j=1}^t \alpha_j\cdot\bigl(P(n)\cdot \rk(\A_j)-h^0(\A_j(n))\cdot r\bigr).
\end{equation}
%%%
We need a little more notation. For $j=0,\dots,t+1$, we introduce
\begin{eqnarray*}
\ul h(j)&:=&\min\bigl\{\,h=1,\dots,\tau\,|\, \widetilde{\A}_h=\A_j\,\bigr\};\q \ul U_j:=U_{\ul h(j)}
\\
\ol h(j)&:=&\max\bigl\{\,h=1,\dots,\tau\,|\, \widetilde{\A}_h=\A_j\,\bigr\};\q \ol U_j:=U_{\ol h(j)},
\end{eqnarray*}
%%%
as well as
$$
\widetilde{U}_j:=\ul U_{j}/\ol U_{j-1},\q j=1,\dots,t+1.
$$
For an index tuple $(i_1,\dots,i_{s!})\in I$, we find the vector space
$$
\widetilde{U}_{i_1,\dots,i_{s!}}:= \bigl(\widetilde{U}_{i_1}\otimes \cdots \otimes\widetilde{U}_{i_{s!}}\bigr)^{\oplus N}.
$$
Using a basis $\ul u$ of $U$ which consists of eigenvectors for the one parameter subgroup $\la$, we identify these
spaces with subspaces of $(U^{\otimes s!})^{\oplus N}$. Then,
$\widetilde{U}_{i_1,\dots,i_{s!}}^\bullet$ stands for the image of $\widetilde{U}_{i_1,\dots,i_{s!}}$ in
$$
\bigoplus_{(d_1,\dots,d_s):\atop d_i\ge 0, \sum i d_i=s!}
\Bigl({\rm Sym}^{d_1}\bigl(U\otimes_k V\bigr)\otimes \cdots \otimes
{\rm Sym}^{d_s}\bigl({\rm Sym}^s(U\otimes_k V)\bigr)\Bigr)
$$
and $\widetilde{U}_{i_1,\dots,i_{s!}}^\star$ for the intersection of $\widetilde{U}_{i_1,\dots,i_{s!}}^\bullet$
with ${\Bbb V}_s(U)$.
A similar construction, generalizing the one at the end of Section \ref{SomeGIT}, associates to a collection of
subspaces $Y_{1},\dots,Y_{{s!}}$ of $U$ the subspace
$Y_{1}\star\cdots\star Y_{s!}\ \hbox{of}\ (U^{\otimes s!})^{\oplus N}$.
Note that
$$
\la=\la(\ul u,\ul \gamma)\q\hbox{with}\q \ul\gamma=\sum_{h=1}^\tau \beta_h\cdot\gamma_p^{(\dim_k(U_h))}.
$$
We define $\la^h:=\la(\ul u, \gamma_p^{(\dim_k(U_h))})$, $h=1,\dots,\tau$. The effect of our definition is that
the spaces $\widetilde{U}_{i_1,\dots,i_{s!}}^\star$, $(i_1,\dots,i_{s!})\in I$, are weight spaces for $\la$ as well
as for $\la^1,\dots,\la^\tau$. We associate to an index $h\in\{\,1,\dots,\tau\,\}$ the index $j(h)\in \{\,1,\dots,t\,\}$
with $\widetilde{\A}_h=\A_{j(h)}$. Then, $\ul{h}(j)\le h$ holds if and only if $j\le j(h)$, and one verifies that $\la$ acts on $\widetilde{U}_{i^0_1,\dots,i^0_{s!}}^\star$
with the weight
%%%%
\begin{equation}
\label{eq:Estimate1}
-\sum_{h=1}^\tau\beta_h\cdot\bigl(s!\cdot \dim_k(U_h)-\nu_{j(h)}(i^0_1,\dots,i^0_{s!})\cdot p\bigr)
\ge
-\sum_{j=1}^t\alpha_j\cdot\bigl(s!\cdot h^0(\A_j(n))-\nu_j(i^0_1,\dots,i^0_{s!})\cdot p\bigr).
\end{equation}
%%%%
In view of the estimates (\ref{eq:Sect2}), (\ref{eq:Estimate0}), and (\ref{eq:Estimate1}), it is now sufficient
to ascertain that the restriction of $L$ to $\widetilde{U}_{i^0_1,\dots,i^0_{s!}}^\star$ is non-trivial.
If it were trivial, then there would have to be an index tuple $(i_1',\dots,i_{s!}')$ with $i_l'\le i^0_l$, $l=1,\dots,s!$,
at least one inequality being strict, such that
\begin{equation}
\label{eq:NonZero}
L_{|\ul U_{i_1'}\star\cdots\star \ul U_{i_{s!}'}}\not\equiv 0.
\end{equation}
This is because $L$ restricts to a non-zero map on $\ul U_{i_1^0}\star\cdots\star\ul U_{i_{s!}^0}$,
as $\phi_{|\A_{i_1^0}\star\cdots\star\A_{i^0_{s!}}}$ is non-trivial.
Now, if (\ref{eq:NonZero}) holds true, then we must also have
\begin{equation}
\phi_{|\A_{i_1'}\star\cdots\star \A_{i_{s!}'}}\not\equiv 0.
\end{equation}
(Compare the arguments at the end of the proof of Theorem \ref{thm:Part1}.)
But, then the tuple $(i_1^0,\dots,i_{s!}^0)$ would not give the minimum in (\ref{eq:Remindi}), a contradiction.
\end{proof}
%%%
By Theorem \ref{thm:CharSemStab}, the subsets ${\frak D}^{\delta\hbox{-}\rm (s)s}$ of $\delta$-(semi)stable dispo sheaves
are the preimages of the sets of GIT-(semi)stable points in ${\Bbb G}$ under the Gieseker morphism.
Therefore, they are open subsets.
%%%
\begin{Prop}
\label{prop:Proper}
The restricted Gieseker morphism
$$
{\rm Gies}_{|\frak D^{\delta\hbox{-}\rm ss}}\colon {\frak D}^{\delta\hbox{-}\rm ss}\lra {\Bbb G}^{\rm ss}
$$
is proper. Since it is also injective, it is finite.
\end{Prop}
%%%
\begin{proof}
This is pretty standard, so we can be a bit sketchy. We apply the valuative criterion of properness.
Let $\ol{\frak Q}$ be the closure of $\frak Q$ in the
quot scheme of $U\otimes{\Oh}_X(-n)$. Then, the parameter space $\frak D$ may also be compactified to
$\ol{\frak D}\lra\ol{\frak Q}$. Given a discrete valuation ring $R$, a morphism $\eta\colon C:={\rm Spec}(R)\lra
{\Bbb G}^{\rm ss}$ which lifts over $C^\star:={\rm Spec}({\rm Quot}(R))$ to a morphism $\eta^\star
\colon C^\star\lra \frak D^{\delta\hbox{-}\rm ss}$, we may first
extend $\eta^\star$ to a morphism $\ol\eta\colon C\lra \ol{\frak D}$. This morphism is associated to
a family
$$
\bigl(q_C\colon U\otimes\pi_X^\star({\Oh}_X(-n))\lra \A_C,\phi_C\colon {\Bbb V}_s(\A_C)\lra {\Oh}_{C\times X}\bigr)
$$
on $C\times X$ where the restriction of $\A_C$ to the special fiber $\{0\}\times X$ may have torsion.
(As usual, one gets $\phi_C$ first on the open subset where $\A_C$ is locally free and then extends it to $X$.)
Let $Z$ be the support of that torsion. Then, the family $q_C$ may be altered to a family
$\widetilde{q}_C\colon U\otimes \pi_X^\star({\Oh}_X(-n))\lra \widetilde{\A}_C$ where $\widetilde{\A}_C$ is now
a $C$-flat family of torsion free sheaves, $\widetilde{q}_C$ agrees with $q_{C}$ on $(C\times X)\setminus Z$,
but $\widetilde{q}_{C|\{0\}\times X}$ may fail to be surjective in points of $Z$.
Let $\iota\colon (C\times X)\setminus Z\lra C\times X$ be the inclusion. Define $\widetilde{\phi}_C$ as the composition
$$
{\Bbb V}_s(\widetilde{\A}_C)\lra \iota_\star\bigl(
{\Bbb V}_s(\widetilde{\A}_{C|(C\times X)\setminus Z})\bigr)=
\iota_\star\bigl({\Bbb V}_s({\A}_{C|(C\times X)\setminus Z})\bigr)
\stackrel{\iota_\star(\phi_{C|(C\times X)\setminus Z})}{\lra} \iota_\star\bigl({\Oh}_{(C\times X)\setminus Z}\bigr)={\Oh}_{C\times X}.
$$
The family $(\widetilde{q}_C,\widetilde{\phi}_C)$ also defines a morphism to ${\Bbb G}$ which coincides with
$\eta$. Let $(\widetilde{q}\colon U\otimes {\Oh}_X(-n)\lra \widetilde{\A},\widetilde{\phi})$ be the restriction of the
new family to $\{0\}\times X$.
One checks the following results:
\begin{itemize}
\item $H^0(\widetilde{q}(n))$ must be injective;
\item Since $(\widetilde{q},\widetilde{\phi})$ defines a semistable point, $\widetilde{\A}$ belongs to a bounded
family (this is an easy adaptation of the proof of Proposition \ref{prop:Bounded007});
\item The techniques of the proof of Theorem \ref{thm:Part1} may also be used to show that $(\widetilde{\A},\widetilde{\phi})$
must be $\delta$-semistable. In particular, the higher cohomology groups of $\widetilde{\A}(n)$ vanish, so that
$H^0(\widetilde{q}(n))$ is indeed an isomorphism.
\end{itemize}
The family $(\widetilde{q},\widetilde{\phi})$ is thus induced by a morphism $\widetilde{\eta}$ which lifts $\eta$ and
extends $\eta^\star$. This finishes the argument.
\end{proof}
%%%
Since ${\Bbb G}^{\rm ss}$ possesses a projective quotient, Proposition \ref{prop:Proper} and Lemma
\ref{lem:Gieseker/Ramanathan}
show that the good quotient
$$
{\rm M}^{\delta\hbox{-}\rm ss}_P(\rho):={\frak D}^{\delta\hbox{-}\rm ss}\catqot \SL(U)
$$
exists as a projective scheme. Likewise, the geometric quotient
$$
{\rm M}^{\delta\hbox{-}\rm s}_P(\rho):={\frak D}^{\delta\hbox{-}\rm s}/ \SL(U)
$$
exists as an open subscheme of ${\rm M}^{\delta\hbox{-}\rm ss}_P(\rho)$.
By Propositions \ref{LUP0} and \ref{Glue0} and the universal property
of a categorical quotient, the space ${\rm M}^{\delta\hbox{-}\rm ss}_P(\rho)$ is indeed a coarse
moduli space. \qed
%%%
\begin{Rem}[S-equivalence]
\label{rem:SEqAndOrbEq}
Recall that two points in ${\frak D}^{\delta\hbox{-}\rm ss}$ are mapped to the same point in the quotient
if and only if the closures of their orbits intersect in ${\frak D}^{\delta\hbox{-}\rm ss}$. Given a point
$y\in {\frak D}^{\delta\hbox{-}\rm ss}$, let $y'\in {\frak D}^{\delta\hbox{-}\rm ss}$ be the point whose
orbit is the unique closed orbit in $\ol{\SL(U)\cdot y}$($\subseteq{\frak D}^{\delta\hbox{-}\rm ss}$).
Then, there is a one parameter subgroup $\la\colon {\Bbb G}_m(k)\lra\SL(U)$ with
$\lim_{z\rightarrow\infty}\la(z)\cdot y\in \SL(U)\cdot y'$. For this one parameter subgroup, one has of course
$\mu(\la,y)=0$. Thus, the equivalence relation that we have to consider on the closed points of
${\frak D}^{\delta\hbox{-}\rm ss}$ is generated by $y\sim \lim_{z\rightarrow\infty}\la(z)\cdot y$ for all
one parameter subgroups $\la\colon {\Bbb G}_m(k)\lra\SL(U)$ with $\mu(\la,y)=0$.
\par
If one looks carefully at the arguments given in the proofs of Theorems \ref{thm:Part1} and \ref{thm:Part2}, one sees
that, for a point $y=([q\colon U\otimes{\Oh}_X(-n)\lra \A],\phi)\in {\frak D}^{\delta\hbox{-}\rm ss}$, the
following observations hold true:
\begin{itemize}
\item If $\la\colon {\Bbb G}_m(k)\lra\SL(U)$ verifies $\mu(\la,y)=0$, its weighted
      flag $(U_\bullet(\la),\alpha_\bullet(\la))$ has the property that the weighted filtration
      $(\A_\bullet,\alpha_\bullet(\la))$ with $\A_j:=q(U_j\otimes{\Oh}_X(-n))$, $j=1,\dots,t$, satisfies
      $$
      M\bigl(\A_\bullet,\alpha_\bullet(\la)\bigr)+\delta\cdot \mu(\A_\bullet,\alpha_\bullet(\la);\phi)\equiv 0.
      $$
\item Given a weighted filtration
      $(\A_\bullet,\alpha_\bullet)$ of $\A$ with
      $$
      M(\A_\bullet,\alpha_\bullet)+\delta\cdot \mu(\A_\bullet,\alpha_\bullet;\phi)\equiv 0,
      $$
      one can assume $h^i(\A_j(n))=0$, $i>0$, and that $\A_j(n)$ is globally generated, $j=1,\dots,t$.
      Hence, there is a unique flag $U_\bullet$ in $U$, such that $H^0(q(U_j))$ maps $U_j$ onto
      $H^0(\A_j(n))$, $j=1,\dots,t$.
      Then, any one parameter subgroup $\la\colon {\Bbb G}_m(k)\lra\SL(U)$ with weighted
      flag $(U_\bullet,\alpha_\bullet)$ satisfies $\mu(\la,y)=0$.
\item For a one parameter subgroup $\la$ with $\mu(\la,y)=0$, $y':=\lim_{z\rightarrow\infty}\la(z)\cdot y$,
      and induced weighted filtration $(\A_\bullet,\alpha_\bullet)$ on $\A$,
      the dispo sheaf $(\A_{y'},\phi_{y'})$ is isomorphic to ${\rm df}_{(\A_\bullet,\alpha_\bullet)}(\A,\phi)$.
\end{itemize}
%%%
This shows that the equivalence relation induced by the GIT process on the closed points of ${\frak D}^{\delta\hbox{-}\rm ss}$
is just S-equivalence of dispo sheaves as introduced in Section \ref{subsec:SEquivDispo}.
\end{Rem}
%%%
\begin{Rem}[Decorated vector bundles on curves]
\label{rem:DecVb}
In \cite{Schmitt0}, given a homogeneous representation $\kappa\colon \GL_r(k)\lra\GL(V)$ and a smooth
projective curve $X$ over the complex numbers,
the moduli problem of classifying triples $(E,L,\phi)$ consisting of a vector bundle
of rank $r$ on $X$, a line bundle $L$ on $X$, and a non-trivial homomorphism $\phi\colon E_\kappa\lra L$,
$E_\kappa$ being associated to $E$ via $\kappa$, was solved
by a GIT procedure similar to the one presented above. Write $V=\widetilde{V}\otimes \det(V)^{\otimes -w}$
where $\widetilde{V}$ is a homogeneous polynomial $\GL_r(k)$-module, say, of degree $u$.
The only characteristic zero issue that is necessary for the construction in \cite{Schmitt0} is the fact that
$\widetilde{V}$ can be written as the quotient of $(V^{\otimes u})^{\oplus v}$ for an appropriate
positive integer $v$. In characteristic $p>0$, this is only true when $p>u$. However, we may use the results of Section \ref{subsec:Poly}.
They imply that $\widetilde{V}$ is a quotient of ${\Bbb D}^{u,v}(V)$ for an appropriate integer $v>0$.
Now, ${\Bbb D}^{u,v}(V)$ is a subrepresentation of $(V^{\otimes u})^{\oplus N}$. This shows that the arguments
given in the present paper may be used to deal with decorated vector bundles on smooth projective curves in
any characteristic (see \cite{HeinS} for the analogous case of decorated parabolic vector bundles).
\end{Rem}
%%%
\section{The proof of the main theorem}
%%%
This section is devoted to the proof of the main theorem. In fact, we will prove a slightly stronger
theorem which is the exact analog to the main result of \cite{Schmitt} in arbitrary characteristic.
To do so, we recall the necessary notions of $\widetilde{\delta}$-semistable pseudo $G$-bundles and so on.
%%%
\subsection{Associated dispo sheaves}
\label{subsec:AssDispos}
%%%
The notion of a ``pseudo $G$-bundle" has been recalled in Section \ref{subsec:Pseudo}. Now,
we relate pseudo $G$-bundles to dispo sheaves.
\par
Let $S$ be a scheme, and $(\A_S,\tau_S)$ a family of pseudo $G$-bundles parameterized by $S$.
Let $\iota\colon U\subseteq S\times X$ be the maximal open subset where $\A_S$ is locally free.
The locally free sheaf $\A_{S|U}$ and the $\GL_r(k)$-module ${\Bbb V}_s$ give rise to the vector bundle
${\Bbb V}_s(\A_{S|U})$, and there is a surjection
$$
{{\cal S}ym}^\star\bigl({\Bbb V}_s(\A_{S|U})\bigr)\lra {{\cal S}ym}^{(s!)}(\A_{S|U}\otimes V)^G.
$$
Define $\widetilde{\tau}_s$ as the restriction of $\tau_{S|U}$ to the subalgebra
${{\cal S}ym}^{(s!)}(\A_{S|U}\otimes V)^G$.
Then, $\widetilde{\tau}_s$ is determined by a homomorphism
$$
\phi^\p\colon {\Bbb V}_s(\A_{S|U}) \lra {\Oh}_U.
$$
Thus, $\tau_{S|U}$ gives rise to the homomorphism
$$
\phi_S\colon {\Bbb V}_s(\A_{S})\lra \iota_\star\bigl({\Bbb V}_s(\A_{S|U})\bigr)\lra
\iota_\star({\Oh}_U)={\Oh}_{S\times X},
$$
by Corollary \ref{NatTrans}.
Therefore, we can associate to the family $(\A_S,\tau_S)$ of pseudo $G$-bundles the family
$(\A_S,\phi_S)$ of dispo sheaves of type $(\rho,s)$.
\par
The map which associates to a pseudo $G$-bundle a dispo sheaf
is injective on isomorphy classes. More precisely, we find
%%%
\begin{Lem}
\label{lem:inj}
Suppose that $(\A,\tau)$ and $(\A,\tau')$ are two pseudo $G$-bundles, such that
the associated dispo sheaves are equal. Then, there is a root of unity $\zeta\in k$,
such that $\zeta\cdot\id_\A$ yields an isomorphism between $(\A,\tau)$ and $(\A,\tau')$.
\end{Lem}
%%%
\begin{proof}
For $d>0$, let
$$
\tau_d,\tau_d'\colon {\cal S}ym^d(\A\otimes V)^G\lra {\Oh}_X
$$
be the degree $d$ component of $\tau$ and $\tau'$, respectively. Note that
$\tau$ is determined by $\bigoplus_{d=1}^s\tau_d$.
Let
$$
\widehat{\tau}_s:\bigoplus_{(d_1,\dots,d_s);\atop d_i\ge 0,\sum i d_i=s!}\Bigl({\cal S}ym^{d_1}\bigl((\A\otimes V)^G\bigr)\otimes\cdots\otimes{\cal S}ym^{d_s}\bigl({\cal S}ym^{s}(\A\otimes V)^G\bigr)\Bigr)\lra{\Oh}_X
$$
be the map induced by $\tau_1$,\dots,$\tau_s$, and define $\widehat{\tau}_s'$ in a similar way.
By definition, $\widehat{\tau}_{s|U}=\widetilde{\tau}_s$.
Our assumption thus grants that $(\A,\widehat{\tau}_s)$ and $(\A,\widehat{\tau}_s')$ are equal.
This implies that, for $1\le d\le s$,
$$
{\cal S}ym^{\frac{s!}{d}}(\tau_d)={\cal S}ym^{\frac{s!}{d}}(\tau_d').
$$
Restricting this equality to the generic point, it follows that
there is an $(s!/d)$th root of unity $\zeta_d$ with
$$
\tau_d'=\zeta_d\cdot \tau_d,\q d=1,\dots,s.
$$
It remains to show that there is an $s!$th root of unity $\zeta$,
such that $\zeta_d=\zeta^d$. To see this, let ${\Bbb A}$ be the
restriction of $\A$ to the generic point. Then, $\widehat{\tau}_s$
and $\widehat{\tau}_s'$, restricted to the generic point, define the
same point
$$
x\in {\Pe}:={\Pe}\Bigl(\bigoplus_{(d_1,\dots,d_s);\atop d_i\ge 0,\sum i d_i=s!}{\rm Sym}^{d_1}\bigl(({\Bbb A}\otimes V)^G\bigr)\otimes\cdots\otimes
{\rm Sym}^{d_s}\bigl({\rm Sym}^{s}({\Bbb A}\otimes V)^G\bigr)\Bigr).
$$
On the other hand, $\bigoplus_{d=1}^s\tau_d$ and $\bigoplus_{d=1}^s\tau'_d$
define points
$$
y,y'\in {\Bbb B}:=\Bigl(\bigoplus_{d=1}^s{\rm Sym}^d({\Bbb A}\otimes V)^G\Bigr)^\vee.
$$
By our assumption, $y$ and $y'$ map both to $x$ under the quotient map
followed by the Veronese embedding
$$
{\Bbb B}\setminus\{0\}\lra \bigl({\Bbb B}\setminus\{0\}\bigr)\catqot {\Bbb G}_m(K)\hookrightarrow\P.
$$
Putting all the information we have gathered so far together, we find the claim
about the $\zeta_i$ and from that the one of the lemma.
\end{proof}
%%%
Let $\widetilde{\delta}\in\Q[x]$ be a positive polynomial of degree at most $\dim(X)-1$.
We choose an $s$ as before and define $\delta:=\widetilde{\delta}/s!$.
A pseudo $G$-bundle is said to be \it $\widetilde{\delta}$-(semi)stable\rm, if the associated
dispo sheaf $(\A,\phi)$ of type $(\rho,s)$ is $\delta$-(semi)stable.
Similarly, given a non-negative rational number $\delta^\star$, we define the pseudo $G$-bundle
$(\A,\tau)$ to be $\delta^\star$-slope (semi)stable, if the associated dispo sheaf
$(\A,\phi)$ is $(\delta^\star/s!)$-slope (semi)stable.
%%%
\begin{Rem}
i) The definition of $\widetilde{\delta}$-(semi)stability is the same as the one given in
\cite{Schmitt}, p.\ 1192.
\par
ii) Using (\ref{eq1:TensorPower}),
it follows that the notion of $\widetilde{\delta}$-(semi)stability does not depend on the choice of
$s$. This is why we threw in the factor $1/s!$.
\end{Rem}
%%%
\subsection{S-equivalence}
\label{subsec:SEqPseudo}
%%%
We fix a stability parameter $\widetilde{\delta}$, i.e., a positive rational polynomial of degree at most $\dim(X)-1$.
Suppose $(\A,\tau)$ is a $\widetilde{\delta}$-semistable pseudo $G$-bundle with
associated dispo sheaf $(\A,\phi)$ and $(\A_\bullet,{\alpha}_\bullet)$
is a weighted filtration with
$$
M(\A_\bullet,{\alpha}_\bullet)+\frac{\widetilde{\delta}}{s!}\cdot \mu(\A_\bullet,{\alpha}_\bullet;\phi)\equiv 0.
$$
The construction used for defining an associated admissible deformation of a dispo sheaf can be easily extended
to give the construction of the associated admissible deformation ${\rm df}_{(\A_\bullet,{\alpha}_\bullet)}(\A,\tau)$.\
As before, we let \it S-equivalence \rm be the equivalence relation ``$\sim_{\rm S}$" on
$\widetilde{\delta}$-semistable pseudo $G$-bundles $(\A,\tau)$ generated by
$$
(\A,\tau)\sim_{\rm S}{\rm df}_{(\A_\bullet,{\alpha}_\bullet)}(\A,\tau).
$$
The injectivity of the map which assigns to the isomorphy class of a pseudo $G$-bundle the isomorphy class of
the associated dispo sheaf (Lemma \ref{lem:inj}) and the definitions of semistability for the respective objects
easily imply that for two pseudo $G$-bundles $(\A,\tau)$ and $(\A^\p,\tau^\p)$ with associated dispo sheaves
$(\A,\phi)$ and $(\A^\p,\phi^\p)$
one has:
%%%%
\begin{equation}
\label{SEquivCompatible}
(\A,\tau)\sim_{\rm S}(\A^\p,\tau^\p)\q\Longleftrightarrow\q
(\A,\phi)\sim_{\rm S}(\A^\p,\phi^\p).
\end{equation}
%%%
In Remark \ref{rem:SEqSingPrinz} below, we will give a nice description of S-equivalence on semistable singular principal $G$-bundles.
%%%
\subsection{Moduli spaces for $\widetilde{\delta}$-semistable pseudo $G$-bundles}
\label{sub:Moduli1}
%%%
An immediate consequence of the definition of semistability of pseudo $G$-bundles and Theorem \ref{thm:DispShBound}
is that, for a given Hilbert polynomial $P$, the set of torsion free sheaves $\A$ with Hilbert polynomial $P$ for which
there exists a $\widetilde{\delta}$-(semi)stable pseudo $G$-bundle $(\A,\tau)$ is bounded. Finally,
the construction carried out in Section \ref{subsec:AssDispos} and Corollary \ref{NatTrans} give a natural transformation
$$
{\rm AD}\colon \ul{\rm M}^{\widetilde{\delta}\rm\hbox{-}(s)s}_P(\rho)\lra
\ul{\rm M}_P^{\delta\rm\hbox{-}(s)s}(\rho,s)
$$
of the functor $\ul{\rm M}^{\widetilde{\delta}\rm\hbox{-}(s)s}_P(\rho)$
which assigns to a scheme $S$ the set of isomorphy classes of families of $\widetilde{\delta}$-(semi)stable pseudo
$G$-bundles with Hilbert polynomial $P$ parameterized by $S$ into the functor
$\ul{\rm M}^{{\delta}\rm\hbox{-}(s)s}_P(\rho,s)$ which assigns to a scheme $S$ the set of isomorphy classes of
families of ${\delta}$-(semi)stable dispo sheaves of type $(\rho,s)$ with Hilbert polynomial $P$ parameterized by $S$.
%%%
\begin{Thm}
\label{ModDeltaSemStab}
Fix the stability parameter $\widetilde{\delta}$ and the Hilbert polynomial $P$.
Then, there is a {\bfseries projective} scheme ${\rm M}^{\widetilde{\delta}\rm\hbox{-}ss}_P(\rho)$ which is
a coarse moduli space for the functors $\ul{\rm M}_P^{\widetilde{\delta}\hbox{-}\rm ss}(\rho)$.
\end{Thm}
%%%
\begin{Rem}
This theorem generalizes the main theorem of \cite{Schmitt} to arbitrary characteristic.
\end{Rem}
%%%
\subsubsection{Construction of the parameter space}
%%%
There is a constant $C$, such that $\mu_{\max}(\A)\le C$ for every
$\widetilde{\delta}$-semistable pseudo $G$-bundle $(\A,\tau)$ with
$P(\A)=P$, i.e., $\A$ lives in a bounded family. Thus, we may
choose the integer $s$ in such a way that ${\cal
S}ym^\star(\A\otimes V)^G$ is generated by elements in degree at
most $s$ for all such $\A$. We choose an $n_0\gg 0$ with the following
properties: For every sheaf $\A$ with Hilbert polynomial $P$ and
$\mu_{\max}(\A)\le C$ and every $n\ge n_0$, one has
%%%
\begin{itemize}
\item $H^i(\A(n))=\{0\}$ for $i>0$;
\item $\A(n)$ is globally generated;
\item The construction of the moduli space of $(\widetilde{\delta}/s!)$-semistable dispo sheaves
      of type $(\rho,s)$ can be performed with respect to $n$.
\end{itemize}
%%%
We choose a $k$-vector space $U$ of dimension $P(n)$. Let ${\frak Q}$ be the quasi-projective scheme which parameterizes quotients $q\colon U\otimes{\Oh}_X(-n)\lra \A$ where $\A$ is a torsion free sheaf with Hilbert polynomial $P$ and $\mu_{\max}(\A)\le C$ (so that ${\cal S}ym^\star(\A\otimes V)^G$ is generated by elements of degree at most $s$) and $H^0(q(n))$ is an isomorphism. Let
$$
{\frak q}_{{\frak Q}}\colon U\otimes\pi_X^\star\bigl({\Oh}_X(-n)\bigr)\lra {\cal A}_{{\frak Q}}
$$
be the universal quotient. By Lemma \ref{SurJect}, there is a homomorphism
$$
{\cal S}ym^\star\bigl( U\otimes\pi_X^\star({\Oh}_X(-n))\otimes
V\bigr)^G \lra {\cal S}ym^\star\bigl(\A_{{\frak Q}}\otimes
V\bigr)^G
$$
which is surjective where $\A_{\frak Q}$ is locally free.
For a point $[q\colon U\otimes{\Oh}_X(-n)\lra \A]\in {\frak Q}$, any
homomorphism $\tau\colon {\cal S}ym^\star(\A\otimes V)^G \lra{\Oh}_X$
of ${\Oh}_X$-algebras is determined by the composite homomorphism
$$
\bigoplus_{i=1}^s {\cal S}ym^i\bigl( U\otimes\pi_X^\star({\Oh}_X(-n))\otimes V\bigr)^G\lra {\Oh}_X
$$
of ${\Oh}_X$-modules. Noting that
$$
{\cal S}ym^i\bigl( U\otimes\pi_X^\star({\Oh}_X(-n))\otimes V\bigr)^G
\cong
{\rm Sym}^i(U\otimes V)^G\otimes \pi_X^\star\bigl({\Oh}_X(-in)\bigr),
$$
$\tau$ is determined by a collection of homomorphisms
$$
\phi_i\colon {\rm Sym}^i(U\otimes V)^G\otimes  {\Oh}_X\lra {\Oh}_X(in),\q i=1,\dots,s.
$$
Since $\phi_i$ is determined by the induced linear map on global sections, we will
construct the parameter space inside
$$
{\frak Y}^0:={\frak Q} \times \bigoplus_{i=1}^{s}{\rm Hom}
\Bigl({\rm Sym}^i (U\otimes V)^G, H^0 \bigl({\cal O}_X
(in)\bigr)\Bigr).
$$
Note that, over ${\frak Y}^0\times X$,
there are universal homomorphisms
$$
\widetilde{\varphi}^i:{\rm Sym}^i (U\otimes V)^G \otimes {\cal O}_
{{\frak Y}^0 \times X} \lra H^0 \bigl({\cal O}_X (in)\bigr) \otimes
{\cal O}_{{\frak Y}^0 \times X},\q i=1,\dots, s.
$$
Let $\varphi^i = {\rm ev} \circ \tilde{\varphi}^i$ be the composition of
$\widetilde{\phi}^i$ with the evaluation map
${\rm ev}\colon H^0({\Oh}_X(in))\otimes{\Oh}_{{\frak Y}^0\times X}\lra \pi_X^\star({\cal O}_X (in))$, $i=1,\dots,s$.
We twist $\varphi^i$ by $\id_{\pi_X^\star({\Oh}_X(-in))}$ and put the resulting maps together to  the  homomorphism
$$
\varphi\colon {\cal V}_{{\frak Y}^0} :=
\bigoplus_{i=1}^s {\cal S}ym^i\bigl(U\otimes \pi_X^\star({\cal O}_X(- n))\otimes V\bigr)^G \lra
{\cal O}_{{\frak Y}^0\times X}.
$$
Next, $\phi$ yields a homomorphism of ${\cal O}_{{\frak Y}^0\times X}$-algebras
$$
\widetilde{\tau}_{{\frak Y}^0}\colon {\cal S}ym^\star ({\cal V}_{{\frak Y}^0}) \lra
{\cal O}_{ {\frak Y}^0 \times X}.
$$
On the other hand, there is a surjective homomorphism
$$
\beta\colon {\cal S}ym^\star ({\cal V}_{{\frak Y}^0}) \lra {\cal
S}ym^\star(\pi^\star({\cal A}_{{\frak  Q}})\otimes V)^G
$$
of graded algebras where the left-hand algebra is graded by
assigning the weight $i$ to the elements in ${\cal S}ym^i(U\otimes
\pi_X^\star({\Oh}_X(-n))\otimes V)^G$. Here, $\pi\colon {\frak
Y}^0\times X\lra {\frak Q}\times X$ is the natural projection. The
parameter space ${\frak Y}$ is defined by the condition that
$\widetilde{\tau}_{{\frak Y}^0}$ factorizes over $\beta$, i.e.,
setting $ {\cal A}_{{\frak Y}}:=(\pi^\star({\cal A}_{{\frak Q}}))_{|{\frak Y}\times X} $, there is a homomorphism
$$
{\tau}_{{\frak Y}} \colon {\cal S}ym^\star ({\cal A}_{{\frak Y}}\otimes V)^G \lra
{\cal O}_{{\frak Y}\times X}
$$
with $\widetilde{\tau}_{{\frak Y}^0|{\frak Y}\times X}={\tau}_{{\frak Y}}\circ\beta$.
Formally, ${\frak Y}$ is defined as the scheme theoretic intersection of the closed subschemes
$$
{\frak Y}_d:=\bigl\{\,y\in {\frak Y}^0\,|\, \widetilde{\tau}^d_{{\frak Y}^0|\{y\}\times X}
\colon \ker\bigl(\beta_{|\{y\}\times X}^d\bigr)\lra {\Oh}_X\ \hbox{\rm is trivial}\,\bigr\}, \q d\ge 0.
$$
The family
$(\A_{{\frak Y}}, \tau_{{\frak Y}})$ is the \it universal family of pseudo $G$-bundles parameterized by ${\frak Y}$\rm.
(In all these constructions, one needs to use Lemma \ref{lem:AllowsBaseChange}.)
%%%
\begin{Prop}[Local universal property]
\label{LUP} Let $S$ be a scheme and $(\A_S,\tau_S)$ a family of
$\widetilde{\delta}$-semi\-stable pseudo $G$-bundles with Hilbert
polynomial $P$ parameterized by $S$. Then, there exist a covering
of $S$ by open subschemes $S_i$, $i\in I$, and morphisms
$\beta_i\colon S_i\lra {\frak Y}$, $i\in I$, such that the family
$(\A_{S|S_i},\tau_{S|S_i})$ is isomorphic to the pullback of the
universal family on ${\frak Y}\times X$ by $\beta_i\times\id_X$
for all  $i\in I$.
\end{Prop}
%%%
\subsubsection{The group action}
%%%
There is a natural action of $\GL(U)$ on the quot scheme ${\frak
Q}$ and on ${\frak Y}^0$. This action leaves the closed subscheme
${\frak Y}$ invariant, and therefore yields an action
$$
\Gamma\colon \GL(U)\times {\frak Y}\lra {\frak Y}.
$$
\begin{Prop}[Gluing property]
\label{Glue} Let $S$ be a scheme and $\beta_i\colon S\lra {\frak
Y}$, $i=1,2$, two morphisms, such that the pullback of the
universal family via $\beta_1\times \id_X$ is isomorphic to its
pullback via $\beta_2\times \id_X$. Then, there is a morphism
$\Xi\colon S\lra \GL(U)$, such that $\beta_2$ equals the morphism
$$
S\stackrel{\Xi\times\beta_1}{\lra}\GL(U)\times {\frak Y}\stackrel{\Gamma}{\lra} \frak Y.
$$
\end{Prop}
%%%
\begin{Rem}
The universal family is equipped with a $\GL(U)$-linearization. If
one fixes, in the above proposition, an isomorphism between its
pullbacks via $\beta_1\times \id_X$ and $\beta_2\times \id_X$,
then there is a unique morphism $\Xi\colon S\lra \GL(U)$ which
satisfies the stated properties and, in addition, that the given
isomorphism is induced by pullback via $(\Xi\times\beta_1\times
\id_X)$ from the linearization of $(\A_{{\frak Y}}, \tau_{{\frak
Y}})$. This fact simply expresses that the moduli stack for
$\widetilde{\delta}$-(semi)stable pseudo $G$-bundles will be the quotient
stack of an appropriate open subscheme of the parameter space
${\frak Y}$.
\end{Rem}
%%%
\subsubsection{Conclusion of the proof}
%%%
Suppose we knew that the points $([q\colon U\otimes{\Oh}_X(-n)\lra\A],\tau)$ in the parameter space
${\frak Y}$ for which $(\A,\tau)$ is $\widetilde{\delta}$-semistable form an open subscheme
${\frak Y}^{\widetilde{\delta}\rm\hbox{-}ss}$. Then, it suffices to show that
${\frak Y}^{\widetilde{\delta}\rm\hbox{-}ss}$
possesses a (good, uniform)
categorical quotient by the action of $\GL(U)$. Indeed, Propositions \ref{LUP} and \ref{Glue}
and the universal property  of the categorical quotient then
imply that ${\rm M}_P^{\widetilde{\delta}\rm\hbox{-}ss}(\rho):={\frak Y}^{\widetilde{\delta}\rm\hbox{-}ss}\catqot \GL(U)$
has the desired properties.
We have the natural surjection ${\Bbb G}_m(k)\times \SL(U)\lra \GL(U)$, $(z,m)\lma z\cdot m$, and obviously
$$
{\frak Y}^{\widetilde{\delta}\rm\hbox{-}ss}\catqot \GL(U)=
{\frak Y}^{\widetilde{\delta}\rm\hbox{-}ss}\catqot \bigl({\Bbb G}_m(k)\times \SL(U)\bigr).
$$
By Example \ref{uniform}, ii), we may first form
$$
\ol{\frak Y}^{\widetilde{\delta}\rm\hbox{-}ss}:={\frak Y}^{\widetilde{\delta}\rm\hbox{-}ss}\catqot {\Bbb G_m}(k)
$$
and then
$$
\ol{\frak Y}^{\widetilde{\delta}\rm\hbox{-}ss}\catqot \SL(U).
$$
We can easily form the quotient $\ol{\frak Y}:=
{\frak Y}\catqot {\Bbb G}_m(k)$. Since ${\Bbb G}_m(k)$ is linearly reductive,
$\ol{\frak Y}$ is a closed subscheme of
$$
{\frak Q}\times \biggl( \bigoplus_{i=1}^{s}{\rm Hom} \Bigl({\rm Sym}^i (U\otimes V)^G, H^0 \bigl({\cal O}_X (in)\bigr)\Bigr)
\big/\hskip-2pt\big/ {\Bbb G}_m(k)\biggr).
$$
In particular, $\ol{\frak Y}$ is projective over ${\frak Q}$.
Let ${\frak D}\lra {\frak Q}$ be the parameter space for dispo sheaves of type $(\rho,s)$
constructed above. If we apply the construction described in Section
\ref{subsec:AssDispos} to the universal family $(\A_{{\frak Y}},\tau_{{\frak Y}})$, we get an
$\SL(U)$-equivariant and ${\Bbb G}_m(k)$-invariant morphism
$$
\widetilde{\psi}\colon {\frak Y}\lra {\frak D}
$$
and, thus, a proper $\SL(U)$-equivariant morphism
$$
\psi\colon \ol{\frak Y}\lra {\frak D}.
$$
By Lemma \ref{lem:inj}, $\psi$ is injective, so that it is even {\bfseries finite}.
Now, there are open subsets ${\frak D}^{\delta\rm\hbox{-}(s)s}$, $\delta:=\widetilde{\delta}/s!$,
which parameterize the $\delta$-(semi)stable dispo sheaves of type
$(\rho,s)$, such that the good, uniform categorical quotient
$$
{\rm M}^{\delta\rm\hbox{-}ss}_P(\rho,s)={\frak
D}^{\delta\rm\hbox{-}ss}\catqot\SL(U)
$$
exists as a projective scheme and the geometric, uniform
categorical quotient
$$
{\rm M}^{\delta\rm\hbox{-}s}_P(\rho,s)={\frak D}^{\delta\rm\hbox{-}s}/\SL(U)
$$
as an open subscheme of ${\rm M}^{\delta\rm\hbox{-}ss}_P(\rho,s)$. By definition of semistability,
$$
\widetilde{\psi}^{-1}({\frak D}^{\delta\rm\hbox{-}ss})
=
{\frak Y}^{\widetilde{\delta}\rm\hbox{-}ss},
$$
whence
$$
{\psi}^{-1}({\frak D}^{\delta\rm\rm\hbox{-}ss})
=
{\frak Y}^{\widetilde{\delta}\rm\hbox{-}ss}\catqot {\Bbb G}_m(k).
$$
Now, Lemma \ref{lem:Gieseker/Ramanathan} implies that the quotient
$$
{\rm M}^{\widetilde{\delta}\hbox{-}\rm ss}_P(\rho):={\frak Y}^{\widetilde{\delta}\rm\hbox{-}ss}\catqot \GL(U)
=
{\psi}^{-1}({\frak D}^{\delta\rm\hbox{-}ss})\catqot \SL(U)
$$
exists as a projective scheme.
Likewise, the open subscheme
$$
{\rm M}^{\widetilde{\delta}\rm\hbox{-}s}_P(\rho):={\frak Y}^{\widetilde{\delta}\rm\hbox{-}s}/ \GL(U)
=
{\psi}^{-1}({\frak D}^{\delta\rm\hbox{-}s})/ \SL(U)
$$
is a uniform (universal) geometric quotient and an open subscheme of
${\rm M}^{\widetilde{\delta}\rm\hbox{-}ss}_P(\rho)$.
\qed
%%%
\subsection{Semistable singular principal bundles}
%%%
\begin{Thm}
\label{HonestAlwaysSemStab}
Fix a Hilbert polynomial $P$, and let $\delta_\infty$ be as in Corollary {\rm\ref{cor:AsymBound}}. For every
polynomial $\widetilde{\delta}\succ s!\cdot {\delta}_\infty$ and every singular
principal $G$-bundle $(\A,\tau)$ with $P(\A)=P$, the following properties are equivalent:
\begin{itemize}
\item[{\rm i)}] $(\A,\tau)$ is (semi)stable.
\item[{\rm ii)}] $(\A,\tau)$ is $\widetilde{\delta}$-(semi)stable.
\end{itemize}
\end{Thm}
%%%
Taking into account Corollary \ref{cor:AsymBound}, the theorem reduces to:
%%%
\begin{Lem}
\label{RecognizeFiltration}
Let $(\A,\tau)$ be a singular principal $G$-bundle with associated dispo sheaf $(\A,\phi)$.
Then, for a weighted filtration $(\A_\bullet,{\alpha}_\bullet)$ of $\A$, the condition
$$
\mu(\A_\bullet,{\alpha}_\bullet; \phi)=0
$$
is satisfied if and only if
$$
(\A_\bullet,{\alpha}_\bullet)=\bigl(\A_\bullet(\beta), {\alpha}_\bullet(\beta)\bigr)
$$
for some reduction $\beta$ of $(\A,\tau)$ to a one parameter subgroup $\la$ of $G$.
\end{Lem}
%%%
\begin{proof}
We show that the first condition implies the second one, the converse being an easy exercise.
Let $\la^\p\colon {\Bbb G}_m(k)\lra \SL_r(V)$ be a one parameter subgroup, such that the associated weighted flag
$$
\bigl(V_\bullet(\la')\colon 0\subsetneq
V_1\subsetneq\cdots\subsetneq V_{t^\p}\subsetneq
V,{\alpha}_\bullet(\la')\bigr)
$$
satisfies $t^\p=t$, $\dim_k(V_i)=\rk\A^\p_i$,
$\A^\p_i=\ker(\A^\vee\lra \A_{t+1-i}^\vee)$, $i=1,\dots,t$, and
${\alpha}_\bullet(\la')=(\alpha_t,\dots,\allowbreak \alpha_1)$, if
${\alpha}_\bullet=(\alpha_1,\dots,\alpha_t)$. Then, the weighted
filtration $(\A_\bullet,{\alpha}_\bullet)$ is associated to a
reduction $\beta^\p$ of the principal $\GL(V)$-bundle ${\cal
I}som(V\otimes\Oh_{U'}, \A^\vee_{|U^\p})$ to $\la'$ with $U^\p$ the maximal
open subset where $\A$ is locally free and all the $\A^\p_i$ are
subbundles. We may choose an open subset $\widetilde{U}\subseteq
X$, such that there is a trivialization $\psi\colon
\A^\vee_{|\widetilde{U}}\lra V\otimes {\Oh}_{\widetilde{U}}$ with
$\psi(\A^\p_i)=V_i\otimes {\Oh}_{\widetilde{U}}$, $i=1,\dots,t$. By definition of
the number $\mu(\A_\bullet,{\alpha}_\bullet;\allowbreak \phi)$,
(\ref{eq1:NochEinKomp0}), and
Proposition \ref{prop1:Recognize1}, we see that there is a one
parameter subgroup $\la\colon {\Bbb G}_m(k)\lra G$, such that
$$
\bigl(V_\bullet(\la),{\alpha}_\bullet(\la)\bigr)= \bigl(V_\bullet(\la^\p),{\alpha}_\bullet(\la^\p)\bigr).
$$
To the principal bundles ${\cal P}(\A,\tau)$ and ${\cal I}som(V\otimes\Oh_U, \A^\vee_{|U})$, we may associate
group schemes ${\cal G}\subset {\cal G\cal L}(V)$ over $U$. Now, ${\cal G\cal L}(V)$ acts on
${\cal I}som(V\otimes\Oh_U, \A^\vee_{|U})/{\cal Q}_{\GL(V)}(\la)$, and the stabilizer of the section
$\beta^\p\colon U^\p\lra {\cal I}som(V\otimes\Oh_{U'}, \A^\vee_{|U'})/Q_{\GL(V)}(\la)$ is a parabo\-lic subgroup
${\cal Q}\subset {\cal G\cal L}(V)_{|U^\p}$, such that
$$
{\cal G\cal L}(V)_{|U^\p}/{\cal Q}={\cal I}som(V\otimes\Oh_{U'}, \A^\vee_{|U^\p})/Q_{\GL(V)}(\la).
$$
The intersection ${\cal Q}_{\cal G}:={\cal Q}\cap {\cal
G}_{|U^\p}$ is a parabolic subgroup. This follows if one applies the above reasoning on weighted
flags to the geometric fibers of ${\cal G}\subset{\cal G\cal L}(V)$ over $U'$. Furthermore, ${\cal G}_{|U^\p}/{\cal
Q}_{\cal G}={\cal P}(\A,\tau)_{|U^\p}/Q_G(\la)$. This can be seen as follows:
Let $C$ be the set of conjugacy classes of parabolic subgroups of $G$. There is a scheme ${\cal P}ar_p ({\cal
G}_{|U^\p}/U^\p)$ over $U^\p$, such that giving a parabolic
subgroup ${\cal Q}_{\cal G}$ of ${\cal G}_{|U^\p}$ the fibers of
which belong to $p\in C$ is the same as giving a section $U^\p\lra
{\cal P}ar_p({\cal G}_{|U^\p}/U^\p)$ (\cite{Dem}, p.\ 443ff). It is
easy to see that ${\cal P}ar_p ({\cal G}_{|U^\p}/U^\p)\cong {\cal
P}(\A,\tau)_{|U^\p}/Q_p$, $Q_p$ being a representative for $p$
(compare \cite{RamRam}, p.\ 281). Finally, ${\cal G}_{|U^\p}/{\cal
Q}_{\cal G}\cong {\cal P}ar_p({\cal G}_{|U^\p}/U^\p)$ (\cite{Dem},
Corollaire 3.6, page 445). Therefore, we have the commutative
diagram
$$
\xymatrix{
{{\cal Q}_{\cal G}} \ar@{^{(}->}[r]\ar@{^{(}->}[d] & {{\cal Q}}_{\phantom{\cal G}} \ar@{^{(}->}[d]
\\
{{\cal G}_{|U^\p}} \ar@{^{(}->}[r] & {{\cal G\cal L}(V)_{|U^\p}}.
}
$$
Taking ${\cal Q}_{\cal G}$-quotients in the left-hand column and
${\cal Q}$-quotients in the right-hand column yields the commutative diagram
$$
\xymatrix{
U^\p \ar@{=}[r] \ar[d]_{\beta} & U^\p \ar[d]^{\beta'}
\\
{{\cal P}(\A,\tau)_{|U^\p}/Q_G(\la)}\ \ar@{^{(}->}[r] & {{\cal I}som(V\otimes\Oh_{U'}, \A^\vee_{|U'})/Q_{\GL(V)}(\la)}
}
$$
and settles the claim.
\end{proof}
%%%
\begin{Rem}[S-equivalence for semistable singular principal $G$-bundles]
\label{rem:SEqSingPrinz}
Let $(\A,\tau)$ be a semi\-sta\-ble singular principal $G$-bundle. By Lemma \ref{RecognizeFiltration},
an admissible deformation is associated to a reduction $\beta\colon U^\p\lra {\cal P}(\A,\tau)_{|U^\p}/Q_G(\la)$ to a one parameter subgroup, such that
$M(\A_\bullet(\beta),\alpha_\bullet(\beta))\equiv 0$.
The structure of the rational principal $G$-bundle
${\cal P}({\rm df}_{(\A_\bullet,{\alpha}_\bullet)}(\A,\allowbreak \tau))$ may be described in the following way:
The reduction $\beta$ defines a principal $Q_G(\la)$-bundle ${\cal Q}$ over $U^\p$, such that
${\cal P}(\A,\tau)_{|U^\p}$ is obtained from ${\cal Q}$ by means of extending the structure group via
$Q_G(\la)\subset G$. Extending the structure group of ${\cal Q}$ via $Q_G(\la)\lra L_G(\la)\subset G$ yields the
principal bundle ${\cal P}({\rm df}_{(\A_\bullet,{\alpha}_\bullet)}(\A,\allowbreak \tau))_{|U^\p}$.
Thus, our notion of S-equivalence naturally extends the one considered by Ramanathan (see, e.g., \cite{Ramanathan}).
\par
Fix a Hilbert polynomial $P$ and a stability parameter $\widetilde{\delta}$.
The most important basic fact which has to be kept in mind is that
a $\widetilde{\delta}$-semistable pseudo $G$-bundle $(\A,\tau)$ with $P(\A)=P$ which is S-equivalent to a semistable
singular principal $G$-bundle $(\A',\tau')$ is itself a semistable singular principal $G$-bundle.
In other words, the class of semistable singular principal $G$-bundles with Hilbert polynomial $P$ is closed
under S-equivalence inside the class of $\widetilde{\delta}$-semistable pseudo $G$-bundles with Hilbert polynomial $P$.
\end{Rem}
%%%
We now come to the statement which grants semistable reduction
theorem and, in particular, projectivity of the moduli spaces of
semistable singular principal $G$-bundles.
%%%
\begin{Thm}
Assume that either $\rho\colon G\lra \GL(V)$ is of low separable index or $G$ is an adjoint group, $\rho$ is the adjoint
representation and it is of low height. Then, for every polynomial $\widetilde{\delta}$ with $\widetilde{\delta}\succ s!\cdot {\delta}_\infty$, a $\widetilde{\delta}$-semistable pseudo $G$-bundle $(\A,\tau)$ with $P(\A)=P$ is a singular principal $G$-bundle.
\end{Thm}
%%%%
\begin{proof}
Let $(\A,\tau)$ be a pseudo $G$-bundle with associated dispo sheaf $(\A,\phi)$.
Write ${\Bbb A}$ for the restriction of $\A$ to the generic point of $X$. 
As in Section \ref{subsub:Spec}, we set
\begin{eqnarray*}
{\Bbb V}_s({\Bbb A})&:=&\bigoplus_{(d_1,\dots,d_s):\atop d_i\ge 0, \sum i d_i=s!}
\Bigl({\rm Sym}^{d_1}\bigl(({\Bbb A}\otimes_k V)^G\bigr)\otimes \cdots \otimes
{\rm Sym}^{d_s}\bigl({\rm Sym}^s({\Bbb A}\otimes_k V)^G\bigr)\Bigr)
\\
{\Bbb W}_s({\Bbb A})&:=&\bigoplus_{i=1}^s \bigl({\rm Sym}^i({\Bbb A}\otimes_k V)^G\bigr)^\vee.
\end{eqnarray*}
The restriction of $\phi$ to the generic point yields an element
$v\in \P({\Bbb V}_s({\Bbb A}))$ and the restriction of $\tau$ to the generic point an element $w\in {\Bbb W}_s({\Bbb A})$.
Note that there is the surjection  
$$
{\Bbb W}_s({\Bbb A})\setminus\{0\}\lra \P\bigl({\Bbb V}_s({\Bbb A})\bigr),
$$
such that the point $w$ maps to $v$.
\par
Let $\la\colon {\Bbb G}_m(K)\lra \SL({\Bbb A})$ be the one parameter subgroup from Corollary \ref{cor:InstabFlag} with $\mu(\la,w)<0$. According to Lemma \ref{lem1:SemStab1002}, we also have $\mu(\la, v)<0$. Let $({\Bbb A}_\bullet,\alpha_\bullet)$ be the weighted flag of $\la$ in ${\Bbb A}$. We may find a weighted filtration $(\A_\bullet,\alpha_\bullet)$ of $\A$ whose restriction to the generic point yields $({\Bbb A}_\bullet,\alpha_\bullet)$. For this weighted filtration, we find
$$
\mu(\A_\bullet,\alpha_\bullet;\tau)=\frac{1}{s!}\cdot \mu(\la,v)<0.
$$
By Corollary \ref{cor:AsymBound}, the weighted filtration $(\A_\bullet,\alpha_\bullet)$ contradicts $\widetilde{\delta}$-semistability.
\end{proof}
%%%
\begin{Ex}
\label{ex:Classical}
The above proof can also be used for classical groups with their standard representations. Assume, for example, that
$G={\rm Sp}_n(k)\subset \SL_{2n}(k)$. If $\A$ is a torsion free sheaf, then giving $\tau$ is the same as giving $\tau_{|U}$, $U$ being the maximal open subset where $\A$ is locally free. Now, giving $\tau_{|U}$ is equivalent to giving a non-trivial anti-symmetric form $\phi_U\colon \A_{|U}\lra \A_{|U}^\vee$. Since $\A^\vee$ is reflexive, the datum of $\phi_U$ is the same as the datum of an anti-symmetric form $\phi\colon \A\lra\A^\vee$. Assume that $(\A,\phi)$ is a singular principal ${\rm Sp}_n(k)$-bundle. Then, ${\cal B}:=\ker(\phi)$ is a proper saturated subsheaf. If it is non-trivial, then the restriction of the weighted filtration $(0\subsetneq {\cal B}\subsetneq \A, (1))$ to the generic point will come from a one parameter subgroup $\la\colon {\Bbb G}_m(K)\lra \SL({\Bbb A})$ with $\mu(\la, v)<0$. Therefore, the theorem holds for ${\rm Sp}_n(k)$ with its standard representation in any characteristic.
\par
A similar reasoning can be applied to ${\rm SO}_n(k)$, if the characteristic of $k$ is not two. It works also for ${\rm GO}_n(k)$ and ${\rm GSp}_n(k)$, if one uses the moduli construction suggested in the introduction.
\end{Ex}
%%%
\subsection{Proof of the semistable reduction theorem}
%%%
Before going into the proof, we need to recall the following result
of Seshadri\footnote{A proof will be given in the appendix.} which can be thought of as the semistable reduction
theorem for GIT quotients.
%%%
\begin{Thm}[Seshadri \cite{Sesh2}, Theorem 4.1]\label{Seshadri's thm}
Let $(X, L)$ be a polarized projective scheme over the field $k$ on which the reductive
group $G$ acts. Then, given a $K$-valued point $x$ of $X^{\rm ss}(L)$, where
$K$ is the quotient field of the complete discrete valuation ring $R$,
there exist a finite extension $R\subseteq R'$ and $g\in G(K')$,
$K'$ being the fraction field of $R'$, such that $g\cdot x$ is an
$R'$-valued point of $X^{\rm ss}(L)$.
\end{Thm}
%%%
So, to prove the semistable reduction theorem it is sufficient to
show that the constructed moduli space is a GIT quotient of a
projective scheme.

We fix a stability parameter $\widetilde{\delta}\succ\widetilde{\delta}_0$ (see Theorem \ref{HonestAlwaysSemStab})
and use the notation of Section \ref{sub:Moduli1}.
By elimination theory, the points in the parameter space ${\frak Y}$ corresponding to singular principal $G$-bundles
form an open subset ${\frak H}$. By Theorem \ref{HonestAlwaysSemStab},
${\frak H}^{\rm (s)s}:={\frak Y}^{\widetilde{\delta}\rm\hbox{-}ss}\cap {\frak H}$ is the open subset corresponding to
(semi)stable singular principal $G$-bundles.
It suffices to show that the quotients
$$
{\rm M}_P^{\rm(s)s}(\rho):={\frak H}^{\rm (s)s}\catqot \GL(U)
$$
exist as open subschemes of ${\rm M}_P^{\widetilde{\delta}\rm\hbox{-}(s)s}(\rho)$.
This will follow immediately, if we show that ${\frak H}^{\rm ss}$ is a \bfseries $\GL(U)$-saturated \rm
open subset of ${\frak Y}^{\widetilde{\delta}\rm\hbox{-}ss}$. This means that for every point
$x\in {\frak H}^{\rm ss}$
the closure of the $\GL(U)$-orbit in ${\frak Y}^{\widetilde{\delta}\rm\hbox{-}ss}$ is entirely contained in
${\frak H}^{\rm ss}$.
Since the points with closed $\GL(U)$-orbit in ${\frak Y}^{\widetilde{\delta}\rm\hbox{-}ss}$ are mapped to the
points with closed $\SL(U)$-orbit in $\ol{\frak Y}^{\widetilde{\delta}\rm\hbox{-}ss}$ (see \cite{OST}, proof of
Proposition 1.3.2), it suffices to show that
$$
\ol{\frak H}^{\rm ss}:={\frak H}^{\rm ss}\catqot {\Bbb G}_m(k)
$$
is an $\SL(U)$-saturated open subset of $\ol{\frak Y}^{\widetilde{\delta}\rm\hbox{-}ss}$.
Let $y,y^\p\in \ol{\frak Y}^{\widetilde{\delta}\rm\hbox{-}ss}$, such that $y^\p$ lies in the closure
of the $\SL(U)$-orbit of $y$. Then, $\psi(y^\p)$ lies in the closure of the $\SL(U)$-orbit of $\psi(y)$.
We may assume that the orbit of $y'$ and hence of $\psi(y')$ is closed.
By the Hilbert--Mumford criterion, one knows that there exists a one parameter subgroup
$\la\colon {\Bbb G}_m(k)\lra \SL(U)$ with $\lim_{x\rightarrow \infty}\la(z)\cdot \psi(y)\in \SL(U)\cdot \psi(y^\p)$.
Note that the injectivity of $\psi$ thus implies $\lim_{z\rightarrow \infty}\la(z)\cdot y\in \SL(U)\cdot y^\p$.
Suppose that $y$ and $\psi(y)$ represent $(\A,\tau)$ and $(\A,\phi)$, respectively.
Now, from the GIT constructions in Section \ref{subsec:ConstModDispSh}, in particular Remark \ref{rem:SEqAndOrbEq}
and Section \ref{sub:Moduli1}, one infers that $\la$
corresponds to a filtration $(\A_\bullet,{\alpha}_\bullet)$ with
$$
M(\A_\bullet,{\alpha}_\bullet)+\frac{\widetilde{\delta}}{s!}\cdot \mu(\A_\bullet,{\alpha}_\bullet;\phi)
\equiv 0
$$ and that a point in the orbit of $\psi(y^\p)$ represents ${\rm
df}_{(\A_\bullet,{\alpha}_\bullet)}(\A,\phi)$, so that a point in the
orbit of $y^\p$ represents ${\rm
df}_{(\A_\bullet,{\alpha}_\bullet)}(\A,\tau)$, by Section
\ref{subsec:SEqPseudo}.  Together with Remark \ref{rem:SEqSingPrinz},
this shows
$$
y\in \ol{\frak H}^{\rm ss}\q\Longrightarrow\q y^\p\in \ol{\frak H}^{\rm ss},
$$
and this is what we wanted to prove.\qed
%%%
\section{Appendix: Semistable reduction for good quotients}
\setcounter{subsection}{1}
%%%
In this appendix, we provide a short proof of Seshadri's theorem \ref{Seshadri's thm} (following his ideas) used in the proof of the semistable reduction theorem for singular principal $G$-bundles. As is well known to the experts (e.g., \cite{BS}), Seshadri's theorem together with the GIT construction of the moduli spaces gives the respective semistable reduction theorem. As an illustration, we show how we can recover the semistable reduction theorem of Langton for semistable sheaves and the semistable reduction theorem for curves. Even if one has constructed the moduli space as a projective scheme, the semistable reduction theorem remains of interest, because it has implications on the moduli stack or related stacks (see, e.g., \cite{HeinS} and \cite{Ngo}).
\par
Let $X$ be a scheme over some scheme $S$ and let $G$ be a
smooth affine $S$-group scheme acting on $X$.
As usual, for an $S$-scheme $T$, we set $X_T:=X\times_S T$.
In the whole section, $K$ denotes the quotient field of a discrete
valuation ring $R$.
Let us recall that an $S$-morphism $\pi\colon X\lra Y$ is called a
\emph{good quotient}, if $\pi$ is an affine $G$-invariant morphism, such that $\pi_\star({\cal O}_X)^G\simeq
{\cal O}_Y$.
%%%
\begin{Lem} 
Assume that there exists a good quotient $\pi\colon X\lra Y$, that $Y$ is
proper over $S$, and that there is a commutative diagram
$$
\xymatrix{
& \spec(K)\ar[r]^-{x}\ar[d] & X\ar[d]
\\
& \spec(R) \ar[r] & S.
}
$$
Then, there exist a finite extension $R\subseteq R'$ and $g\in G(K')$, $K'$ being
the fraction field of $R'$, such that $g\cdot x$ is an $R'$-valued point of $X$.
\end{Lem}
%%%%
\begin{proof}
Let $Z_K$ be the closure in $X_K$ of the $G_K$-orbit of the graph of $x\colon \spec(K)\lra X$.
Then, there exists a uniquely determined closed subscheme $Z_R$ of $X_R$, such that 
$Z_R\lra \spec(R)$ is flat and the generic fiber is isomorphic to $Z_K$.
It is the closure of the $G_R$-orbit of the graph of $\spec(K)\lra X$ in $X_R$. 
Let us remark that $Z_R\lra \spec(R)$ is faithfully flat, i.e., that the fiber $T$ over the closed
point of $\spec(R)$ is non-empty. This follows from the fact that $X_R\lra Y_R$
is a closed surjective map and the $K$-point $\pi(x)$ of $ Y$ can be extended to an $R$-point of $Y$. 
Now, the lemma follows from the existence of quasi-sections of faithfully flat morphisms.             
\end{proof}
%%%
The above lemma is a slight strengthening of a reformulation of \cite{BB-Sw}, Lemma 2.9.
It implies a generalization of Theorem \ref{Seshadri's thm} by the following remark:
By Seshadri's generalization of Mumford's GIT (see \cite[Theorem
4]{SeshGeomRed}), the assumptions of the lemma are satisfied, if $S$ is
of finite type over a universally Japanese ring, $G/S$ is a reductive
group scheme, acting on a projective scheme with a linearization in an ample line bundle on it, and 
$X$ is the open subset of $G$-semistable points.
%%%
\begin{Thm}[Stable reduction for curves]
The Deligne--Mumford stack of stable curves is proper over $\mathbb
Z$. More precisely, if $X\lra \spec(K)$ is a stable curve, then there
exist a finite extension $K\subset K'$  and a (unique) stable
family $X'\lra \spec(R')$, where $R'$ is the normalization of $R$ in
$K'$, such that the restriction of $X'$ to $\spec(K')$ is
isomorphic to $X\times _KK'$.
\end{Thm}
%%%
For the history and references concerning this theorem, we refer to
\cite{Edidin}.
%%%
\begin{proof}
The moduli scheme of stable curves is constructed as a GIT
quotient of the scheme $\overline{H}_g$ that parameterizes stable
curves of genus $g$ together with their $n$-canonical embeddings into
some ${\mathbb P}^N$ by an action of $\PGl (N+1)$ (see \cite{GiesCurve}). Since the GIT
quotient is projective, we can use the above lemma. A curve
$X\lra \spec(K)$, after choosing an embedding into ${\mathbb P}^N$,
gives rise to a map $\spec(K)\lra {\overline H}_g$. Then, 
after possibly changing the map with a group action, we can extend it to a
map $\spec(R') \lra {\overline H} _g$. This gives the required family, because there is a
universal family over ${\overline H}_g$, 
\end{proof}
%%%
\begin{Thm}[Langton's theorem; see \cite{HL}, Theorem 2.B.1]
Let $X$ be a projective $\mathbb Z$-scheme with geometrically
connected fibers and let ${\cal O}_X(1)$ be an ample line bundle on $X$.
Let $\F_K$ be a Gieseker semistable sheaf on $X\times \spec(K)$.
Then, there exist a finite extension $K\subset K'$  and a family
$\F'_{R'}$ of Gieseker semistable sheaves on $X$ parameterized by
$\spec(R')$, where $R'$ is the normalization of $R$ in $K'$, such
that $\F'_{K'}\simeq \F_K\otimes _K K'$.
\end{Thm}
%%%
In fact, Langton proved the slightly stronger assertion that in the
above theorem one can always take $K'=K$, but we need to start
with an $R$-flat family of sheaves on $X$. Langton's algorithm works also for slope semistable sheaves for which there is no moduli space in general.
%%%
\begin{proof}
The theorem follows from the above lemma and the GIT construction of the moduli space of Gieseker semistable sheaves (see \cite{La3}).
\end{proof}
%%%%
\section*{Acknowledgments}
%%%
We thank Professor Mehta for suggesting this research while
pointing out, in a conversation with the last author, that, since
our previous work in characteristic 0 deals with the semistability
of a tensor rather than with the semistability of a vector bundle,
it was more suited for attacking the case of positive
characteristic. Our thanks go also to N.\ Fakhruddin, J.\
Heinloth, Y.\ Holla, S.\ Koenig, V.\ Mehta, J.S.\ Milne, A.\
Premet, and S.\ Subramanian for discussions related to this work.
\par
Tom\'as G\'omez and Ignacio Sols are supported by grant number
MTM2004-07090-C03-02 from the Spanish ``Ministerio de Educaci\'on
y Ciencia". T.G.\ thanks the Tata
Institute of Fundamental Research for the hospitality during his
visit in December 2004. He also thanks the Polish Academy of
Science (Warsaw) and the University of Duisburg--Essen for
hospitality during his visits where parts of this work were
discussed.
\par
Ignacio Sols is supported by grant number
MTM2004-07090-C03-02 from the Spanish ``Ministerio de Educaci\'on
y Ciencia".
\par
Adrian Langer was partially supported by a Polish KBN grant
(contract number 1P03A03027). A.L.\ also thanks the DFG
Schwerpunkt ``Globale Methoden in der Komplexen Geometrie" for
supporting his visit to the University of Duisburg--Essen.
\par
Alexander Schmitt acknowledges support by the DFG via a Heisenberg fellowship, via the
Schwerpunkt program ``Globale Methoden in der Komplexen Geometrie---Global Methods in Complex Geometry'' and via SFB/TR 45 ``Periods, moduli spaces and arithmetic of algebraic varieties''.
A.S.\ was also supported by the DAAD via the ``Acciones Integradas Hispano--Alemanas"
program, contract number D/04/42257. In the framework of this program,
the third author visited the Consejo Superior de Investigaciones Cient\'\i ficas (CSIC) in Madrid
where some details of the paper were discussed. During the preparation of the revised version, A.S.\ stayed at the
Institut des Hautes \'Etudes Scientifiques where he benefited from support of the European Commission
through its 6th Framework Program ``Structuring the European Research Area" and contract No.\ RITA-CT-2004-505493
for the provision of Transnational Access implemented as Specific Support Action.
%%%

%%%
\subsection*{Addresses}
%%%
T.G.:
\\ 
Instituto de Ciencias Matemáticas (CSIC-UAM-UCM-UC3M)
\\ 
C/Serrano 113bis
\\ 
E-28006 Madrid
\\
Spain
\\ 
and
 \\
Departamento de Algebra
\\
Facultad de Ciencias Matem\'aticas
\\
Pza.\ de las Ciencias, 3
\\
Universidad Complutense de Madrid
\\
E-28040 Madrid
\\
Spain
\\
E-mail: \tt tomas.gomez@mat.csic.es\rm.
\\
\\
A.L.:  
Institute of Mathematics
\\ 
Warsaw University
\\ 
Ul.\ Banacha 2
\\ 
PL-02-097 Warszawa
\\
Poland
\\
and
\\
Institute of Mathematics
\\
Polish Academy of Sciences
\\
 Ul.\ \'Sniadeckich 8
\\
PO\ Box 21
\\
PL-00-956 Warszawa
\\
Poland
\\
E-mail: \tt alan@mimuw.edu.pl\rm.
\\
\\
A.H.W.S.: 
\\
Freie Universit\"at Berlin
\\ 
Institut f\"ur Mathematik
\\ 
Arnimallee 3
\\ 
D-14195 Berlin
\\
Germany
\\ 
E-mail: \tt alexander.schmitt@fu-berlin.de\rm.
\\
\\
I.S.: 
\\
Departamento de Algebra
\\
Facultad de Ciencias Matem\'aticas
\\
Pza.\ de las Ciencias, 3
\\
Universidad Complutense de Madrid
\\
E-28040 Madrid
\\
Spain
\\ 
E-mail: \tt isols@mat.ucm.es\rm.
\end{document}